\magnification=\magstep1
\input amstex
\UseAMSsymbols
\input pictex
\vsize=23truecm
\NoBlackBoxes
\parindent=18pt

   \font\rmk=cmr8    \font\itk=cmti8


\def\op{{\text{\rm op}}}

\def\image{\operatorname{Im}}

\def\mod{\operatorname{mod}}

\def\Hom{\operatorname{Hom}}
\def\sub{\operatorname{sub}}

\def\End{\operatorname{End}}
\def\Ext{\operatorname{Ext}}

\def\rad{\operatorname{rad}}
\def\add{\operatorname{add}}
\def\Ker{\operatorname{Ker}}

\def\soc{\operatorname{soc}}
\def\Tr{\operatorname{Tr}}

\def\bdim{\operatorname{\bold{dim}}}

\def\top{\operatorname{top}}

\def\ch{\operatorname{char}}
  \def\ss{\ssize }
\def\arr#1#2{\arrow <1.5mm> [0.25,0.75] from #1 to #2}

\def\s{\hfill \square}

\vglue0truecm
\centerline{\bf  Gorenstein-projective modules over short local algebras.}
                     \bigskip
\centerline{Claus Michael Ringel, Pu Zhang}
                \bigskip\medskip

\noindent {\narrower Abstract: \rmk Following the well-established
terminology in commutative algebra, any (not necessarily commutative)
finite-dimensional local algebra $\ss A$ with radical $\ss J$
will be said to be {\itk short} provided $\ss J^3 = 0$. As in the commutative case, we show:
If a short local algebra $\ss A$ has an indecomposable  non-projective Gorenstein-projective module $\ss M$,
then either $\ss A$ is self-injective (so that all modules are Gorenstein-projective)
and then, of course, $\ss |J^2| \le 1$,
or else $\ss |J^2| = |J/J^2| - 1$ and $\ss |JM| = |J^2||M/JM|.$
More generally, we focus the attention to
semi-Gorenstein-projective and $\ss\infty$-torsionfree modules, even to $\ss\mho$-paths
of length 2, 3 and 4. In particular, we show that the existence of a non-projective
reflexive module implies that $\ss |J^2| < |J/J^2|$ and further restrictions.
In addition, we consider exact complexes of projective modules with a non-projective image.
Again, as in the
commutative case, we see that if such a complex exists, then $\ss A$ is self-injective
or satisfies the condition $\ss |J^2| = |J/J^2| - 1.$
Also, we show that any non-projective semi-Gorenstein-projective
module $\ss M$ satisfies
$\ss \Ext^1(M,M) \neq 0$. In this way, we prove the
Auslander-Reiten conjecture (one of the classical homological
conjectures) for arbitrary short local algebras.

Many arguments used in the commutative case actually
work in general, but there are interesting differences and
some of our results may be new also in the commutative case.
	\medskip
\noindent
2010 Math Subject classification. Primary 16G10, Secondary 13D07, 16E65, 16G50, 20G42.
	\smallskip
\noindent
Supported by NSFC 12131015, \ 11971304.
\par}
	\bigskip
{\bf 1. Introduction.}
	\medskip
{\bf 1.1. The algebras and their modules} 
	\smallskip
Let $A$ be a finite-dimensional algebra with radical $J = J(A).$
The modules to be considered are left $A$-modules of finite length (if not
otherwise asserted). We denote by $|M|$ the length of the module $M$.
If $M$ is a module, let $p\colon PM \to M$ be a projective cover of $M$ and
$\Omega M$ the kernel of $p$. The modules $\Omega^i M$ with $i\ge 1$ are
the syzygy modules of $M$. The module $\top M = M/JM$ will be called
the top of $M$ and we write $t(M) = |\top M|.$

All algebras $A$ considered here will be local finite-dimensional $k$-algebras,
where $k$ is a field,  and
for simplicity, we will assume that $A/J = k$. The module $A/J$ will always be denoted
by $S$; it is the unique simple module. 
Let $e = e(A) = |J/J^2|$. 
A local algebra $A$ is said to be {\it short} provided $J^3 = 0.$
Usually, we will assume that $A$ is short and then we write
$a = a(A) = |J^2|$ and call $(e(A),a(A))$ the {\it Hilbert-type} of $A$.

If $M$ is a module with Loewy length at most 2, we call
$\bdim M = (t(M),|JM|)$ (or its transpose) the {\it dimension vector} of $M$
(note that $\bdim M$ is only defined for modules $M$ of Loewy length at most 2;
we have $\bdim S = (1,0)$ and there is no module
with dimension vector $(0,1)$). We call a module $M$ {\it bipartite} provided $\soc M = J M$.
{\it A module has Loewy length at most $2$ if and only if it is the direct sum of
a bipartite and a semisimple module.}
	\bigskip
{\bf 1.2. Complexes and $A$-duality} 
	\smallskip
For any module $M$, let
$M^* = \Hom(M,{}_AA)$ be the {\it $A$-dual} of $M$
(it is a right $A$-module, thus an $A^{\op}$-module), and
$\phi_M\colon  M \to M^{**}$ the canonical map defined by $\phi(m)(\alpha) = \alpha(m)$
for $m\in M$ and $\alpha\in M^*$.
A module $M$ is  {\it torsionless} if  $M$ is a submodule
of a projective module, or, equivalently, if $\phi_M$
is a monomorphism.
The module $M$ is said to be {\it reflexive} if $\phi_M$
is bijective. Note that {\it an indecomposable module which is torsionless and not
projective has Loewy length at most $2.$}
	
We will consider exact complexes of projective modules, they are of the form
$P_\bullet = (P_i,d_i\colon P_i \to P_{i-1})_i$, thus
$$
 \cdots  @>>> P_{1} @>d_{1}>> P_{0} @>d_0>> P_{-1} @>d_{-1}>> P_{-2} @>>> \cdots,
$$
with projective modules $P_i$ such that $\image d_i = \Ker d_{i-1}$, for all $i\in \Bbb Z$.
A module $M$ is said to be an {\it image in}\enspace
$P_\bullet$, provided $M = \image d_i$ for some $i \in \Bbb Z$.
The exact complex $P_\bullet$ is said to be
{\it minimal} provided that any map
$d_i$ maps into the radical of $P_{i-1}$.
Given a complex $P_\bullet$ of projective modules,
we may form the $A$-dual complex $P_\bullet^*$, forming the $A$-dual of the modules
$P_i$ and of the maps $d_i.$
	
A module $M$ is  {\it Gorenstein-projective} provided it is an image in
an exact complex $P_\bullet$ of projective modules, with $P_\bullet^*$
again being exact; $M$ is
{\it semi-Gorenstein-projective} provided $\Ext^i(M,A) = 0$
for all $i\ge 1,$ and $M$ is {\it $\infty$-torsionfree,} provided $\Tr M$ is
semi-Gorenstein-projective,
where $\Tr$ is Auslander's transpose operator.
Note that {\it a module $M$ is Gorenstein projective iff $M$ is both
semi-Gorenstein-projective and $\infty$-torsionfree.}
	\bigskip
{\bf 1.3. The topics to be considered} 
	\smallskip
The topics to be discussed in this paper (and its sequel [RZ3])
concern the module theory for a
short local algebra $A.$  The main results of the
present paper are stated in Sections 1.4 to 1.8.
The central question concerns the existence of non-projective
Gorenstein-projective modules or of related ones, and properties of such modules.
The case of $A$ being commutative has been
considered before in several papers published between 1980 and 2010
(in particular, see [L, Y, HSV, CV, AIS]).
Our aim is to extend the results known for commutative rings
to general rings. Some of our observations may be
new also in the commutative case.

Two properties of Gorenstein-projective modules are important:
Gorenstein-projective modules are reflexive,
and they are images in exact complexes of projective modules.
Theorems 1.1 and 1.2 announced in Section 1.4 deal with the existence
of non-projective reflexive modules. Theorem 1.3 stated in Section 1.5
concerns the existence of
non-projective images in exact complexes of projective modules.
	\bigskip
{\bf 1.4. Existence of reflexive modules}
	\smallskip
We say that a non-zero module $M$ of Loewy length at most 2 is {\it solid} provided
any endomorphism of $M$ is a scalar multiplication on $\soc M$
(as a consequence, any non-invertible
endomorphism vanishes on the socle). A solid module
is of course indecomposable (a characterization of the solid modules
using covering theory will be given in Proposition A.3 of Appendix A).
	\medskip
{\bf Theorem 1.1.} {\it Let $A$ be a short local algebra which is not
self-injective. If
there exists a reflexive module which is not projective, then
$2 \le a \le e-1$. Also, ${}_AJ$ and the right module $J_A$ are solid.}
	\medskip
The bound $a \le e-1$ cannot be improved, 
since Proposition 15.1 shows that for any $a$ with 
$1\le a \le e-1$, there exists an algebra of
Hilbert type $(e,a)$ with non-projective reflexive modules.

Also note that in general, ${}_AJ$ may be solid, whereas
$J_A$ is not solid, as Example 4.8 shows.
	
Theorem 1.1 can be rephrased:
If $A$ is any artin algebra, there exists a non-projective reflexive module
iff there exists a non-projective module $N$ with $\Ext^i(N,A) = 0$
for $i=1,2.$  Namely, there is the following recipe: If $N$ is
a non-projective module
with $\Ext^i(N,A) = 0$ for $i = 1,2$, then the module $\Omega^2N$ is
non-projective and reflexive. The reverse construction is given by the agemo-functor
$\mho = \Tr\Omega\Tr$: If $M$ is reflexive, then
$\Ext^i(\mho^2 M,A) = 0$, for $i= 1,2.$
This recipe is part of general considerations outlined in Section 2,
which focus the attention to what we call the $\mho$-paths of $A$.
With reference to $\mho$-paths of $A$, the existence assumption in Theorem 1.1
just says that there exists an $\mho$-path of length 2.
	
Theorem 1.1 assumes that there exists a non-projective reflexive module, thus
an $\mho$-path of length 2,
or equivalently, a non-projective module $M$ with $\Ext^i(M,A) = 0,$ for $1\le i \le 2.$
The next theorem shows that
the existence of a non-projective module $M$ with
$\Ext^i(M,A) = 0,$ for $1\le i \le 4$, yields a stronger assertion. Again
using Section 2.4, there are several  reformulations.
The existence of a
non-projective module $M$ with $\Ext^i(M,A) = 0,$ for $1\le i \le 4$ is
equivalent to the existence of an $\mho$-path of length 4, and also to
the existence of a non-projective reflexive
module $M$ with $\Ext^i(M,A) = 0$ for $i = 1,2$.

	\medskip
{\bf Theorem 1.2.} {\it Let $A$ be a short local algebra which is not
self-injective. Assume that $M$ is an indecomposable, reflexive and
non-projective module with $\Ext^i(M,A) = 0$ for $1=1,2.$
Then $2 \le a = e-1$. If $t = t(M)$, then
$\bdim X = (t,at)$ for $X \in \{\Omega^2M,\Omega M, M, \mho M\}.$}
	\bigskip
{\bf 1.5. Existence of exact complexes of projective modules}
	\medskip
{\bf Theorem 1.3.} {\it Let $A$ be a short local algebra which is not self-injective, with
a non-zero minimal exact complex $P_\bullet = (P_i,d_i)_i$
of projective modules. Then  $1\le a = e-1.$
	\smallskip
Also, $M_i = \image d_i$ is bipartite for $i\ll 0$. Let $t_i = t(P_i) = t(M_i).$
There is $v\in \Bbb Z$ such that for $i\le v$, we have $t_i = t$ and
$\bdim M_i = (t,at)$. And, there are just two possibilities:
	\smallskip
{\rm Type I.} For all $i\in \Bbb Z,$ the module
 $M_i$ is bipartite with $\bdim M_i = (t,at)$ (thus $t_i = t$).
	\smallskip
{\rm Type II.} We can choose $v$ in such a way that first, $t_{i+1} > t_i$ for $i\ge v$, second, the module $M_{v+1}$ is not bipartite, and third, $|JM_i| < at_i$ for $i > v$.}
	\medskip
For commutative rings, Theorem 1.3 is
due to Christensen-Veliche [CV]; here,
the case $a = 1$ does not occur. But in general, the case $a = 1$ is possible, see 
Example 9.2.
Also, for $A$ commutative, and $P_\bullet$ a complex of type II,
all the modules $M_i$ with $i\le v$ are bipartite, whereas we do not
know whether this holds true in general.
For $A$ commutative, the existence of a non-zero minimal exact
complex $P_\bullet$ of projective modules
implies that $J^2 = \soc A$, whereas in general, it
does neither imply that $J^2 = \soc {}_AA$, nor that
$J^2 = \soc A_A$, see Examples 9.2 and 9.3. If $J^2 = \soc {}_AA$,
then all the modules $M_i$ with $i\le v$ are bipartite, see Corollary 13.3.

For a typical example of a complex $P_\bullet = (P_i,d_i)$
of type I, see Proposition 10.7:
Let $A$ be  of Hilbert type $(e,e-1)$, with $e\ge 2$ and $x\in A$
with $x^2 = 0$ and $Jx = J^2$ (a left Conca element). For all $i$, let $P_i = {}_AA$ and $d_i$
the right multiplication by $x$. Then all images in $P_\bullet$
are equal to $Ax$. If $x$ is also right Conca, then also $P_\bullet^*$ is exact, thus
$Ax$ is Gorenstein-projective.

Theorem 1.3 describes the structure of a minimal exact complex of projective modules,
if $A$ is not self-injective. For $A$
being self-injective, see Corollary A.8 in Appendix A.
	\bigskip
{\bf 1.6. Semi-Gorenstein-projective and $\infty$-torsionfree modules}
	\smallskip
Both Theorems 1.2 and 1.3 imply:
{\it If $A$ is a short local algebra which is not
self-injective, with a Gorenstein-projective module which is not
projective, then $2 \le a = e-1$.} There is the following strengthening.
		\medskip
{\bf Theorem 1.4.} {\it Let $A$ be a short local algebra which is not
self-injective.  Assume that there exists a non-projective indecomposable module $M$ which is
semi-Gorenstein-projective or $\infty$-torsionfree. Then
$2 \le a = e-1$ and $J^2 = \soc {}_AA = \soc A_A.$  
Moreover, let $t = t(M).$ We have in addition:

\item{\rm (1)} 
If $M$ is torsionless and semi-Gorenstein-projective, then
$\bdim \Omega^iM  = (t,at)$ for all $i\ge 0.$

\item{\rm (2)}  
If $M$ is $\infty$-torsionfree, then
$\bdim \mho^iM  = (t,at)$ for all $i\ge 0.$

\item{\rm (3)} 
If $M$ is reflexive and semi-Gorenstein-projective, or if $M$ is
$\infty$-torsionfree, then also $\bdim M^* = (t,at).$

\item{\rm (4)}
If $M$ is Gorenstein-projective, then
$\bdim X = (t,at)$ for $X =  \Omega^iM$ and $X = \mho^iM$, where $i\ge 0$, as well
as for $X = M^*.$
\par}
	\medskip
{\bf Remark.} In general, if $A$ is a short local algebra and 
$M$ is semi-Gorenstein-projective, its Loewy length
may be 3; and if it is 2, we may have $\bdim M^* \neq \bdim M$
(see the algebra $A$ mentioned in Example 9.5: the right $A$-module $M = m_1A$ is 
semi-Gorenstein-projective and has $\bdim M = (1,2),$
whereas $\bdim M^* = (2,1)$, the right $A$-module $\mho(m_1A)$ is
also semi-Gorenstein-projective and its Loewy length is 3).
	\bigskip
{\bf 1.7. The Auslander-Reiten conjecture} 
	\smallskip
Using Theorem 1.4 as well as Proposition A,5 in Appendix A
we get the following result.
	\medskip
{\bf Theorem 1.5.} {\it Let $A$ be a short local algebra and
$M$ a non-projective semi-Gorenstein-projective module. Then
$\Ext^1(M,M) \neq 0$. Moreover,
if $A$ is not self-injective, then $\Ext^i(M,M)\neq 0$
for all $i\ge 1.$}
	\medskip
Recall that the Auslander-Reiten conjecture [AR] for
an artin algebra $A$ asserts: {\it If $M$ is a non-projective
semi-Gorenstein-projective module, then $\Ext^i(M,M)\neq 0$
for some $i\ge 1$.} Thus, Theorem 1.5 shows that
the Auslander-Reiten conjecture holds true for short local algebras
in a stronger from.
For $A$ self-injective, Theorem 1.5 is due to Hoshino [Ho1], 1982.
For commutative short local rings, 
a proof of the Auslander-Reiten conjecture was given
by Huneke-\c Sega-Vraciu in 2004. 
	
Let $A$ be a short local algebra which is self-injective. Let $M$ 
be a non-projective module.
Then Theorem 1.5 asserts that $\Ext^1(M,M) \neq 0$ (over a self-injective algebra,
all modules are semi-Gorenstein projective). If $A$ is, in addition, commutative, then
we even have $\Ext^i(M,M) \neq 0$ for  all $i\ge 1$, see Huneke-\c Sega-Vracio [HSV].
However, for $A$ non-commutative, this is not true: we
may have $\Ext^i(M,M) = 0$ for some or even for all $i \ge 2,$ see 
Proposition A.19 in Appendix A.
	\bigskip
{\bf 1.8. Existence of $\mho$-paths of length 3} 
	\smallskip
We have mentioned that
for any pair $(e,a)$ with $1 \le a \le e-1$,
there are short local algebras of Hilbert type $(e,a)$ with
a non-projective reflexive module (see Proposition 15.1), thus with an
$\mho$-path of length 2, whereas 
the existence of an $\mho$-path of length 4 implies that $a = e-1$ (see Theorem 1.1).
Proposition 15.2 provides an example of a short local algebra of Hilbert type $(6,2)$
with a non-projective $3$-torsionfree module, thus with an $\mho$-path of length 3.

We do not know whether for any pair $(e,a)$ with $2 \le a \le e-2$
there is a short local algebra which has non-projective $3$-torsionfree modules,
thus $\mho$-paths of length 3.
	\bigskip  
{\bf 1.9. Summary} 
	\smallskip
The short local algebras $A$ with $e\le 1$ are self-injective
Nakayama algebras, thus let us restrict to the short local algebras $A$ with $e\ge 2$.
They can be separated as follows:
	\smallskip
\item{$(1)$} $a = 0$ (thus $A$ is a radical-square-zero algebra).
   Reflexive modules and images in exact complexes of projective modules are projective.
   (Theorems 1.1 and 1.3).
\item{$(2)$} $a = 1$ (this includes the self-injective algebras).
   There may be $\mho$-paths of
   arbitrary length. There may be non-projective images in exact complexes of projective modules.
\item{$(3)$}  $2 \le a \le e-2$. There may be $\mho$-paths of length 3, but never of length 4
   (Theorem 1.2).
   Images in exact complexes of projective modules are projective (Theorem 1.3).
\item{$(4)$} $2 \le a = e -1$. There may be $\mho$-paths of
   arbitrary length,
   and there may be non-projective images in exact complexes of projective modules.

\item{$(5)$} $e \le a.$
   Reflexive modules and images in exact complexes of projective modules are projective
   (Theorems 1.1 and 1.3).
	\smallskip
The paper [RZ3] shows a further separation, namely
between $a\le \frac14\,e^2$ and $\frac14\,e^2 < a.$
	\smallskip
It has turned out that several short local algebras with $a = e -1$
are of great interest, see Gasharov-Peeva [GP, 1990],
Avramov-Gasharov-Peeva [AGP, 1997], Veliche [V, 2002], Yoshino [Y, 2002],
Jorgensen-\c Sega [JS2, 2006], Christensen-Veliche [CV, 2007],
Hughes-Jorgensen-\c Sega [HJS, 2009], all dealing with commutative rings.
A non-commutative  algebra $A$ of Hilbert type $(3,2)$
has been analyzed in [RZ1,RZ2]; the construction will be
generalized in Section 11. 
	\smallskip
For $e\ge 3$, Section 11 exhibits a short local
algebra $\Lambda$ with $a = e-1$ which has a non-projective Gorenstein-projective
module $M$, a semi-Gorenstein-projective module $M'$ which is not torsionless,
and an $\infty$-torsionfree module $M''$ with $\Ext^1(M'',A) \neq 0$.
In particular, $\Lambda$ has complexes of type I and II.
	\bigskip
\vfill\eject
{\bf 1.10. Outline of the paper}
	\smallskip
The proofs of Theorems 1.1 and 1.2 are given in Section 4 and 6, respectively.
The proofs of Theorems 1.3 and 1.4 can be found
in Section 9. The proof of Theorem 1.5 is given in Section 12.

Many arguments used in the commutative case work in general, but there are
also some decisive differences. For the convenience of the reader, we will
provide complete proofs, the only exceptions
are the use of the appendix of [CV], see Lemma 9.1 below, and of
basic properties of the $\mho$-quiver and the $\mho$-paths, see Section 2
(here we follow [RZ1]).

Throughout the paper, $L(e)$
denotes the local $k$-algebra with $J^2 = 0$, $|J| = e$ and $L(e)/J = k.$
If $A$ is any local algebra with $e(A) = e$ (and $A/J = k$), then
$A/J^2 = L(e)$ and we will interpret the
$L(e)$-modules as the $A$-modules annihilated by $J^2$, thus as the
$A$-modules of Loewy length at most 2.

Often, we will assume that $A$ is not self-injective. After all, over a
self-injective algebra, all modules are Gorenstein-projective.
Appendix A provides an overview over the module theory of
self-injective (equivalently, Gorenstein)
short local algebras and the local radical-square-zero algebras $L(e)$,
based on the relationship between these algebras and the Kronecker quivers
$K(e)$. 

The essence of Sections 6 and 7 is:
If one is interested in exact complexes of projective modules, or in long
$\mho$-paths,
then the cases $e\le a$ and $2\le a\le e-2$ can be discarded, and
one has to look at the case $a = e-1$.
This case is considered in Sections 7, 10, 11, 12 and in
the examples 9.3 and 9.4.
In particular, we show in Corollary 10.5 that a commutative short local algebra
of Hilbert type $(e, e-1)$ has no complex of type II which involves
a projective module of rank 1.

Examples of algebras with or without non-projective modules which are reflexive or
are images in exact complexes of projective modules
are constructed in Sections 14 and 15. In particular, we show that for any
pair $(e,a)$ of integers with $2\le a \le e-1$, there is an algebra of Hilbert
type $(e,a)$ with a non-projective reflexive module. Also, we provide
an example of an algebra of Hilbert type $(6,2)$ with a non-projective
$3$-torsionfree module.

Sections 3 and 8 are devoted to the simple module $S$, its syzygies
and the $\mho$-component which contains $S$. We stress
that for any local algebra, if $S$ is reflexive or the image in an exact complex
of projective modules, then $A$ is self-injective, see Lemma 3.2.

The main tool in the paper
will be the use of the transformation $\omega^e_a$ on $\Bbb Z^2$ as defined
in Section 5:
It describes for suitable modules $M$ in which way $\bdim M$
is changed when we apply $\Omega_A$ (see the Main Lemma 5.4 and 13.1, 
but also [RZ3]). The Main Lemma draws the attention
to the possible equality $t(\Omega^2 M) = et(\Omega M)-at(M)$.
 Appendix B
is devoted to the numbers $b_n = b(e,a)_n$ defined recursively by the corresponding
rule $b_{n+1} = eb_n-ab_{n-1}$, starting with $b_{-1} = 0,\ b_0 = 1$. It presents an
explicit formula for these numbers $b_n$  due to Avramov, Iyengar, \c Sega, provided
$a < \frac 14e^2$.
	\smallskip
We hope that the use of two independent numbering systems as suggested by the journal
does not lead to confusion: the sections and subsections are numbered
consecutively; independently, the assertions, examples and some of the remarks
are also numbered consecutively. 
	\bigskip\medskip
\vfill\eject
{\bf 2. The $\mho$-quiver and the $\mho$-paths}
	\medskip
In this section, $A$ will be an arbitrary artin algebra.
We provide a survey on the $\mho$-quiver (and the $\mho$-paths), following
[RZ1, Sections 1.5, 4.4 (and also 1.9)].

The $\mho$-quiver was introduced in [RZ1] in order to formalize ideas
which are due to Auslander
(1968), Auslander-Bridger (1969) and Auslander-Reiten (1996)
(and which were further elaborated by many others) building up
 the realm of Gorenstein-projective modules (for historical references, in particular for
the bibliographical data of relevant papers, see [RZ1]). The $\mho$-quiver
provides the general frame for several important module theoretical
concepts which carry deviating names: torsionless modules, reflexive modules,
(semi-)Gorenstein-projective modules, $n$-torsionfree modules
(a module $M$ is said to be {\it $n$-torsionfree,} provided $\Ext^i(\Tr M,A) = 0$
for $1\le i \le n$), $\infty$-torsionfree modules, and so on; in particular,
it explains the wording ``totally reflexive" used by
Avramov-Martsinkovsky (2002) for the Gorenstein-projective
modules. Finally, it highlights the
duality between semi-Gorenstein-projective modules
and $\infty$-torsionfree modules.
	\bigskip
{\bf 2.1. The operator $\mho$}
	\smallskip
Let $M$ be a module.
We denote by $\mho M$ 
the cokernel of a minimal left $\add A$-approximation
of $M$ (equivalently, we may define $\mho M = \Tr\Omega\Tr M$,
where $\Tr$ is Auslander's transpose operator, see [RZ1], Lemma 4.4);
the operator $\mho$ is called the {\it agemo} operator. 
	\medskip
{\it Let $M$ be a module.}
\item{$\bullet$}
{\it The module $\mho M$ has no indecomposable projective
direct summands.}
\item{$\bullet$}
{\it If $M$ is indecomposable, not projective and torsionless, then
$\mho M$ is indecomposable (and not projective).}
\item{$\bullet$} {\it
The module $M$ is reflexive iff both $M$ and $\mho M$ are torsionless.}
	\bigskip
{\bf 2.2. The $\mho$-quiver}
	\smallskip
The vertices of the {\it $\mho$-quiver} are the isomorphism classes $[X]$
of the indecomposable non-projective modules $X$ and there is an arrow
$$
{\beginpicture
    \setcoordinatesystem units <2cm,1cm>
\put{$[X]$} at 0 0
\put{$[Z]$} at 1 0
\setdashes <1mm>
\plot 0.7 0  0.3 0 /
\setsolid
\arr{0.31 0}{0.3 0}
\endpicture}
$$
provided $X = \Omega Z$ and $\Ext^1(Z,A) = 0$, or, equivalently,
provided $X$ is torsionless and $Z = \mho X$.
If $X$ is torsionless (and indecomposable and non-projective), then there is a
(uniquely determined) exact
sequence $0 \to X \to P \to \mho X \to 0$ with $P$ projective (thus $X\to P$
is a minimal left $\add A$-approximation); such a sequence
is called an {\it $\mho$-sequence.} In this way, the arrows of the $\mho$-quiver just
correspond to the $\mho$-sequences.
This explains the direction of the arrow $[X] \leftarrow [\mho X]$ used here: The
usual convention for using arrows in order to draw attention to
short exact sequences $0 \to X \to Y \to Z \to 0$ is to draw an 
arrow $[X] \leftarrow [Z]$ (and often one uses a dashed arrow).

Paths in the $\mho$-quiver will be
called {\it $\mho$-paths,} the connected components of the $\mho$-quiver will be
called {\it $\mho$-components.}

A decisive feature of the $\mho$-quiver is the following: Any module $M$ is the start of
{\bf at most
one} arrow in the $\mho$-quiver (and then this arrow ends in $\Omega M$)
and also the end of {\bf at most
one} arrow in the $\mho$-quiver (and then this arrow starts in $\mho M$).
Thus any $\mho$-component is a linearly oriented quiver
$\Bbb A_n$ with $n\ge 1$ vertices, or an oriented cycle $\widetilde{\Bbb A}_n$ with
$n\!+\!1\ge 1$ vertices, or of the form $-\Bbb N,$ or $\Bbb N,$  or $\Bbb Z$.
(Note that we consider any subset $I$ of $\Bbb Z$
as a quiver, with an arrow from $z$ to $z\!-\!1$
provided that both $z\!-\!1$ and $z$ belong to $I$.)
	\bigskip

{\bf 2.3. Dictionary} 
	\smallskip
{\it Let $M$ be an indecomposable non-projective module.
	\medskip
\item{$\bullet$}
$M$ is torsionless iff $M$ is the end of an $\mho$-path of length $1$.
\item{$\bullet$}
$M$ is reflexive iff $M$ is the end of an $\mho$-path of length $2$.
\item{$\bullet$}
$M$ is $n$-torsionfree iff $M$ is the end of an $\mho$-path of length $n$.
\item{$\bullet$}
$M$ is $\infty$-torsionfree iff $M$ is the end of an infinite $\mho$-path.
	\medskip
\item{$\bullet$}
$\Ext^i(M,A) = 0$ for $1\le i \le t$ iff $M$ is the start of an $\mho$-path of length t.
\item{$\bullet$}
$M$ is semi-Gorenstein-projective iff $M$ is the start of an infinite $\mho$-path.
	\medskip
\item{$\bullet$}
$M$ is Gorenstein-projective iff $M$ is the start of an infinite $\mho$-path
and the end of an infinite $\mho$-path (thus iff the $\mho$-component containing $M$
is an oriented cycle $\widetilde{\Bbb A}_n$, or of the form $\Bbb Z$).\par}
	\bigskip
{\bf 2.4. Some bijections} 
	\smallskip
The operators
$\Omega^2$ and $\mho^2$ provide
inverse bijections between isomorphism classes as follows:
$$
{\beginpicture
\setcoordinatesystem units <1.9cm,1cm>
\put{$\left\{ \matrix \text{\rm indecomposable} \cr
                          \text{\rm non-projective modules $M$}\cr
                   \text{\rm which are reflexive}
                   \endmatrix\right\}$} at -1 0
\put{$\left\{ \matrix \text{\rm indecomposable} \cr
                          \text{\rm non-projective modules $M$}\cr
                 \text{\rm with  $\Ext^i(M,A)=0$ for $i=1,2$}
                   \endmatrix\right\}$} at 3.25 0
\arr{0.5 0.1}{1.5 0.1}
\arr{1.5 -.1}{0.5 -.1}
\put{$\Omega^2$} at 1 -.35
\put{$\mho^2$} at 1 .35
\endpicture}
$$
this is the bijection between the end and the start of the $\mho$-paths of length 2.

In the same way, we may look at the $\mho$-paths of length 4. We obtain the
following bijections (again, all modules $M$ are assumed to be indecomposable and non-projective) looking at
the end, the middle  and the start, respectively, of any $\mho$-path of length 4.
$$
{\beginpicture
\setcoordinatesystem units <2.5cm,1cm>
\put{$\biggl\{ \matrix M \cr
 \text{\rm $4$-torsionfree}
                   \endmatrix \biggr\}$} at 0 0
\put{$\Biggl\{ \matrix \text{\rm $M$ reflexive,} \cr
                          \text{\rm $\Ext^i(M,A)=0$}\cr
                 \text{\rm for $i=1,2$}
                   \endmatrix\Biggr\}$} at 2 0
\put{$\biggl\{ \matrix \text{\rm $\Ext^i(M,A)=0$ } \cr
                 \text{\rm for $i=1,2,3,4$}
                   \endmatrix\biggr\}$} at 4 0
\arr{0.7 0.1}{1.3 0.1}
\arr{1.3 -.1}{0.7 -.1}
\put{$\Omega^2$} at 1 -.35
\put{$\mho^2$} at 1 .35
\arr{2.7 0.1}{3.3 0.1}
\arr{3.3 -.1}{2.7 -.1}
\put{$\Omega^2$} at 3 -.35
\put{$\mho^2$} at 3 .35
\endpicture}
$$
	\medskip
{\bf 2.5. $A$-duality} 
	\smallskip

{\it Let $0 \to X \to P \to Z \to 0$ be an $\mho$-sequence.
If $Z$ is reflexive, then also $0 \to Z^* \to P^* \to X^*\to 0$ is an
$\mho$-sequence and $X$ (thus also $X^*$) is reflexive.}
	\medskip
Proof. For the first
assertion, see [RZ1], 4.2(b). If $Z$ is reflexive, then $\mho X = Z$
and $\mho^2 X = \mho Z$
both are torsionless, thus $X$ is reflexive. $\s$
	\medskip
In terms of $\mho$-paths, the assumption that $Z$ is reflexive
means that there is an $\mho$-path (in $\mod A$)
of length 3 as
shown below on the left, the conclusion that $X^*$ is reflexive
concerns the existence of the $\mho$-path of length 3 in $\mod A^{\op}$ shown on the right.
$$
{\beginpicture
\setcoordinatesystem units <1.2cm,1cm>
\put{\beginpicture
\put{$X$} at 0 0
\put{$Z$} at 1 0
\multiput{$\circ$} at 2 0  3 0 /
\setdashes <1mm>
\arr{0.7 0}{0.3 0}
\arr{1.7 0}{1.3 0}
\arr{2.7 0}{2.3 0}
\endpicture} at 0 0
\put{\beginpicture
\put{$Z^*$} at 0 0
\put{$X^*$} at 1 0
\multiput{$\circ$} at 2 0  3 0 /
\setdashes <1mm>
\arr{0.7 0}{0.3 0}
\arr{1.7 0}{1.25 0}
\arr{2.7 0}{2.3 0}
\endpicture} at 5 0
\endpicture}
$$
	\bigskip
\vfill\eject
{\bf 3. The $\mho$-component of $S$}
	\medskip
We collect some general observations concerning finite-dimensional local algebras $A$
which are not necessarily short, mostly well-known. 
	\medskip
{\bf Lemma 3.1.} {\it Let $A$ be a local artinian ring. Any module of
finite projective dimension is projective.}
	\medskip
Proof. Let $m$ be the Loewy length of $A$.
Assume that $M$ is a module with finite projective dimension $t\ge 1.$ Let
$$
   0 \to P_t \to \cdots \to P_0 \to M \to 0
$$
be a minimal projective resolution, thus $P_t \neq 0$.
Now $P_t$ is a submodule of $\rad P_{t-1}.$ But
$P_t$ has Loewy length  $m$, whereas $\rad P_{t-1}$ 
has Loewy length $m-1$, impossible. 
$\s$
	\medskip
Several characterizations of finite dimensional local algebras which are
self-injective:
	\medskip
{\bf Lemma 3.2.} {\it Let $A$ be a finite-dimensional local algebra. 
The following assertions are equivalent:
\item{\rm(i)} $\soc{}_AA$ is simple.
\item{\rm(ii)} $A$ is self-injective.
\item{\rm(iii)} All modules are Gorenstein-projective.
\item{\rm(iv)} All modules are reflexive.
\item{\rm(v)} $S$ is reflexive.
\item{\rm(vi)} $\mho S$ has Loewy length at most $m-1$, where $m$ is the Loewy length of $A$.
\item{\rm(vii)} All modules are images in exact complexes of projective modules.
\item{\rm(viii)}
 $S$ is the image in an exact complex of projective modules.
\item{\rm(ix)} $S$ is the kernel of a map $g\colon P \to P'$ with $P,\ P'$ projective.
\item{\rm(x)} {\it $I({}_AA)$ has finite projective dimension.}
\item{\rm(xi)} $\soc{}A_A$ is simple.\par}
	\medskip
The left-right-symmetry of the assertions (i) and (xi)
means that {\it we
may also add the right module versions of the assertions} (ii) {\it to} (x).
	\medskip
Proof. (i) $\implies$ (ii): If
 $\soc{}_AA$ is simple, then the injective envelop of ${}_AA$
is indecomposable. But the indecomposable injective $A$-module has the same dimension as
$A$, thus ${}_AA$ is injective. (ii) $\implies$ (iii) is well-known. Any
Gorenstein-projective module is reflexive and is an image in an exact complex
of projective modules. Thus we have (iii) $\implies$ (iv) $\implies$ (v),
as well as (iii) $\implies$ (vii) $\implies$ (viii).
Of course, there are the obvious implications
(viii) $\implies$ (ix), then (v) $\implies$ (ix), and also (ii) $\implies$ (x). 
	\medskip
(v) $\implies$ (vi). If $M$ is indecomposable, reflexive and not projective,
then $\mho M$ is indecomposable, torsionless and not projective (see Section 2.1), thus
there is an embedding $\mho M \subseteq JP$ with $P$ projective. Therefore,
the Loewy length of $\mho M$ is at most $m-1.$
	\medskip
(vi) $\implies$ (i). We
assume that $\mho S$ has Loewy length at most $m-1$.
Let $a = |J^{m-1}|.$  By assumption, $a\ge 1.$
Since $S$ is torsionless, there is an $\mho$-sequence $0 \to S @>u>> P @>p>> \mho S \to 0.$
Let $P$ be of rank $t$. Thus $t\ge 1$ and
$|J^{m-1}P| = at$.
Since $\mho S$ has Loewy length at most $m-1$,
$J^{m-1}P$ is contained in the kernel of $p$, thus $at \le 1,$ and therefore
$a = 1$ and $t = 1.$

Assume now that there is a simple submodule $U$ of ${}_AA$ which is not contained in
$J^{m-1}$.
Let $v\colon U \to A$ be the inclusion map. Let $f\colon S \to U$ be an isomorphism. Since
$u$ is a left $\add(A)$-approximation, there is $f'\colon P \to A$ with $f'u = vf.$

Let us assume that $f'$ is not surjective. Then the image of $f'$ is a module
of Loewy length at most ${m-1}$, thus $J^{m-1}P$
is contained in the kernel of $f'$. We have $J^{m-1} \neq 0$. Since $J^{m-1}P
\subseteq \Ker(p) = \image (u)$ and $\image(u)$ is simple, we see that $J^{m-1}P = \image(u).$
It follows that $f'u = 0$ in contrast to $vf \neq 0.$

Thus we see that $f'$ is
surjective. There is $f''\colon  \mho S \to A/U$ such that the following diagram commutes:
$$
\CD
 0 @>>> S @>u>>     P @>p>> \mho S @>>> 0 \cr
 @.   @VVf V           @VVf'V  @VVf''V \cr
 0 @>>> U @>v>> {}_AA @>>> A/U @>>> 0.
\endCD
$$
Since $f'$ is surjective, also $f''$ is surjective.
Since $J^{m-1}$ is not contained in $U$,
the module $A/U$ has Loewy length $m$. Therefore also $\mho S$ has Loewy length $m$,
a contradiction. This shows that $\soc{}_AA \subseteq J^{m-1}$. Since $a = 1$, it follows that
$\soc{}_AA$ is simple.
	\medskip
(ix) $\implies$ (xi). Let $S$ be the kernel of
a map $g\colon P \to P'$ with $P,P'$ projective. Write both
$P$ and $P'$ as direct sums
of copies of ${}_AA$, thus $g$ is given by a matrix with entries $g_{ij} \in
\End({}_AA) = A$ and we can assume that all entries belong to $J$.
But this implies that $\bigoplus \soc A_A$ is contained in the kernel of $g$.
Since the kernel of $g$ is simple, we see that ($P = {}_AA$ and that)
$\soc A_A$ is simple.
	\medskip
(xi) $\implies$ (ii). In the previous parts of the proof, we have seen that $ii)$
implies (xi). If we apply this to the opposite algebra of $A$, we see that (xi) implies
(ii). 
	\medskip
(x) $\implies$ (ii). If $I({}_AA)$ has finite projective dimension, then Lemma 3.2
asserts that $I({}_AA)$ is projective. 
$\s$
	\medskip
{\bf Remark 3.3.} We recall that an algebra $A$ is said to be {\it Gorenstein} provided
both modules ${}_AA$ and $A_A$ have finite injective dimension, or, equivalently, 
provided both modules $I({}_AA)$ and $I(A_A)$ have finite projective dimension. 
The equivalence
of (ii) and (x) shows that {\it a finite dimensional local algebra is Gorenstein iff
it is self-injective} (in commutative algebra, it is customary to refer to these algebras
as Gorenstein algebras). 
	\smallskip
{\bf Remark 3.4.} Both (vi) and (ix) imply that $S$ is the kernel of a map
$g\colon P \to Z$ with $P$ projective and $Z$ of Loewy length at most $m-1$.
However, $S$ may be the kernel of a map
$g\colon P \to Z$ with $P$ projective and $Z$ of Loewy length at most $m-1$,
whereas $A$ is not self-injective: Take the algebra
$A = k\langle x,y\rangle/\langle x^2,y^2,xy\rangle$ and $Z = A/Ayx.$
	\smallskip
{\bf Remark 3.5.}  The implication (v) $\implies$ (ii) has been shown by Marczinzik in [M1],
and he used this opportunity to ask whether any finite-dimensional algebra is self-injective provided
all simple modules are reflexive. This is not true, see [R2].
	\smallskip
{\bf Remark 3.6.} According to Theorems 1.1 and 1.3, the existence of a non-projective
reflexive module or a non-projective image in an exact complex of projective
modules, implies severe restrictions on the algebra $A$, however there do exist
many algebras which are not self-injective with such modules. As we see
in Lemma 3.2, the situation is different, if $S$ itself is reflexive, or is
the image in an exact complex of projective: This can happen only if $A$
is self-injective.
	\bigskip
{\bf Lemma 3.7.} {\it Let $A$ be a local algebra.
The following conditions are equivalent.
\item{\rm(i)} $\Ext^1(S,{}_AA) = 0.$
\item{\rm(ii)} $A$ is self-injective.\par}
	\medskip
Proof. Of course, (ii) implies (i). Conversely, assume that $\Ext^1(S,{}_AA) = 0.$
Then $\Ext^1(M,{}_AA) = 0$ for all $A$-modules $M$, thus ${}_AA$ is injective.
$\s$
	\medskip
{\bf Corollary 3.8.} {\it Let $A$ be a local algebra which is
not self-injective. Then the $\mho$-component
of $A$ which contains $S$ is of type $\Bbb A_2$ with $[S]$ as its sink.}
	\medskip
Proof. Since $S$ is torsionless, there is an arrow ending in $S$. Since $S$ is not
reflexive, there is no $\mho$-path
of length 2 ending in $S$. Since $\Ext^1(S,A) \neq 0$, no
arrow starts in $S$. $\s$
	\bigskip
{\bf Proposition 3.9.}
{\it A short local algebra is self-injective if and only if
either $a = 0$ and $e\le 1$ or else $a = 1$ and $J^2 = \soc {}_AA$.}
	\medskip
Proof. According to Lemma 3.2, $A$ is self-injective if and only if $\soc{}_AA$ is simple.
If $a = 0$ and $e\le 1$ or if $a = 1$ and $J^2 = \soc {}_AA$,
then $\soc{}_AA$ is simple, thus $A$ is self-injective. Conversely,
assume that $A$ is self-injective.
If $J^2 = 0$, and $J \neq 0$,
then  $\soc {}_AA = J$, thus $a = 0,\ e = 1$; if $J^2 \neq 0,$  then
$J^2 \subseteq \soc{}_AA$, thus we must have $a = 1$ and $J^2 = \soc {}_AA$.
 $\s$
	\bigskip\bigskip
{\bf 4. Reflexive modules and the proof of Theorem 1.1}
	\medskip
{\bf 4.1. Reflexive modules}
	\smallskip
We assume that $A$ is a short local algebra. We want to analyze the structure of
reflexive modules. As we will see, the existence of a reflexive module which is not
projective puts severe restrictions on $A$.
	\medskip
{\bf Lemma 4.1.} {\it Let $A$ be a short local algebra.
Let $M$ be indecomposable, torsionless and not projective. Then $M$ is bipartite or simple.}
	\medskip
Proof. Since $M$ is torsionless, there is an embedding $u\colon M \to P$ with $P$ projective.
Let $P = \bigoplus_i P_i$ with $P_i = {}_AA$ for all $i$. The composition of $u$
with any projection $P \to P_i$
cannot be surjective, since otherwise it would split and $M$ would have a direct summand
isomorphic to ${}_AA.$ Thus the image of $u$ is contained in $JP$ and therefore
of Loewy length at most $2$. It follows that $M$ is the direct sum of a bipartite
module and a semi-simple module. Since $M$ is indecomposable, it is bipartite
or simple. $\s$
	\medskip
{\bf Lemma 4.2.} {\it Let $A$ be a short local algebra which is not self-injective.
Let $M$ be a module which is indecomposable, reflexive and not projective. Then both modules
$M$ and $\mho M$ are bipartite.}
	\medskip
Proof. According to Lemma 4.1, 
$M$ is bipartite or simple. According to Lemma 3.2, $M$ cannot be
simple, thus $M$ is bipartite.
Since $M$ is indecomposable, reflexive, and not projective, $\mho M$ is
indecomposable, torsionless, and not projective. Using again 4.1, we see that
$\mho M$ has Loewy length at most 2. Since $\Ext^1(\mho M,A) = 0$ (see Section 2.2),
Lemma 3.2 asserts that $\mho M$ cannot be simple. Thus also $\mho M$ is bipartite.
$\s$
	\medskip
{\bf Lemma 4.3.} {\it Let  $A$ be a short local algebra which is not self-injective.
Then, the module ${}_AB = {}_AA/J^2$ is not reflexive.}
	\medskip
Proof. Any map ${}_AB \to {}_AA$ maps into $J$, thus
$\Hom({}_AB,{}_AA) = \Hom({}_AB,{}_AJ),$ and we can identify
$\Hom({}_AB,{}_AJ)$ with $J_A$, sending $\phi\colon {}_AB \to {}_AJ$ to
$\phi(1).$ Thus, ${}_AB^* = J_A$. It is sufficient to show that $\dim\ (J_A)^* =
\dim \Hom(J_A,A_A) > 1+e$, since $\dim B = 1+e.$
But $\Hom(J_A,A_A)$ has the proper subspace $\Hom(\top J_A,\soc A_A)$,
and this subspace has dimension at least $2e$, since $\soc A_A$ is not
simple.
$\s$
	\medskip
{\bf 4.2. Proof of Theorem 1.1}
	\medskip
Let us repeat the assertion. 
	\medskip
{\it Let  $A$ be a short local algebra which is not self-injective.
Let $M$ be a module which is reflexive and not projective. 
Assume that 
there exists a reflexive module which is not projective, then
$2 \le a \le e-1$.
Also, the module ${}_AJ$ and the right $A$-module $J_A$ are solid.}
	\medskip
Proof. By assumption, there is a reflexive module $M$ which is not projective.
In addition, we can assume that $M$ is indecomposable. 
Let $\bdim M = (t,s)$ and $z = |\top \mho M|.$
Let $0 \to M @>u>> P @>p>> \mho M \to 0$ be an $\mho$-sequence, where $P$ is
projective of rank $z$, thus also $z = |\top \mho M|.$
Of course, we can assume that $u$ is an inclusion map.

In the following, we denote by $X^{(z)}$
the direct sum of $z$ copies of a module $X$.
We have $P = {}_AA^{(z)},$ with
$JP = {}_AJ^{(z)}$ and $M$ is a submodule of $JP.$
	\medskip
(1) {\it We have $s < et.$}
	\medskip
Proof. Let $B = A/J^2$, thus $B = L(e)$. Since $A$ is not self-injective, we have $e\ge 2.$
Since $M$ is bipartite, it is a $B$-module. Its projective cover as a $B$-module
is of the form $p'\colon P' \to M$ with $\bdim P' = (t,et)$. Since $p'$ is surjective,
we have $s \le et$.

Now assume that $s = et$. Then $p'$ is an isomorphism, thus $M$
is a projective $B$-module.
Since $M$ is indecomposable, $M$ is
the projective left $B$-module of rank $1$. However, according to Lemma 4.3,
the module ${}_AB$ is not reflexive.
	\medskip
(2) {\it We have $\soc M = J^2P$ and therefore $s = az.$}
	\medskip
Proof. Since $\mho M$ has Loewy length at most 2, we have
$J^2P \subseteq \Ker(p) = M.$ Since $J^2P$ is semisimple, it follows that $J^2P \subseteq
\soc M.$ On the other hand, $M \subseteq JP$ implies that
$\soc M = JM \subseteq J^2P$, thus $\soc M = J^2P.$

By definition, $a = |J^2|.$ Altogether, $s = |JM| = |\soc M| = |J^2P| = |J^2|z = az.$
$\s$
	\medskip
(3) {\it We have $J^2 = \soc {}_AA$ and therefore $a \ge 2.$}
	\medskip
Proof.
If $J^2 \neq \soc{}_AA$, there is a simple submodule $U$ of ${}_AA$ which is not
contained in $J^2$. Let $f\colon M \to U$ be a homomorphism with image $f(M) = U,$
and $v\colon U \to {}_AA$ the inclusion map.
Since $u\colon M \to P$ is a left $\add(A)$-approximation,
there is $f'\colon P \to {}_AA$ such that $vf = f'u.$
If $f'$ is not surjective, then $f'(P) \subseteq J$, thus $f'(JP) \subseteq J^2$
and therefore $f'u(M) \subseteq J^2$. But $f'u = vf$ and $vf(M) = v(U) = U$
is not contained in $J^2.$ This shows that $f'$ is surjective. There is the
following commutative diagram
$$
\CD
 0 @>>> M @>u>>     P @>p>> \mho M @>>> 0 \cr
 @.   @VVf V           @VVf'V  @VVf''V \cr
 0 @>>> U @>v>> {}_AA @>>> A/U @>>> 0.
\endCD
$$
Since $f'$ is surjective, also $f''$ is surjective. Since $U$ is not contained in $J^2$,
the module $A/U$ has Loewy length 3. Thus, also $\mho M$ has Loewy length 3.
But $\mho M$ has Loewy length at most 2, see Lemma 4.2.
This contradiction shows that $J^2 = \soc {}_AA$.

By definition, $a = |J^2|,$ thus $a = |\soc{}_AA|.$
Since $A$ is not self-injective, we have $|\soc{}_AA| \ge 2,$ see Lemma 3.2.
$\s$
	\medskip
(4) {\it ${}_AJ$ is solid.}
	\medskip
Proof. Let $\phi$ be an endomorphism of ${}_AJ.$

Recall that $P = {}_AA^{(z)}$ with inclusion map $u\colon M \to P,$ thus we can write
$u$ as the transpose of $[u_1,\dots,u_z]$, where $u_i\colon M \to {}_AA = A_i.$
As we know, $\soc M = J^2P = \bigoplus_{i=1}^z J^2A_i,$ thus $J^2A_i = A_i\cap \soc M.$

We denote the inclusion map $JA_1 \subset A_1$ by $v_1$ and write $u_1 = v_1u'_1$, where
$u'_1\colon M \to JA_1.$ Let $f\colon M \to {}_AA$ be the composition
$$
 M @>u'_1>> JA_1 @>\phi>> JA_1 @>v_1>> A_1 = {}_AA.
$$
Since $u$ is an $\add({}_AA)$-approximation, there are maps
$g_i\colon {}_AA\to {}_AA$ such that $g = [g_1,\dots,g_z]$ satisfies
$f = gu = \sum g_iu_i$. The map $g_1\colon {}_AA\to {}_AA$ is the right multiplication
by some element $\lambda\in A$.

Given $x\in A_1\cap\ \soc M = J^2A_1$, we consider the element $[x,0,\dots,0] \in M$ and
apply the map $f = \sum g_iu_i$
to it.
Since $f = v_1\phi u'_1$, we have $f([x,0\dots,0]) = \phi(x).$
On the other hand, we have $u_i([x,0,\dots,0]) = 0$ for $i\ge 2$, thus
$\sum g_iu_i([x,0,\dots,0]) = g_1(x) = x\lambda.$
This shows that
$$
 \phi(x) = f([x,0\dots,0]) = \sum g_iu_i([x,0,\dots,0])= x\lambda
$$ for all $x\in J^2A_1$.
Now $J^2A_1$ is annihilated from the right by $J$, thus $x\lambda = \overline \lambda x$,
where $\overline \lambda = \lambda+J$ is an element of $A/J = k$.
This shows that the restriction
of $\phi$ to $J^2A_1 = J^2$ is the scalar multiplication by $\overline \lambda.$
By (3), $J^2 = \soc {}_AA$. It follows that ${}_AJ$ is solid.
$\s$
	\medskip
(5) {\it We have $ez \ge at.$}
	\medskip
Proof. We use again the decomposition $P = {}_AA^{(z)}.$ We have
$JP = {}_AJ^{(z)}$ and   $M$ is a submodule of $JP.$
Let $u'\colon M \to J^{(z)}$, $v\colon J \to A$ and $w\colon J^2\to A$ be
the canonical inclusion maps. Thus $u = v^{(z)}u'.$ Given $a\in A$, we denote by
$r(a)\colon {}_AA \to {}_AA$ the right multiplication by $c$.
If $c\in J$, then $r(c)$ maps $J$
into $J^2$ and the map $r(c)\colon J \to J^2$ depends only on the residue class $\overline c$
of $c$ modulo $J^2.$ Thus we may write $r(\overline c) = r(c)\colon  J \to J^2$ and there is the following commutative diagram
$$
\CD J @>v>> A \cr
  @Vr(\overline c)VV     @VVr(c)V \cr
    J^2 @>w>> A
\endCD
$$
 In this way, we obtain the following linear map
$$
 \Phi\colon (J/J^2)^{(z)} \to \Hom(M,J^2), \quad \text{defined by}\quad
\Phi(\overline c_1,\dots,\overline c_z) = [r(\overline c_1),\dots,
r(\overline c_z)]u'.
$$

Let us show that $\Phi$ is surjective.
Let $f\colon M \to J^2$ be any homomorphism.
By assumption, the inclusion map $u= v^{(z)}u'\colon M \to A^{(z)}$ is a
left $\add({}_AA)$-approximation.
Thus, there is $f'\colon A^{(z)} \to {}_AA$ such that $wf = f'u.$
We write $f'$ as $[r(c_1),\dots,r(c_z)]$ with elements $c_i\in A$.
Since $f$ vanishes on $\soc M = (J^2)^{(z)}$, we have $(J^2)c_i = 0$, thus $c_i\in J$,
for all $1\le i \le z$.

Thus, we have the following diagram.
$$
{\beginpicture
    \setcoordinatesystem units <2cm,2cm>
\put{$M$} at 0 1
\put{$J^{(z)}$} at 1 1
\put{${}_AA^{(z)}$} at 2 1
\put{$J^2$} at 0 0
\put{${}_AA$} at 2 0
\arr{0 0.8}{0 0.2}
\arr{2 0.8}{2 0.2}
\arr{0.2 1}{0.8 1}
\arr{1.2 1}{1.8 1}
\arr{0.2 0}{1.8 0}
\arr{0.85 0.85}{0.15 0.15}
\put{$f$} at -.1 0.5
\put{$f'= [r(c_1),\dots,r(c_z)]$} [l] at 2.05 0.5
\put{$\ss [r(\overline c_1),\dots,r(\overline c_z)]$} [l] at 0.6 0.45
\put{$u'\strut$} at 0.5 1.1
\put{$v^{(z)}$\strut} at 1.5 1.1
\put{$w$} at 1 .1
\endpicture}
$$
Here, the outer rectangle commutes by the choice of $f'$.
Since $c_i\in J$, we have $r(c_i)v = wr(\overline c_i)$, thus
$[r(c_1),\dots,r(c_z)]v^{(z)} = w[r(\overline c_1),\dots,r(\overline c_z)]$. Since $w$
is a monomorphism, it follows that also the triangle on the left commutes:
$f = [r(\overline c_1),\dots,r(\overline c_z)]u'$. Thus, we see that
$$
 f = [r(\overline c_1),\dots,r(\overline c_z)]u'
 =
 \Phi(\overline c_1,\dots,\overline c_z).
$$
In this way, we see that $\Phi$ is surjective, thus
$\dim\ (J/J^2)^{(z)} \ge \dim\Hom(M,J^2)$

Now, $\dim\ (J/J^2)^{(z)} = ez$.
Second, any map $M \to J^2$ factors through the projection
$M \to \top M$, thus  $\dim\Hom(M,J^2) = \dim \Hom(\top M,J^2) = ta$.
Therefore $ez \ge ta$.
$\s$
	\medskip
(6) {\it We have $a < e$.}
	\medskip
Proof. Assume for the contrary, that $e \le a$. Using (2) and (1), we have
$az = s < et \le at,$  and therefore $z < t$.
Using (5), we have $at \le ez \le az,$ thus $t\le z.$ Thus, we obtain a
contradiction.
$\s$
	\medskip
(7) {\it The right $A$-module $J_A$ is solid.}
	\medskip
Proof. If $M$ is a reflexive and non-projective module, then $M^*$ is a
reflexive and non-projective $A^\op$-module. Thus (3) asserts that $J_A$ is solid.
$\s$
	\medskip 
The assertions (4), (6) and (3) and (7) are as required.
This completes the proof.
$\s$

	\bigskip
{\bf Corollary 4.4.} {\it Let $A$ be a short local algebra which is not
self-injective. If there exists a reflexive module which is not projective,
then both modules ${}_AJ$ and $J_A$ are solid and of Loewy length $2$. In particular,
$\soc {}_AA = J^2 = \soc A_A.$}
	\medskip
Proof. Theorem 1.1 asserts that ${}_AJ$ is solid, and that $a\ge 2$, thus
${}_AJ$ is bipartite and of Loewy length 2.
It follows that $\soc {}_AA = \soc{}_AJ = J^2.$
If $M$ is reflexive and not-projective, then $M^*$
is a reflexive and non-projective $A^{\op}$-module, thus $A^{\op}$ satisfies also
the assumptions of Theorem 1.1.
$\s$
	\bigskip
{\bf Remark 4.5.}
Note that an element $c\in J$ belongs to $\soc A_A = \soc J_A$ if and only if $cJ = 0$.
As a consequence, {\it $J^2 = \soc A_A$ if and only if ${}_AJ$ is a faithful $A/J^2$-module.}
	\bigskip
{\bf 4.3. Some Examples}
	\smallskip 

{\bf Example 4.6.} {\it A short local algebra with $J^2 = \soc {}_AA \subset \soc A_A$.}
Since our general assumption is $J^3 = 0$,
we always have $J^2 \subseteq \soc{}_AA$
as well as $J^2 \subseteq \soc A_A.$ We may have $J^2 = \soc {}_AA$ and
$J^2 \neq \soc A_A$ as the following example shows.
Let $A$ be the $k$-algebra with radical generators $x,y$ and relations
$$
{\beginpicture
    \setcoordinatesystem units <1cm,1cm>
\put{$yx,\ y^2,\ x^3,\ x^2y.$} at -5 0.5
\put{$x$} at 0 1
\put{$y$} at 1 1
\put{$x^2$\strut} at 0 0
\put{$xy$\strut} at 1 0
\arr{0 0.8}{0 0.25}
\arr{1 0.8}{1 0.25}
\multiput{$\ss x$} at -.2 .6  0.8 .6  /
\put{${}_AJ$} at -1 1
\endpicture}
$$
Here, $J^2 =
Ax^2 + Axy = \soc {}_AA$ is of length 2, whereas $\soc A_A = x^2A + yA + xyA$
is of length 3.
	\bigskip
{\bf Examples 4.7.} {\it Short local algebras
with ${}_AJ$ indecomposable, but not
solid.} First example: Here, ${}_AJ$ has a non-zero nilpotent endomorphism.

Let $A$ be generated by $x,y,z$ with relations
$$
  z^2,\ xy,\ yx,\ yz,\ zy,\ zx-xz,\ y^2-xz,\ x^3.
$$
$$
{\beginpicture
    \setcoordinatesystem units <.8cm,1cm>
\put{$x$} at 0 1
\put{$y$} at 2 1
\put{$z$} at 4 1
\put{$x^2$} at 0 0
\put{$y^2$} at 2 0
\arr{0 0.8}{0 0.2}
\arr{0.2 0.8}{1.8 0.2}
\arr{2 0.8}{2 0.2}
\arr{3.8 0.8}{2.2 0.2}
\multiput{$\ss x$\strut} at -.2 0.6  3 .7 /
\multiput{$\ss y$\strut} at 2.2 0.6 /
\multiput{$\ss z$\strut} at 1 .7 /
\put{${}_AJ$} at -1.3 .8
\endpicture}
$$
There is the endomorphism $f$ of ${}_AJ$
given by $f(y) = f(z) = 0$ and $f(x) = z$.

Second example: Here we exhibit an $\Bbb R$-algebra such that $\End({}_AJ) \sim \Bbb C.$
We consider the $\Bbb R$-algebra with generators $x,y$, and the relations are
$$
{\beginpicture
    \setcoordinatesystem units <.8cm,1cm>
\put{$ xy-yx,\ x^2+y^2.$} at -6 0.7

\put{$x$} at 0 1
\put{$y$} at 2 1
\put{$x^2$} at 0 0
\put{$y^2$} at 2 0
\arr{0 0.8}{0 0.2}
\arr{0.2 0.8}{1.8 0.2}
\arr{2 0.8}{2 0.2}
\setdashes <1mm>
\arr{1.8 0.8}{0.2 0.2}
\multiput{$\ss x$\strut} at -.2 0.6  2.2 0.6 /
\multiput{$\ss y$\strut} at 0.7 0.9  1.3 0.9 /
\put{${}_AJ$} at -1 .8
\endpicture}
$$
(Note that the 2-Kronecker module $\widetilde J$ as mentioned in Section A.2 of Appendix A
is $(\Bbb C,\Bbb C;1,i)$, where we write $1$ for the
identity map $\Bbb C \to \Bbb C$ and $i\colon \Bbb C \to \Bbb C$ for the multiplication by $i$;
of course, $\End(\Bbb C,\Bbb C;1,i) = \Bbb C.$)

Note that both algebras are commutative.
	\medskip
{\bf Example 4.8.} {\it A short local algebra with ${}_AJ$ solid, whereas
$J_A$ is not solid.}
Let $A$ be generated by $x,y,z$ with relations
$$
  x^2,\ y^2,\ z^2,\ yx,\ yz,\ zx-xy,\ zy-xz.
$$
$$
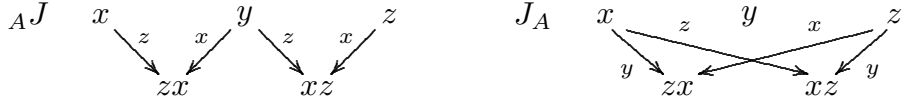
{\beginpicture
    \setcoordinatesystem units <.8cm,.8cm>
\put{\beginpicture
\put{$x$} at 0 1
\put{$y$} at 2 1
\put{$z$} at 4 1
\put{$zx$} at 1 0
\put{$xz$} at 3 0
\arr{0.2 0.8}{0.8 0.2}
\arr{1.8 0.8}{1.2 0.2}
\arr{2.2 0.8}{2.8 0.2}
\arr{3.8 0.8}{3.2 0.2}
\multiput{$\ss x$} at 1.4 0.7  3.4 0.7 /
\multiput{$\ss y$} at /
\multiput{$\ss z$} at 0.6 0.7  2.6 0.7 /
\put{${}_AJ$} at -1 1
\endpicture} at 0 0
\put{\beginpicture
\put{$x$} at 0 1
\put{$y$} at 2 1
\put{$z$} at 4 1
\put{$zx$} at 1 0
\put{$xz$} at 3 0
\arr{0.1 0.8}{0.8 0.2}
\arr{0.3 0.8}{2.8 0.2}
\arr{3.7 0.8}{1.3 0.2}
\arr{3.9 0.8}{3.2 0.2}
\multiput{$\ss x$\strut} at 2.9 0.9 /
\multiput{$\ss y$\strut} at 0.3 0.3  3.7 0.3 /
\multiput{$\ss z$\strut} at 1.1 0.9  /
\put{$J_A$} at -1 1
\endpicture} at 7 0
\endpicture}
$$
Here, ${}_AJ$ is solid, whereas $J_A$ is the direct sum of a module with dimension
vector $(2,2)$ and a simple module (generated by $y$). Note that ${}_AJ$ is solid,
but not faithful.
	\medskip
Note that Theorem 1.1 and its Corollary 4.4 assert that the algebras exhibited in
Examples 4.6, 4.7 and 4.8 do not
have non-projective reflexive modules, thus all semi-Gorenstein-projective and
all $\infty$-torsionfree modules are projective.
	\bigskip
{\bf 4.4. Short local algebras $A$ with $a \le 1$}
	\smallskip
Recall that a module of finite length is said to be {\it uniform} provided it has
a simple socle. If the
module $M$ has Loewy length at most 2, then $JM$ is simple if and only
if $M$ is the direct sum of a uniform module and a semisimple module.
Thus, if $A$ is a short local algebra, then $a \le 1$ and
$e\ge 1$ if and only if ${}_AJ$ is the direct sum of a uniform module and a semisimple module.
	\medskip
{\bf Lemma 4.9.} {\it Let $A$ be a short local algebra with $a \le 1.$
The following assertions are equivalent:
\item{\rm(i)} $A$ is self-injective and $J\neq 0.$
\item{\rm(ii)} There exists a non-projective reflexive module.
\item{\rm(iii)} ${}_AJ$ is solid.
\item{\rm(iv)} ${}_AJ$ is indecomposable.
\item{\rm(v)} ${}_AJ$ is uniform.
\item{\rm(vi)} ${}_AJ$ is simple or bipartite.
\item{\rm(vii)} Either $a = 0$ and $e = 1$, or else $a = 1$ and $J^2 = \soc{}_AA$.}
	\medskip
The proof is straightforward:
(i) $\implies$ (ii): If $J\neq 0$, then there are non-projective modules.
For $A$ self-injective, all modules are reflexive. (ii) $\implies$ (iii): See Theorem 1.1.
(iii) $\implies$ (iv): Solid modules are indecomposable. (iv) $\implies$ (v): An
indecomposable module $M$ with $|JM| \le 1$ is uniform.
(v) $\implies$ (vi): Clear. (vi) $\implies$ (vii):
If ${}_AJ$ is simple, then $a = 0,\ e = 1$. Otherwise $J^2$ is the socle of
${}_AJ$, and thus $a = 1$.
(vii) $\implies$ (i): See Lemma 3.2.
$\s$
	\bigskip
{\bf 4.5. Short local algebras $A$ with $e \le 2$}
	\medskip
{\bf Lemma 4.10.}
{\it Let $A$ be a short local algebra with $e \le 2.$
The following assertions are equivalent:
\item{\rm(i)} $A$ is self-injective and $J\neq 0.$
\item{\rm(ii)} There exists a non-projective reflexive module.
\item{\rm(iii)} ${}_AJ$ is uniform.
\item{\rm(iv)} Either $a = 0$ and $e = 1$, or else $a = 1$ and $J^2 = \soc{}_AJ$.}
	\medskip
Again, the proof is straightforward:
(i) $\implies$ (ii): If $J\neq 0$, then there are non-projective modules.
For $A$ self-injective, all modules are reflexive. (ii) $\implies$ (iii): Since there exists
a non-projective reflexive module, $e\ge 1.$
If $e = 1$, then $a = 0$ or $a = 1$ and in both cases ${}_AJ$ is of course uniform. Thus, according to Theorem 1.1, we can assume that $a < e = 2$  and that $M = {}_AJ$ is solid.
Since $M$ is indecomposable, it follows that $a \neq 0$. But
$|JM| = a = 1$ implies that $M = {}_AJ$ is uniform. (iii) $\implies$ (iv): Assume that
${}_AJ$ is uniform. Either ${}_AJ$ is simple, then $a = 0$ and $e = 1,$
or else $J^2 = \soc{}_AJ$ and $a = |J^2| = 1.$ (iv) $\implies$ (i): See Lemma 3.2.
$\s$
	\medskip
{\bf Example 4.11.} The algebra $A = k[x,y]/(x,y)^3$  is a short local algebra with $e = 2$
such that ${}_AJ$ is solid, thus indecomposable, but (of course) not uniform.
	\bigskip\bigskip
{\bf 5. Bipartite modules}
	\medskip
{\bf 5.1. $\mho$-sequences over short local algebras}
	\smallskip
First, let us apply the observations of Section 3 to $\mho$-sequences over short local algebras.
	\medskip
{\bf Corollary 5.1.} {\it Let $A$ be a short local algebra and
$0 \to X \to P \to Z \to 0$ an $\mho$-sequence.}
	\smallskip
{\rm (a)} {\it If $A$ is self-injective,
then either $X$ is bipartite, or else
$X$ is simple and then $Z = A/\soc {}_AA$.}
	\smallskip
{\rm (b)} {\it If $A$ is not self-injective,
and $Z$ has Loewy length at most $2$, then $Z$ is bipartite, and either $X$ is
also bipartite or else $X$ is simple and $a = 0,\ e\ge 2.$}
	\medskip
Proof.
(a)
The module $X$ is indecomposable and of Loewy length at most $2$. Thus, if $X$
is not simple, then $X$ is bipartite. If $X = S$ is simple, then $Z = A/\soc {}_AA.$
	
(b) Both $X$ and $Z$ are indecomposable modules of Loewy length at most 2.
Now $Z$ cannot be simple, since otherwise Lemma 3.2 asserts
that $A$ is self-injective. Since $X$ is indecomposable, it is either bipartite or
simple. If $X = S$ is simple, then Lemma 3.1 
shows that the Loewy length of $A$ cannot be 3
(since we assume that $Z = \mho S$ has Loewy length at most 2). Thus $a = 0.$
Since $A$ is not self-injective, we have $e\ge 2.$
$\s$
	\medskip
Let us add also the following observation.
	\medskip
{\bf Lemma 5.2.} {\it Let $A$ be a short local algebra.
If $M$ is a reflexive module which is bipartite, then also $M^*$
is (reflexive and) bipartite.}
	\medskip
Proof. We can assume that $M$ is indecomposable.
If $M$ is projective, then $M = {}_AA$ implies that $A$ has Loewy length 
2, thus also $M^* = A_A$ is bipartite. Thus, we assume that $M$ is not projective. Of
course, $M^*$ is torsionless. If $M^*$ would be projective, also $M$ would be projective.
Thus $M^*$ has Loewy length at most 2. Also $M^*$ cannot be simple, since otherwise
$A$ is self-injective and also $M$ is simple.
$\s$
	\medskip
Proposition 6.1 will provide more
information on the $A$-dual $M^*$ of a bipartite reflexive module $M$.
	\medskip
{\bf Example 5.3.}
{\it If $M$ is torsionless and bipartite, then $M^*$ has
Loewy length at most $2,$ but does not have to be bipartite.}

Namely, if $M$ is bipartite, then $M$ is annihilated by $J^2$,
thus any map $f\colon M \to {}_AA$
maps into $J$. If $x\in J^2$, then the 
right multiplication $r(x)\colon {}_AA\to {}_AA$ 
by $x$
sends $J$ to $0$, thus $r(x)f = 0$. Thus shows that $M^*$ has Loewy length at most $2$.

A typical example is given by the algebra $A = \Lambda_0$ considered in Section 11
(and before in [RZ1]), namely the {\bf right} $A$-module $m_1A = (x-y)A$,
as discussed in [RZ1].
Of course, $m_1A$ is torsionless and bipartite, but
$(m_1A)^* = M(q)^{**} = \Omega M(1)$ (see 6.5 (8) and Theorem 1.6 in [RZ1])
has a simple direct summand.
	\bigskip\medskip
\vfill\eject
{\bf 5.2. The Main Lemma}
	\smallskip
Given arbitrary integers $a$ and $e$, let
$$
 \omega^e_a = \bmatrix e & -1 \cr
                    a & 0 \endbmatrix.
$$
If $M$ is a module of Loewy length at most $2$, we write $\omega^e_a \bdim M$ for
multiplying $\omega^e_a$ with (the transpose of) $\bdim M$. 
	\medskip
{\bf Lemma 5.4.}
{\it If $M$ is a module of Loewy length at most $2$, then there is a natural number $w \ge 0$
such that
$$
 \bdim \Omega M = \omega^e_a \bdim M + (w,-w),
$$
and such that $\Omega M$ has a direct summand of the form $S^w$.
In particular, if $\Omega M$ is
bipartite, then}
$$
 \bdim \Omega M = \omega^e_a \bdim M.
$$
	\medskip
Proof. Let $M' = \Omega M$. There is an exact sequence $0 \to M' \to P \to M \to 0$ with
$P$ projective and we can assume that the map $M' \to P$ is an inclusion map. Let
$U = J^2P.$ Since $M$ has Loewy length at most 2, $U$ is mapped under $P \to M$ to zero,
thus $U \subseteq M'$. Since $U$ is semisimple, we have $U \subseteq \soc M'.$ Also,
$M'$ is a submodule of $JP$, thus
$M'/U$ is a submodule of $JP/J^2P$ and therefore semisimple. This shows that
$JM' \subseteq U.$ Let $w = |U/JM'|.$ Then
$$
 \bdim M' = (|M'/JM'|,|JM'|) =(|M'/U|+w,|U|-w) = (|M'/U|,|U|)+(w,-w).
$$
It remains to calculate $|U|$ and $|M'/U|.$ Let $\bdim M = (t,s)$. Then
$P = {}_AA^t$, thus $|U| = |J^2P| = at.$ Also, $|M'/U| = |JP/J^2P|-|JM| = et-s.$
This shows that $(|M'/U|,|U|) = \omega^e_a\bdim M.$ This completes the proof of the first
formula.

Write $M' = X\oplus Y$ with $X$ bipartite and $Y$ semisimple. Then $\soc M' =
\soc X\oplus \soc Y = JX\oplus Y$ (here we use that $X$ is bipartite),
and $JM' = JX\oplus JY = JX\oplus 0 = JX.$
Since $JM' \subseteq U \subseteq \soc M'$, the direct decomposition
$\soc M' = JX\oplus Y$ yields $U = JX\oplus (Y\cap U).$
As a submodule of $Y$, the module $Y\cap U$ is a direct sum of copies of $S$.
Since $Y\cap U$ is isomorphic to $U/JM'$, we have $|Y\cap U| = |U/JM'| = w,$
thus $Y\cap U$ is isomorphic to $S^w$. Since $Y$ is semisimple, the submodule
$Y\cap U$ is a direct summand of $Y$, thus a direct summand of $M'$. This shows that
$M'$ has a direct summand of the form $S^w$, namely $Y\cap U.$

It remains to show the second assertion:
If $\Omega M$ is bipartite, then $\Omega M$ has no direct summand isomorphic to $S$, thus
$w = 0.$
$\s$
	\medskip
{\bf Remark 5.5.} The Main Lemma 5.4 focuses the attention to a direct summand of $\Omega M$
which is of the form $S^w$. However, we should stress that $S^w$ may not be the largest
semisimple direct summand of $\Omega M$, as already the case $e=1,\ a = 0$ and $M = S$
shows: here is $\Omega M = S$ and $w = 0,$ thus $S^w = 0$ (see also Remark 13.4). 
Section 13 is devoted to a discussion of $\Omega M$ and its semisimple direct summands.
	\bigskip
\vfill\eject
{\bf 5.3.  Aligned modules}
	\smallskip
Let $A$ be a short local algebra of Hilbert type $(e,a).$
We say that a module $M$ of Loewy length at most 2 is
{\it aligned}
provided $\bdim \Omega M = \omega^e_a \bdim M$. Note that {\it if $M$ is aligned,
then $|J\Omega M| =a\cdot t(M).$} Here is a reformulation of part of the Main Lemma 5.4.
	\medskip
{\bf Corollary 5.6.} {\it Let $A$ be a short local algebra and
$M$ a module of Loewy length at most $2$.  If $\Omega M$ is bipartite,
then $M$ is aligned.} $\s$
	\medskip
The converse is not true: We have $\Omega S = {}_AJ$, thus the module $S$ is always
aligned (since  $\omega^e_a\bdim S = \omega^e_a(1,0) = (e,a) = \bdim {}_AJ$), whereas
$\Omega S$ is bipartite iff $J^2 = \soc {}_AA$. In particular, for $a = 1$,
${}_AJ$ is bipartite iff $A$ is self-injective (as mentioned already in 4.10).
	\medskip
{\bf Corollary 5.7.} {\it Let $A$ be a short local algebra which is not self-injective.
Let $M$ be indecomposable, reflexive and not projective. Then $\mho M$ is aligned.}
	\medskip
Proof. According to 4.2, $M$ has Loewy length at most 2 and $\mho M$ is bipartite.
Since $M$ is torsionless, $\Omega(\mho M) = M.$
$\s$
	\bigskip
{\bf Remark 5.8.} The subsequent paper [RZ3] will provide several characterizations of the
aligned modules.
	\bigskip
{\bf 5.4. The module class $\Cal Z(q)$, where $q$ is a rational number}
	\smallskip
If $A$ is a short local algebra with $a\ge 1$,
and $q$ is a non-negative rational number, let $\Cal Z(q) = \Cal Z_A(q)$
be the class of all
indecomposable modules $M$ with Loewy length at most 2 such that
$|JM| = q\cdot a\cdot t(M).$ Note that $\Cal Z(0) = \{S\}$.
	\medskip
{\bf Lemma 5.9.}
{\it If $M \in \Cal Z(q)$ is aligned, then $\Omega M \in \Cal Z(\frac{1}{e-qa}).$}
	\medskip
Proof. If $\bdim M = (t,qat)$, then $\bdim \Omega M = (et-qat,at)= ((e-qa)t,at)$ and thus
$|J(\Omega M)| = at = \frac{1}{e-qa}\cdot a \cdot t(\Omega M).$
$\s$
	\bigskip
{\bf 5.5.  Bipartite sequences and bipartite syzygy modules}
	\smallskip
We say that an exact sequence
$$
 \epsilon\colon  \quad 0 @>>> X @>>> P @>p>> Z @>>> 0
$$
is {\it bipartite,} provided $P$ is projective, both $X,Z$ have Loewy length at most 2
and $X$ is bipartite, or, equivalently, provided $Z$ has Loewy length at most 2, $p$ is
a projective cover, and $S$ is not a direct summand of $X$.
Note that if $M$ has Loewy length at most 2, then
$\Omega M$ is bipartite if and only if the projective cover
$p\colon PM \to M$ yields a bipartite sequence $0\to \Omega M @>>> PM @>p>> M \to 0$.

Starting with a module $M$ of Loewy length at most 2, we look at all
its syzygy modules $\Omega^iM$
with $i\ge 1$. Of particular interest will be the case that
the modules $\Omega^iM$ with $1\le i \le n$ are
bipartite (thus $S$ is not a direct summand of $\Omega^iM$ for all $1\le i \le n$).
	\medskip
{\bf Corollary 5.10.} {\it Let
$M$ be of Loewy length at most $2$ and assume that there is $n\ge 1$ such that
the modules $\Omega^iM$ with $1\le i \le n$ are
bipartite. Then }
$$
 \bdim \Omega^nM = (\omega^e_r)^n \bdim M.
$$
\vglue-.5truecm
$\s$
	\bigskip
{\bf 5.6. The recursion formula}
	\smallskip
Let $M$ be a module. Recall that we write $t(M) = |\top M|.$
For $i\in \Bbb N$, let $\beta_i(M) = t(\Omega^iM)$. 
As in commutative algebra [BH, L],
one may call these numbers $\beta_i(M)$ the {\it Betti numbers} of $M$.
	\medskip
{\bf Proposition 5.11.} {\it Let $M$ be of Loewy length at most $2$.
If the modules $M$ and $\Omega M$ are aligned, then
$$
 \beta_2(M) = e\beta_1(M)-a\beta_0(M);
$$
thus either $a = 0$, or else $\beta_0(M) = \frac1a(e\beta_1(M)-\beta_2(M))$.

In particular, if $M$ is a module such that both modules
$\Omega M$ and $\Omega^2 M$ are bipartite, then $M$ and $\Omega M$ are aligned.}
	\medskip
Proof. We write $t_i = \beta_i(M) = t(\Omega^iM)$ for $0\le i \le 2$.
Let $s_1 = |J \Omega M|.$ Since $M$ is aligned, $s_1 = at_0$. Since $\Omega M$ is aligned,
$t_2 = et_1-s_1$.
Thus $t_2 = et_1-s_1 = et_1 - at_0.$

The last sentence follows from Corollary 1 in 5.5.
$\s$

	\medskip
{\bf Remark 5.12.}
In Lescot [L], modules with Loewy length at most $2$ such that the modules
$\Omega^iM$ with $1\le i \le n$ are bipartite, are called ``$n$-exceptional'' modules;
the modules which are $n$-exceptional for all $n\ge 1$ are called ``exceptional''. See
[RZ3] for a further discussion of these ``exceptional'' modules.

	\bigskip\bigskip
{\bf 6. More on reflexive modules and the proof of Theorem 1.2}
	\medskip
{\bf 6.1. The module class $\Cal Z(q)$}
	\smallskip
Recall that for $q\in \Bbb Q$, the class $\Cal Z(q)$ consists of
all the indecomposable modules $M$ of Loewy length at most $2$ such that
$|JM| = q\cdot a \cdot t(M).$
	\medskip 
{\bf Proposition 6.1.} {\it Let $A$ be a short local algebra of Hilbert type $(e,a)$.
Let $M$ be a reflexive bipartite module with $\bdim M = (t,s)$.
Then $a$ divides $s$ and $\bdim M^* = (s/a,a  t).$ Thus, if
$M \in \Cal Z(q)$, then $M^* \in \Cal Z(q^{-1}).$}
	\medskip
Proof.  We have seen already in the proof of Theorem 1.1 that $a$
divides $s$; the essential assertion is the formula for $\bdim M^*$ (but it implies,
of course, that $a$ divides $s$).

Since there exists a non-projective reflexive module $M$, we know that ${}_AJ$ is a solid
$A$-module. Since $M$ is not simple, we also know that $a\ge 1.$ Let $\Cal H$ be the set
of homomorphisms $f\colon M \to {}_AA$ with semi-simple
image (thus, these are the homomorphisms with image in $J^2$, and also the homomorphisms
with kernel containing the socle of $M$).
If $g\colon {}_AA \to {}_AA$ is the right multiplication by some element from $J$, then
$gf = 0$. This shows that $\Cal H$ is contained in the socle of $M^*$. Of course, $|\Cal H| = at.$
On the other hand, if $f\colon M \to {}_AA$ is any element of $M^*$, then
$gf(M) \subseteq g(J) \subseteq J^2$ shows that $gf$ belongs to $\Cal H.$ This shows that
$M^*/\Cal H$ is a semi-simple right $A$-module. Now $M^*$ is indecomposable and has no
simple direct summand, thus $\Cal H = \soc M^*.$

Let $u_i\colon M \to A_i = {}_AA$ be maps such that
$u = [u_1,\dots,u_z]\colon M \to \bigoplus_{i=1}^z A_i$ is
a minimal left $\add(A)$-approximation of $M$.
We can assume that $u$ is an inclusion map.
Since the cokernel of $u$ has Loewy length at most 2, we know that $J^2P$ is contained in the
socle of
$M$ and actually equal to $\soc M$. It follows that $s = |\soc M| = az.$ In particular,
$s$ is divisible by $a$.

We claim that $u_1,\dots,u_z$ is a basis of $M^*/\Cal H.$ First, we show the linear independence.
Thus, let us assume that there are scalars $\lambda_i\in k$ such that
$f = \sum_i\lambda_i u_i$ belongs to $\Cal H.$ We have to show that $\lambda_i = 0$ for all $i$.
Thus, assume that some $\lambda_i$ is non-zero, say let $\lambda_1 \neq 0$.
Let $0 \neq x\in J^2A_1$. We apply $f$ to $[x,0,\dots,0]$ and
get $f([x,0,\dots,0]) = \lambda_1x \neq 0.$ But this means that $f$ does not vanish on $\soc M$,
thus $f\notin \Cal H,$  a contradiction.

Second, we have to show that $u_1,\dots, u_z$ generate $M^*$ modulo $\Cal H.$
Let $f\colon M \to {}_AA$ be any homomorphism.
Since $u$ is a left $\add(A)$-approximation, there are maps $f_i\colon {}_AA \to {}_AA$ such that
$f = \sum_i f_iu_i.$ Write $f_i = \lambda_i\cdot 1_M+g_i$ where $\lambda_i\in k$ and
$g_i$ maps into $J$. Then $f =
\sum_i f_iu_i = \sum_i \lambda_iu_i + g$, with $g = \sum_i g_iu_i$. The image of any $u_i$
is contained in J, thus the image of $g_iu_i$ is contained in $J^2$. This shows that $g\in \Cal H.$

Altogether, we see that $u_1, \dots, u_z$ is a basis of $M^*/\Cal H$.
Since $M^*$ is bipartite, $\top M^* = M^*/\soc M^* = M^*/\Cal H.$ Therefore
$t(M^*) = |M^*/\Cal H| = z = s/a.$
	
Since $M$ is not simple, we have $s\neq 0$. We write $q = \frac s{at}$, so that
$M \in \Cal Z(q)$. Then $\bdim M^* = (s/a,at)$ shows that $M^* \in \Cal Z(q^{-1}).$
$\s$
	\medskip
{\bf Corollary 6.2.} 
{\it Let $A$ be a short local algebra of Hilbert type $(e,a)$.
Let $M$ be a reflexive bipartite module with $\bdim M = (t,at)$. 
Then $\bdim M^* = \bdim M.$} $\s$
	\bigskip
{\bf Proposition 6.3.} {\it Let $A$ be a short local algebra which is not
self-injective. Let $M$ be indecomposable, reflexive, not projective, with $\Ext^1(M,A) = 0.$
Then}
$$
 \Omega M \in \Cal Z(\frac{a+1}{e})  \ \
 \text{and} \quad
 M \in \Cal Z\bigl(\frac{e}{a+1}\bigr).
$$

We may add that we also have $\mho M \in \Cal Z(\tfrac{e^2-a-1}{ae})$.
	\medskip
Proof. Since $A$ is not self-injective,
the modules $\Omega M$ and $M$ are not simple. Also, we know that $a \ge 2$ according to Theorem 1.1.

Let $\bdim M = (z,ay).$ Therefore $\bdim \Omega M = (ez-ay,az)$, according to Lemma 5.4.
By Proposition 6.1, we have $\bdim M^* = (y,az)$ and $\bdim (\Omega M)^* = (z,aez-a^2y).$
According to Section 2.5, the $A$-dual of the
$\mho$-sequence $0 \to \Omega M \to P \to M \to 0$ is the $\mho$-sequence
$0 \to M^* \to P^*\to (\Omega M)^* \to 0$, and Lemma 5.4 asserts that
$$
 \bdim M^* = \omega^e_a \bdim (\Omega M)^* = \omega^e_a (z,aez-a^2y) = (ez-aez + a^2y,az).
$$
Altogether, we see that $(y, az) = (ez-aez + a^2y,az).$
Thus $ez-aez + a^2y = y$ and therefore $e(1-a)z = (1-a^2)y.$ Since $a\neq 1$, we see that
$y = \frac e{a+1} z$ and therefore $|\soc M| = |JM| 
= ay = \frac{ae}{a+1}z = \frac{ae}{a+1}t(M).$ This shows that $M$ belongs
to $\Cal Z(\frac{e}{a+1}).$  $\s$

	\bigskip
{\bf 6.2. Proof of Theorem 1.2} 
	\medskip
Let us recall the assertion.
	\medskip
{\it Let $A$ be a short local algebra which is not
self-injective. Assume that $M$ is an indecomposable, reflexive and
non-projective module with $\Ext^i(M,A) = 0$ for $1=1,2.$
Then $2 \le a = e-1$. If $t = t(M)$, then
$\bdim X = (t,at)$ for $X \in \{\Omega^2M,\Omega M, M, \mho M\}.$}
	\medskip
Proof. We apply Proposition 6.3 to $M$ and to $\Omega M$.
Namely, $M$ is reflexive and $\Ext^1(M,A) = 0$, thus we see that $M$ belongs to
$\Cal Z(\frac{e}{a+1})$, and $\Omega M$ belongs to $\Cal Z(\frac{a+1}{e}).$
Also, $\Omega M$ is reflexive and $\Ext^1(\Omega M,A) = 0.$ Thus we see that
$\Omega M$ belongs to $\Cal Z(\frac{e}{a+1}).$
In this way, we see that $\Omega M$ belongs to
$\Cal Z(\frac{a+1}{e})\cap \Cal Z(\frac{e}{a+1})$

But if $\Cal Z(q) \cap \Cal Z(q')$ is non-empty, then $q = q'$. It follows that
$\frac{e}{a+1} = \frac{a+1}{e},$ therefore $a = e-1$. The inequality $2 \le a$
is mentioned already in Theorem 1.1.

Let $\bdim M = (t,s)$. Since $M$ belongs to $\Cal Z(\frac{e}{a+1})$
and $a+1=e$, it follows that $s = at.$ Since the modules
$\Omega M, M, \mho M$ are aligned, and $\omega^e_a(t,at) = (t,at)$,
it follows that $\bdim X = (t,at)$ for $X \in \{\Omega^2M,\Omega M, M, \mho M\}.$
$\s$
	\bigskip\medskip
{\bf 7. The defect, defined in case $a = e-1$}
	\medskip
Since the case $a = 1$ does not provide any challenge, the interesting cases
are those with $a\ge 2$. But we include the case $a=1$ in order to point out that
the cases $a = e-1$ may be seen as having features which are similar to
the self-injective algebras of Hilbert type $(2,1)$. The self-injective algebras of
Hilbert type $(2,1)$ are exhibited in detail in Sections A.8 -- A.11 in
Appendix A: they are well related to
the Kronecker algebra $K(2)$. 
	\medskip
If $a = e-1$ and $M$ is a module of Loewy length at most 2
with $\bdim M = (t,s),$ let $\delta(M) = at-s$. We call
$\delta(M)$ the {\it defect} of $M$.
	\medskip
{\bf Lemma 7.1.} {\it Let $a = e-1\ge 1$.
Let $0 \to X \to P \to Z \to 0$ a bipartite sequence. Then $\bdim X = \bdim Z + \delta(Z)(1,1)$
and $\delta(X) = a\delta(Z)$.}
	\medskip
Proof. We have
$$
\align
 \bdim X &= (et-s,at) = ((a+1)t-s,at)  \cr &= (t,s)+(at-s,at-s)
 = \bdim Z + \delta(Z)(1,1),
\endalign
$$
and $ \delta(X) = a((a+1)t-s) - at = a^2t-as = a(at-s) = a \delta(Z).$
$\s$
	\medskip
{\bf Lemma 7.2.} {\it Let $1 \le a = e-1$.
Let $0 \to X \to P \to Z \to 0$ be  bipartite. Then the following conditions are equivalent:
\item{\rm (i)} $\delta(X) = 0$.
\item{\rm (ii)} $\delta(Z) = 0$.
\item{\rm (iii)} $\bdim X = \bdim Z$.
\item{\rm (iv)} $t(X) = t(Z)$.
\item{\rm (v)} $|JX| = |JZ|$.}
	\medskip
Proof. Since $\delta(X) = a\delta(Z)$, the conditions (i) and (ii) are equivalent.
Since $\bdim X = \bdim Z + \delta(Z)(1,1)$, the conditions (ii) and (iii) are equivalent.
Of course, (iii) implies both (iv) and (v). Now $\bdim X = \bdim Z + \delta(Z)(1,1)$
means that $t(X) =  t(Z) + \delta(Z)$ and $|JX| = |JZ| + \delta(Z).$
Thus, if (iv) of (v) is satisfied, then $\delta(Z) = 0$, thus (ii) holds.
$\s$
	\bigskip
{\bf Lemma 7.3.} {\it Let $a = e-1 \ge 1$.
If $\delta(M) = 0,$ then either $t(\Omega M) = t(M)$ and
$\delta(\Omega M) = 0$, or else $t(\Omega M) > t(M)$,
$\delta(\Omega M) > 0$ and $\Omega M$ is not bipartite.

If $\delta(M) > 0,$ then $t(\Omega M) > t(M)$ and
$\delta(\Omega M) > 0$.
	
Thus, if $\delta(M) \ge 0$ and $\delta(\Omega (M) > 0,$ then }
$$
 \cdots > \beta_{i+1}(M) > \beta_i(M) > \cdots > \beta_1(M) > \beta_0(M) = t(M).
$$

Proof. The Main Lemma 5.4 asserts that $\bdim \Omega M = \omega^e_a\bdim M+(w,-w)$
for some $w\ge 0.$

First, let $\delta(M) = 0,$ then $\bdim M = (t,at)$ for some $t>0$.
Now, $\omega^e_a(t,at) = (t,at).$ We have $\bdim \Omega M = (t,at)+(w,-w)$
for some $w \ge 0.$ If $w = 0,$ then trivially
$t(\Omega M) = t =  t(M)$ and $\delta(\Omega M) = 0)$. If $w > 0$, then
$t(\Omega M) = t+w > t =  t(M)$ and $\delta(\Omega M) =
a(t+w)-(at-w) = (a+1)w > 0.$ Also, $\Omega M$ is not bipartite,
according to Lemma 5.4.

Second, assume that $at-s = \delta(M) > 0,$ thus $at > s$.
Now $\bdim \Omega M = (et-s+w,at-w)$ for some $w\ge 0.$
Then $t(\Omega M) = et-s+w = at+t-s+w > t+w \ge t = t(M).$
Also, $a(et-s+w) = a(t+at-s+w) > a(t+w) \ge at \ge at-w$, thus
$\delta(\Omega M) > 0.$
	\smallskip
The last assertion follows by induction.
$\s$
	\medskip
{\bf Remark 7.4.} For further considerations concerning short algebras with $a = e-1$,
we refer to Sections 10, 11, 12.
	\bigskip\bigskip
{\bf 8. The syzygy modules of $S$}
	\medskip

{\bf Lemma 8.1.} {\it If $e\le a$, and $0 \to X \to P \to Z \to 0$ is a bipartite sequence, then $|\soc X| > |\soc Z|$.}
	\medskip
Proof. Let $\bdim X = (t,s)$ and $\bdim Z = (t',s').$
The Main Lemma 5.4 asserts that $(t,s) = \omega^e_a(t',s') =
(et'-s',at').$ Thus $|\soc X| = s = at' \ge et' > s'$, since $t = et'-s' > 0.$
$\s$
	\bigskip
If $(a_n)_n$ is a sequence of real numbers, we write (as usual)
$\lim_n a_n = \infty$ provided for every integer $b$ there is $N = N(b)$
such that $a_n > b$ for all $n\ge N.$
	\medskip
{\bf Proposition 8.2.} {\it Let $A$ be a short local algebra with $e\ge 2$.
Then $\lim_n \beta_n(S) = \infty,$ thus also $\lim_n |\Omega^nS| = \infty$.
If, in addition, $a< e$, then the sequence of these Betti numbers $\beta_n(S)$ of $S$
is strictly increasing: $\beta_n(S) < \beta_{n+1}(S)$ for all $n\in \Bbb N.$ }
	\medskip
Proof. For any module $M$, we have $t(M) \le |M| \le (e+a+1)t(M)$, thus
$\lim_n \beta_n(M) = \infty$ if and only if $\lim_n |\Omega^n(M)| = \infty,$
	\medskip
Let $t_n = \beta_n(S) = t(\Omega^n S).$ For $a < e$, we show that
the sequence $(t_n)_n$ is strictly increasing.
	\medskip
First, let $a = 0$. Then $\Omega^nS = S^{e^n}$ for all $n\ge 0$. Since $e\ge 2$, we have
$e^{n+1} > e^n$, thus $t_n < t_{n+1}$.
	
Second, let $1\le a \le e-1.$
We have $t_0 = 1,\ t_1 = e.$ We show by induction
that $t_{n+1} > t_n$
for all $n\ge 0.$ For $n=0,$ this holds true since $e\ge 2.$
Thus, let $n\ge 1.$ We assume that
$t_{n+1} > t_n$.
The Main Lemma 5.4 asserts that $t_{n+2} \ge et_{n+1}-at_n$. Thus
$t_{n+2} - t_{n+1} \ge
et_{n+1}-at_n - t_{n+1} = (e-1)t_{n+1}-at_n \ge at_{n+1}-at_n
= a(t_{n+1}-t_n) > 0$, where we use that $a\ge 1.$
This shows that $t_{n+2} > t_{n+1}$.
	\medskip
Finally, let $e\le a.$ We show that $\lim_n |\Omega^nS| = \infty$.
If all the modules $\Omega^n S$ are bipartite,
then Lemma 8.1 asserts that $|\soc \Omega^{n+1} S| >
|\soc \Omega^n S|$ for all $n\ge 0$, thus
$|\Omega^nS| \ge |\soc \Omega^nS| > n$ for all $n$.

It remains to consider the case that there is some
$\Omega^mS$ which is not bipartite. Let $m$ be minimal.
We claim that $\Omega^m S$ is not simple.

If $m = 1$, then ${}_AJ = \Omega S$ is of course not simple.
Let $m \ge 2.$
The minimality of $m$ implies that $Z = \Omega^{m-1}S$ is bipartite. Let
$p\colon P \to Z$ be a projective cover, thus $\Omega^m S$ is the kernel of $p$.
Since $Z$ is of Loewy length 2, we see that $J^2P$ is contained in the kernel
$\Omega^mS$ of $p$. We have $|J^2P| \ge |J^2| = a \ge 2,$ thus $|\Omega^mS| \ge 2.$
This shows that $\Omega^mS$ is not simple.
	
Thus, there is $m\ge 1$ such that
$\Omega^mS$ is neither bipartite nor simple. We
have $\Omega^mS \simeq S\oplus X$ for some  $X\neq 0.$
By induction, we have $\Omega^{bm}S
\simeq S\oplus \bigoplus_{i=0}^{b-1} \Omega^{im}X,$ for all $b\ge 0,$
thus $\Omega^{bm}S$
is the direct sum of $b+1$ non-zero modules. As a consequence,
$\Omega^{bm+i}S$
is the direct sum of $b+1$ non-zero modules, for all $i\ge 0,$  and therefore
$|\Omega^{bm+i}S| > b$ for all $i\ge 0.$ Thus, let $N(b) = bm.$
$\s$
	\medskip
{\bf Example 8.3.} {\it A short local algebra $A$ with $\beta_1(S) = \beta_2(S)$.}
In general, the Betti numbers are not strictly increasing, as the following
example shows.
Let $A$ be the  $k$-algebra generated by $x,y$ with relations
$$
{\beginpicture
    \setcoordinatesystem units <1cm,1cm>
\put{$yx,\ x^2-y^2,\ x^3$} at -5 0.5
\put{$x$} at 0 1
\put{$y$} at 1 1
\put{$x^2$\strut} at 0 0
\put{$xy$\strut} at 1 0
\arr{0 0.8}{0 0.25}
\arr{1 0.8}{1 0.25}
\arr{0.8 .8}{0.25 0.25}
\multiput{$\ss x$} at -.2 .6  1.2 .6  /
\multiput{$\ss y$} at 0.45 .7  /
\put{${}_AJ$} at -1 1
\endpicture}
$$
It is a short local algebra of Hilbert type $(2,2)$. We have $\Omega S = {}_AJ$
with dimension vector $(2,2)$. As $\Omega ({}_AJ)$ we can take
the submodule of ${}_AA^2$
generated by $[y,x]$ and $[0,y]$, and this is a free $L(2)$-module of rank 2,
thus $\bdim \Omega^2S = (2,4).$ We see that
$\beta_1(S) = 2 = \beta_2(S).$
	\bigskip\bigskip
{\bf 9. Proof of Theorems 1.3 and 1.4}
	\medskip
For the proof of Theorem 1.3, we will use the following result by Christensen and Veliche.
	\medskip
{\bf 9.1. The Christensen-Veliche Lemma} 
	\smallskip
{\bf Lemma 9.1 (Christensen-Veliche [CV]).}
{\it Let $e>0$ and $a > 1$ be integers and let $(c_i)_{i\ge0}$ be
a sequence of positive integers with
$$
  c_i = ec_{i+1}-ac_{i+2} \quad\text{for all} \quad i\ge 0.
$$
Then $a = e-1$ and $c_i = c_0$ for all $i$.}
	\medskip
Proof. See the appendix of [CV]. $\s$
	\bigskip
{\bf 9.2. Proof of Theorem 1.3}
	\smallskip
Let $A$ be a short local algebra which is not self-injective.
Since $A$ is not self-injective, we have $e\ge 2$.
Let $P_\bullet$ be a non-zero minimal complex of projective
modules which is exact.
Let $t_i$ be the rank of $P_i$ and $M_i$ the image of $d_i$.
Since $P_i$ is a projective cover of $M_i$, we have $t( M_i) = t_i$.

Note that we have $a\ge 1$. Namely, if $a = 0,$ then the modules $M_i$ are semisimple and
$\Omega S = S^e$ shows that the sequence $\cdots, t_{i+1}, t_i, \cdots$ is
strictly decreasing. Impossible.

Next, we show that $M_i$ is bipartite for $i\ll 0$.
Let $t = |M_0|$. According to Proposition 8.2, there is $N = N(b)$
such that $|\Omega^n S| > b$ for all $n \ge N.$ Let $n \ge N$ and assume that $S$ is a direct summand
of $M_{-n}$. Then $\Omega^n S$ is a direct summand of $\Omega^n M_{-n} = M_0,$ and therefore
$|\Omega^n S| \le |M_0| = b,$ a contradiction. This shows that all the modules $M_{-n}$ with $n\ge N$
are bipartite.

Using, if necessary, an index shift, we can assume that all the modules $M_i$ with
$i\le 0$ are bipartite. Let $c_i = t_{-i} =  t( M_{-i})$ for $i\ge 0$.
Since all the modules $M_{-i}$ are bipartite, Proposition 5.11 provides
the recursion formula which asserts that
$$
 c_i = ec_{i+1}-ac_{i+2}
$$
for all $i\ge 0.$ Thus we can use the Christensen-Veliche Lemma 9.1 in order to conclude
that $a = e-1$ and
that the sequence $c_0, c_1,\dots$ is constant, thus that the sequence $t_0, t_{-1}, t_{-2},\dots$
is constant.

There are two possibilities: First, all the modules $M_i$ may be bipartite. In this case,
$t_i = t_{i+1}$ for all $i\in \Bbb Z.$

Second, not all modules $M_i$ are bipartite, thus
there is a minimal index $u$ such that $M_{u+1}$ is not bipartite.
As we have seen, this implies that $t_u = t_i$ for all $i\le u.$

Since $S$ is a direct summand of $M_{u+1}$, we use again Proposition 8.2
in order to see that there is some $i\ge u$
such that $t_{i+1} > t_i$. Let $v$ be the minimal index $i$ with this property. Thus we have
$$
 t_{v+1} > t_v = t_{v-1} = \cdots .
$$
We apply Lemma 7.2 to the bipartite sequence $0 \to M_u \to P_{u-1} \to M_{u-1} \to 0$.
Since $t( M_u) = t_u = t_{u-1} = t(M_{u-1})$, it follows that $\delta(M_u) = 0.$
The first part of Lemma 7.3 yields by induction that
$\delta(M_i) = 0$ for $u\le i \le v$ and then that $\delta(M_{u+1}) > 0.$ The last part of Lemma 7.3
asserts that
$$
 \cdots >t( M_{i+1}) > t( M_i) > \cdots > t( M_{v+1}) > t( M_v)
$$
(with $i\ge v)$. This completes the proof.
 $\s$
	\bigskip
{\bf 9.3. Complexes of type I and of type II}
	\smallskip 
We will say that a complex $P_\bullet$ is of {\it type} I, provided it is
a non-zero minimal exact complex of projective
modules, and all the modules $P_i$ have the same rank.

We will say that a complex $P_\bullet$ is of {\it type} II, provided it is
a non-zero minimal exact complex of projective
modules $P_i$, and there is some integer $v$ such that
$$
 \cdots > t_{v+2} > t_{v+1} > t_v = t_{v-1} = t_{v-2} = \cdots,
$$
where $t_i$ is the rank of $P_i$.
	\medskip
{\bf Example 9.2.}
{\it An algebra $A$ of Hilbert type $(2,1)$ with $J^2\subset \soc{}_AA$
and $J^2\subset \soc A_A$ with a complex of type {\rm I.}}
	\smallskip
In contrast to the commutative case, we cannot assert
in Theorem 1.3 that $J^2 = \soc{}_AA$ or that
$J^2 = \soc A_A$, as the following example shows:
Let $A$ be the $k$-algebra with generators $x,y$ and relations $x^2, xy, y^2.$
$$
{\beginpicture
    \setcoordinatesystem units <1cm,1cm>
\put{$x^2,\ xy,\ y^2.$} at -5 0.5
\put{$x$} at 0 1
\put{$y$} at 1 1
\put{$yx$\strut} at 0 0
\arr{0 0.8}{0 0.25}
\multiput{$\ss y$} at -.2 .6   /
\put{${}_AJ$} at -1.5 1
\endpicture}
$$
Note that $y$ belongs to $\soc {}_AA$ and $x$ belongs to $\soc A_A$, but neither $x$ nor
$y$ belong to $J^2$. The ideal $J^2$ is 1-dimensional, 
whereas $\soc {}_AA$ and $\soc A_A$ are 2-dimensional. 

The complex
$$
\CD
 \cdots @>x>> {}_AA @>x>> {}_AA @>x>> \cdots
\endCD
$$
is non-zero, minimal and exact (here, $x$ denotes the right multiplication by $x$,
thus all the images are equal to $M = Ax = A/A x$).
(Note that $x$ is a left Conca element,
as defined in Section 10.3.)
	\bigskip
{\bf Example 9.3.}
{\it An algebra $A$ of Hilbert type $(3,2)$ with $J^2\neq \soc{}_AA$,
with a complex of type {\rm II.}}
	\smallskip
The algebra $A$ will be similar to the
algebra $\Lambda_0$ considered in Section 11 (and before in
[RZ1]), but with the relation $xz = 0$ instead of $xz = zx$.
To be precise: $A$ is generated by
$x, y, z$, subject to the relations:
$$
 x^2,\ y^2,\ z^2,\ xy+qyx,\ xz,\ yz,\ zy-zx,
$$
with $q\in k$ having infinite multiplicative order.
Following [RZ1], we may visualize the algebra as follows:
$$
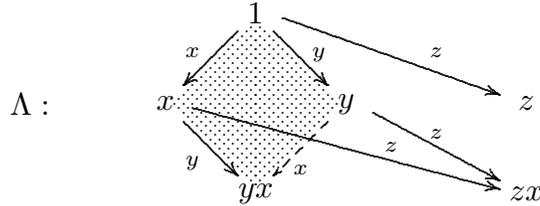
{\beginpicture
      \setcoordinatesystem units <1cm,1cm>
      \put{$\Lambda:$} at -2.5 1
      \put{$1$} at 0 2
      \put{$x$} at -1 1
      \put{$y$} at 1 1
      \put{$yx$} at 0 0
      \arr{0.2 1.8}{0.8 1.2}
      \arr{-.2 1.8}{-.8 1.2}
      \arr{-.8 0.8}{-.2 0.2}
\put{$z$} at 3 1
\put{$zx$} at 3 0
\arr{0.3 1.95}{2.7 1.1}
\arr{1.3 0.9}{2.7 0.15}
\arr{-.7 0.95}{2.7 0.05}
\setdashes <1mm>
\arr{0.8 0.8}{0.2 0.2}
\setlinear
\setshadegrid span <.5mm>
\vshade -1 1 1 <,z,,> 0 0 2 <z,,,> 1 1 1 /
\multiput{$\ssize x$\strut} at -.7  1.6  0.5 0.3 /
\multiput{$\ssize y$\strut} at 0.7  1.6  -.7 0.4  /
\multiput{$\ssize z$\strut} at 2  1.6  2 0.7  1.5 0.55 /
\endpicture}
$$
The algebra $A$ has the basis $1, \ x, \ y, \ z, \ yx,\ zx$.
We have $|\soc{}_AA| = 3$
with basis $yx, zx, z,$  whereas, of course, $|J^2| = 2.$
We get a complex of type II by taking the projective covers of the
modules $A(x-\alpha y)$
where $\alpha = q^{-i}$ with $i\ge 2$, and a minimal projective resolution of
$A(x-q^{-1}y)$. Note that $\Omega A(x-q^{-1}y) = A(x-y)\oplus S.$
	\medskip
{\bf Question 9.4.} Let $P_\bullet = (P_i,d_i)$ be a complex of type II with $M_i = \image d_i$
for all $i\in \Bbb Z$. Let
$u(P_\bullet)$ be the minimum of $i\in \Bbb Z$ with $M_{i+1}$ not bipartite (thus $S$ is a
direct summand of $M_{u+1}$, but not of $M_i$, for $i\le u$) and $v(P_\bullet)$ the minimum
of all $i\in \Bbb Z$ such that $M_i$ is not aligned (thus $M_v$ is not aligned, whereas
$M_i$ is aligned for all $i<v$). According to Corollary 1 in 5.5, we have $u\le v$,
and one may ask whether $u < v$ is possible. As we will see 
in Corollary 13.3, we have $u = v$ provided
$J^2 = \soc{}_AA.$

The question may be rephrased as follows: We look at exact sequences of the form
$$
 0 \to M' \to P \to \cdots \to   P @>>> M \to 0
$$
with $P$ a projective module which occurs $s\ge 1$ times,
with $M$ bipartite, $t(M') > t(P) = t(M)$ and $S$ a direct summand of
$\Omega M$ (namely, $M = M_u$, $P = P_u =\cdots = P_v,$ and $M' = M_{v+1}$).
The question is the following: Does there exist an exact
sequence of this kind with $s\ge 2$ such that, in addition, $M$ has a projective coresolution?
	\bigskip
{\bf 9.4. Proof of Theorem 1.4}
	\smallskip
Proof of the first part.
We assume that $A$ is a short local algebra which is not
self-injective and that there exists a module 
which is indecomposable, non-projective and
either semi-Gorenstein-projective or $\infty$-torsionfree.
Thus, there is a reflexive module which is not projective and therefore
Theorem 1.1 asserts that $a \ge 2.$
Also, there exists an $\mho$-path of length 4, thus Theorem 1.2 asserts that
$a = e-1$ and $J^2 = \soc {}_AA = \soc A_A$. This is the first part of Theorem 1.4.
$\s$
	\medskip 
Proof of (1). Let $M$ be indecomposable, non-projective, semi-Gorenstein-projective 
and torsionless. Let $t_i = t(\Omega^i M)$ for $i \ge 0.$
According to the Dictionary 2.3, $M$ is the start of an 
infinite $\mho$-path and the end of an $\mho$-path of length 1. Thus, there exists the
following $\mho$-path
$$
{\beginpicture
\setcoordinatesystem units <1.2cm,1cm>
\put{$\Omega^2M$} at -.6 0
\put{$\Omega M$} at .8 0
\put{$M$} at 2 0
\put{$\mho M$.} at 3.2 0
\setdashes <1mm>
\arr{-1.1 0}{-1.7 0}
\arr{0.3 0}{-.2 0}
\arr{1.7 0}{1.2 0}
\arr{2.7 0}{2.3 0}
\put{$\cdots$} at -2 0 
\endpicture}
$$
We see that all the modules $N = \Omega^iM$ with $i\ge 1$ are middle modules of 
$\mho$-paths of length 4, but this means that $N$ is indecomposable, reflexive, non-projective,
and satisfies $\Ext^j(N,A) = 0$ for $j = 1,2$. 
According to Theorem 1.2, we have $\bdim N = \bdim \mho N = (t_i,at_i),$ with 
$t_i = t(N)$. 
In particular, $t_{i-1} = t(\Omega^{i-1}M) = t(\mho N) = t_i.$ 

Altogether, we see that all $t_i$ with $i\ge 0$ are equal, thus equal to $t_0 = t(M)$,
and that $\bdim \Omega^i(M) = (t,at)$ for all $i\ge 0.$
$\s$
	\medskip
Proof of (2). Let $M$ be indecomposable, non-projective and $\infty$-torsionfree. 
Let $t_i = t(\mho^i M)$ for $i \ge 0.$
According to the Dictionary 2.3, $M$ is the end of an 
infinite $\mho$-path: 
$$
{\beginpicture
\setcoordinatesystem units <1.2cm,1cm>
\put{$M$} at -2 0 
\put{$\mho M$} at -.6 0
\put{$\mho^2 M$} at .8 0
\put{$\mho^3M$} at 2.3 0
\put{$\cdots$.} at 3.8 0
\setdashes <1mm>
\arr{-1.1 0}{-1.7 0}
\arr{0.3 0}{-.2 0}
\arr{1.8 0}{1.3 0}
\arr{3.3 0}{2.8 0}
\endpicture}
$$
We see that all the modules $N = \mho^iM$ with $i\ge 2$ are middle modules of 
$\mho$-paths of length 4, but this means that
$N$ is indecomposable, reflexive, non-projective,
and satisfies $\Ext^j(N,A) = 0$ for $j = 1,2$. 
According to Theorem 1.2, we have $\bdim \Omega^2N = \bdim \Omega N = (t_i,at_i),$ 
where $t_i = t(N)$. 
Now $\Omega^2 N = \mho^{i-2}N$ and $\Omega N = \mho^{i-1}N$, thus
$t(\mho^{i-2}N) = t_i = t(\mho^{i-1}N).$ 

Altogether, we see that all $t_i$ with $i\ge 0$ are equal, thus equal to $t_0 = t(M)$,
and that $\bdim \mho^i(M) = (t,at)$ for all $i\ge 0.$
$\s$
	\medskip
Proof of (3). This concerns $M^*$. First, assume that $M$ is $\infty$-torsionfree.
There is an $\mho$-sequence $0 \to M \to P \to \mho M \to 0.$ Since both $M$ and
$\mho M$ are reflexive, the $A$-dual sequence $0 \leftarrow M^*  \leftarrow P^*
\leftarrow (\mho M)^* \leftarrow 0$ is also an $\mho$-sequence. By (2), we have
$t(\mho M) = t,$ thus $P$ has rank $t$, therefore $P^*$ has rank $t$. This implies
that $t(M^*) = t.$ Since $M^*$ is bipartite, torsionless and 
semi-Gorenstein-projective, 
it follows from the right version of (1) that $\bdim M^* = (t,at).$

Second, assume that $M$ is semi-Gorenstein-projective and reflexive.
We consider an $\mho$-sequence $0 \to \Omega M \to P \to M \to 0.$ Since both
$M$ and $\Omega M$ are reflexive, also the $A$-dual $0 \leftarrow (\Omega M)^*  \leftarrow P^*
\leftarrow M^* \leftarrow 0$ is an $\mho$-sequence. Now, the rank of $P$ is $t$, thus the rank
of $P^*$ is $t$ and therefore $|\top (\Omega M)^*| = t.$ Now, $(\Omega M)^* = \mho M^*$.
Since $M^*$ is $\infty$-torsionfree, it follows from the right version of (2) that 
$\bdim M^* = \bdim \mho M^* = (t,at).$
	\medskip
Proof of (4). If $M$ is Gorenstein-projective, then $M$ is both semi-Gorenstein projective and reflexive, as well as $\infty$-torsionfree. Thus (4) follows from (1), (2) and (3), 
$\s$
	\medskip
{\bf Example 9.5.} {\it A short local algebra with an
indecomposable module $M$ which is semi-Gorenstein-projective
and torsionless (but not reflexive), with $\bdim M^* \neq \bdim M.$}
Let $A = \Lambda_0$ as discussed in Section 11 (and before in [RZ1])
and let $M$ be
the {\bf right} module $m_1A = (x-y)A$ (as above in 5.3).
The module $M$ is indecomposable, semi-Gorenstein-projective,
and torsionless (but not reflexive). We have
$(m_1A)^* = M(q)^{**} = \Omega M(1)$,
see 6.7 in [RZ1]. Therefore $\bdim\ m_1A = (1,2)$,
whereas $\bdim\ (m_1A)^* = (2,1)$.

Since $m_1A$ is torsionless, also $\mho(m_1A)$ is semi-Gorenstein-projective.
On the other hand, $\mho(m_1A)$ has Loewy length 3, see [RZ1] 7.3.
	\bigskip\bigskip
{\bf 10. Some complexes of type I}
	\medskip
{\bf 10.1. Local modules}
	\smallskip
First, let us consider local modules. Note that a module $M$
with Loewy length at most 2 is local iff $\bdim M = (1,s)$ for some natural number $s\ge 0$.
	\medskip
{\bf Lemma 10.1.} {\it Let $A$ be a short local algebra with $a = e-1$
and assume that $A$ is not self-injective.
If \ $0 \to X \to P \to Z \to 0$ is a bipartite sequence,
with $X$ a local module, then $\bdim X = \bdim Z = (1,a).$
In particular, also $Z$ is local.}
	\medskip
Proof. First, let $e = 2,$ thus $a = 1.$
Since $A$ is not self-injective, Lemma 3.2 asserts that $J^2 \subset \soc{}_AA$, thus
${}_AJ = I\oplus S$,
where $I$ is indecomposable and of length 2. Let $B$ be the factor algebra of $A$
modulo the annihilator of $I$, thus of ${}_AJ$. Then $a(B) = 0,\ e(B) = 1,$ thus
$I$ and $S$ are the only indecomposable $B$-modules. Since $X$ is cogenerated
by ${}_AJ$, it is a $B$-module. Since $X$ is bipartite, we have $X = I$, thus
$\bdim X = (1,1).$
Since the cokernel of the embedding $X \to P$ has Loewy length at most
2, we see that the projective module $P$ has rank $1$, thus
$\bdim Z = (1,1).$

Second, let $e \ge 3.$ Since $X$ is local and not simple,
$\bdim X = (1,s)$ for some $s$ with $1\le s \le e.$
According to the Main Lemma 5.4, $\bdim Z = (\frac sa,-1+\frac sa (a+1))$. It follows that $\frac sa$ has
to be an integer. Since $a \le s \le a+1$ and $a\ge 2$, it follows that $s = a$
and therefore $\bdim X = (1,a) = \bdim Z.$
$\s$
	\medskip
{\bf Remark 10.2.} Let $A$ be a short local algebra with $a = e-1$
and assume that $A$ is not self-injective. Let
$0 \to X \to P \to Z \to 0$ be a bipartite sequence.
{\it If $Z$ is a local module, then $X$ does not have to
be local.} For an example, take an algebra of the form $A = \Lambda_0$
as discussed in Section 11 (and before in [RZ1]).
Let $X$ be the submodule of $P = {}_AA$
generated by $x$ and $y$ and $Z = P/X.$ Then both $X$ and $Z$ are indecomposable of Loewy length $2$.
We have $\bdim X = (2,2),$ and $\bdim Z = (1,1),$ thus $Z$ is local whereas $X$ is not local.
Note that $\delta(X) = 2$, and $\delta(Z) = 1.$
	\medskip
{\bf Corollary 10.3.} {\it  Let $A$ be a short local algebra with $a = e-1$
and assume that $A$ is not self-injective.
If $X$ is a local reflexive module, then
$\bdim X = \bdim \mho X = (1,a).$ }
	\medskip
Proof. Since $X$ is torsionless, there is an exact sequence $\epsilon\colon \
0 \to X \to P \to \mho X \to 0$. Since $X$ is even reflexive, we know that $\mho X$
has Loewy length at most 2, thus $\epsilon$ is a bipartite sequence.
$\s$
	\bigskip
{\bf 10.2. Commutative short local algebras}
	\medskip 
We consider now the case of a commutative short local algebra with $a = e-1$.
First, let $A$ be an arbitrary commutative artinian ring.
	\medskip
{\bf Lemma 10.4.}
{\it Let $A$ be a commutative artinian ring. If
$M,\ \Omega M$ and $\Omega^2M$ are local modules, then $\Omega^3 M \simeq \Omega M$.}
	\medskip
Proof. Let $p\colon A \to M,\ p'\colon  A \to \Omega M,\ p''\colon A \to \Omega^2 M$
be projective covers. Let $u\colon \Omega M \to A$ be the kernel of $p$ and
$u'\colon \Omega^2 M \to A$ be the kernel of $p'$. Then we have $(up')(u'p'') = 0$.
Now $up', u'p''$ are right multiplications by elements of $A$. Since $A$ is
commutative, we have $(u'p'')(up') = 0,$ thus $p''u = 0$ (since $p'$ is
epi and $u'$ is mono). The sequence $0 \to \Omega M @>u>> A @>p''>> \Omega^2M \to 0$
is a short exact sequence, since $u$ is mono, $p''$ epi, $p''u = 0,$ and
$|\Omega M|+|\Omega^2 M| = |{}_AA|.$ Thus $\Omega^3 M = \Omega M.$
$\s$
	\medskip
{\bf Corollary 10.5.} {\it
Let $A$ be a commutative short local algebra.
Then any complex of type {\rm I} involving a projective module of rank $1$
is periodic of
period $2$, and there is no complex of type {\rm II} involving a projective module
of rank $1$.}
	\medskip
If $A$ is a non-commutative short local algebra, then there may exist
non-periodic complexes of type {\rm I} involving a projective module of rank 1,
as well as complexes of type {\rm II} involving a projective module of rank 1.
A typical example is the algebra $A = \Lambda_0$ discussed in Section 11 (and
before in [RZ1]).
	\bigskip
{\bf Proposition 10.6.}
{\it Let $A$ be a commutative short local algebra with $a = e-1$
and assume that $A$ is not self-injective.
If $X$ is a local module
and an $\mho$-path of length $4$ ends in $X$, then
$X$ is Gorenstein-projective with $\Omega$-period $2$
and  $\bdim \Omega X = \bdim X$.}
	\medskip
Proof. The $\mho$-path shows that the modules $X,\ \mho X, \mho^2X$ are reflexive.
Corollary 10.1 shows successively that the modules $\mho X$, then $\mho^2 X$,
finally $\mho^3 X$ are local.
We apply Lemma 10.2 to $M = \mho^3X$ (with $\Omega M = \mho^2X,\ \Omega^2 M = \mho X,\
\Omega^3 M = X$) and see that $X \simeq \mho^2 X.$ This shows that
$X$ is Gorenstein-projective with $\Omega$-period $2$. Also we see that
$\bdim \Omega X = \bdim X$.
$\s$
	\bigskip
{\bf 10.3. Conca elements}
	\smallskip
Here is a simple way for obtaining complexes of type I.
Following [AIS] (but dealing also with non-commutative local algebras),
a non-zero element $x$ will be called a {\it left Conca element}
provided $x^2 = 0$ and $J^2 = Jx$.
And $x$ is called a {\it Conca element},
provided $x^2 = 0$ and $J^2 = Jx = xJ$. If $x$ is a left Conca element,
$Ax$ is bipartite with $\bdim Ax = (1,a)$.
Let $r(x)\colon {}_AA \to {}_AA$ be the right multiplication by $x$,
defined by $r(x)(a) = ax$ for $a\in A$.
Obviously, the existence of a left Conca element implies that
$1 \le a \le e-1$ (namely, $r_x$ maps $J$ onto $J^2$ and has
$Ax = J^2+Ax$ in its kernel, thus we get a surjective map $J/Ax \to J^2$,
and $|J/Ax| = e-1,$ whereas $|J^2| = a$). In 15.1, we will see that if
$1\le a \le e-1$, there are algebras of Hilbert type $(e,a)$ with a
Conca element $x$ such that $Ax$ is reflexive. If $a = e-1$, then
for any Conca element $x$, the module $Ax$ has to be
reflexive, even Gorenstein-projective, as the following proposition shows.
	\bigskip
{\bf Proposition 10.7.} {\it Let $A$ be a short local algebra of Hilbert type $(e,e-1)$
with $e\ge 2$. Let $x$ be a left Conca element in $A$, and
$P_\bullet = (P_i,d_i)$ with $P_i = {}_AA$ and $d_i = r(x)\colon {}_AA \to {}_AA$ for all
$i\in \Bbb Z$.
Then $P_\bullet = (P_i,d_i)$  is a minimal exact complex of projective complexes with all images being equal to $Ax$. In particular, $P_\bullet$ is a complex of type I.
If $x$ is a Conca element, then also $P_\bullet^*$ is exact, thus $Ax$ is Gorenstein
projective.}
	\medskip
Proof. Since $x^2 = 0,$ we have $\image r_x \subseteq \Ker r_x$, thus $P_\bullet$ is
a complex. We have $Ax = \image r_x \subseteq \Ker r_x = \Omega (Ax)$,
and $\dim Ax = a+1$, whereas $\dim \Omega(Ax) = (1+e+a)-(a+1) = e$. Our
assumption $a = e-1$ implies that $Ax = \Omega (Ax)$, thus $P_\bullet$ is exact.
Of course, $P_\bullet$ is minimal, since $x\in J.$ Altogether we see that $P_\bullet$ is
a complex of type I with all images being equal to $Ax$.

The $A$-dual complex $P_\bullet^*$ is $(P_i^*,d_i^*)$ with $P_i^* = A_A$ and
$d_i^*:A_A \to A_A$ the left multiplication defined by $x$
(defined by $l(x)(a) = xa$ for $a\in A$). If we assume that $x$ is a Conca element, then
$x$ is a left Conca element of $A^{\op}$, therefore $P_\bullet^*$ is exact.
Altogether we see: If $x$ is a Conca element, then both  $P_\bullet$ and  $P_\bullet^*$
are exact complexes, thus $Ax$ is Gorenstein-projective.
$\s$
	\medskip
{\bf Remark 10.8.} Of course, a left Conca element is not necessarily a Conca element.
A typical example is the element $x$ in
$A = k\langle x,y\rangle/\langle x^2, xy, y^2\rangle$: Here, $Jx = kyx = J^2$, whereas
$xJ = 0.$
	\bigskip
{\bf 10.4. Exact complexes $P_\bullet^*$ with $H_i(P_\bullet^*) \neq 0$ for all 
$i\in \Bbb Z$}
	\smallskip
Next, let us draw the attention to minimal exact complexes $P_\bullet$
such that $H_i(P_\bullet^*) \neq 0$ for all $i\in \Bbb Z$. Answering questions in [CV],
Hughes-Jorgensen-\c Sega [HJS] provided examples of such complexes over
a commutative ring $A$, namely over a short local algebra of Hilbert type $(5,4)$.
In the non-commutative setting, there are such examples already
over short local algebras of Hilbert type $(2,1)$ and $(3,2).$
	\medskip
{\bf Examples 10.9.} {\it Short local algebras with
minimal exact complexes $P_\bullet$
such that $H_i(P_\bullet^*) \neq 0$ for all $i\in \Bbb Z$.}

As first example, take the algebra $A$ of Hilbert type $(2,1)$ exhibited in 9.3
and the complex $P_\bullet$ mentioned there, where $d_i\colon {}_AA \to {}_AA$ is the multiplication
by $y$ for all $i\in \Bbb Z$. All images are equal to $Ay$, thus 2-dimensional, and
therefore $P_\bullet$ is exact. In the $A$-dual complex $P_\bullet^*$, all images
are $yA,$ thus 1-dimensional. Thus $H_i(P_\bullet^*) \neq 0$ for all
$i\in \Bbb Z.$

An example $A$ with Hilbert type $(3,2)$ is the algebra $A = \Lambda_0$
as discussed in Section 11 (but also in [RZ1]; actually, one may take any
algebra of the form $\Lambda(q)$ as considered
in [RZ1], with arbitrary $q$).
Let $M = Ay.$ Then $\Omega M \simeq M$.
If $P_\bullet$ is the complex with $P_i = {}_AA$ and with
all maps $d_i\colon P_i \to P_{i-1}$ being the right multiplication by $y$, then
$P_\bullet$ is exact and minimal,  all images in $P_\bullet$ are $Ay$ (thus
bipartite), whereas all images in $P_\bullet^*$ are isomorphic to the 2-dimensional
right module $yA$ and therefore $\dim H_i(P_\bullet^*) = 2$ for all $i\in \Bbb Z.$
	\bigskip
{\bf 10.5.  Complexes of type I and of type II} 
	\smallskip
Any $\infty$-torsionfree module $M$ has a projective
coresolution which is the concatenation of $\mho$-sequences,
we may call it its {\it $\mho$-coresolution}. We may concatenate
the $\mho$-coresolution of $M$ with
a minimal projective resolution of $M$ and obtain in this way
a minimal exact complex $P_\bullet(M)$ of projective modules.
Given an $\infty$-torsionfree module $M$, one may ask whether
$P_\bullet(M)$ is a complex of type I or of type II.
	
Let us stress that both cases are possible, as the algebra $A = \Lambda_0$
considered in Section 11 (and before in [RZ1]) shows.
The $\Lambda_0$-module
$M(1)$ is $\infty$-torsionfree, and $\Omega M(1)$
has a simple direct summand, thus the minimal projective resolution
of $M(1)$ consists of projective modules whose rank is not
bounded (see Proposition 8.2), thus $P_\bullet(M(1))$ is a complex of type II.

On the other hand, if $M$ is Gorenstein-projective, then both
$P_\bullet(M)$ and $P_\bullet(M)^*$ are exact complexes of projective modules, thus
$P_\bullet(M)$ has to be a complex of type I.

But there are also $\infty$-torsionfree modules which
are not Gorenstein-projective, such that $P_\bullet(M)$ is
a complex of type I. For example, the
$\Lambda_0^\op$-module $m_{q^2}\Lambda_0$
is $\infty$-torsionfree.
Here, for $\alpha\in k$, we define $m_\alpha = x - \alpha y\in \Lambda_0$.
The syzygies of $m_{q^2}\Lambda_0$ are the modules
$m_{q^i}\Lambda_0$ with $q\le 1,$ thus of rank 1. We see that
$P_\bullet(m_{q^2}\Lambda_0)$ is a complex of type I.
	\medskip
In addition, let us remark that there are complexes $P_\bullet = (P_i,d_i)$
of type I such that the image $M$ of some $d_i$ is semi-Gorenstein-projective, but
not Gorenstein-projective. An example is the $\Lambda_0^\op$-module
$M = m_1\Lambda_0$ in [RZ1].
	\bigskip\bigskip
{\bf 11. Some short local algebras of Hilbert type $(e,e-1)$}
	\medskip
In this section, we are going to construct a short local algebra
of Hilbert type $(e,e-1)$, where $e \ge 3$, with
semi-Gorenstein-projective modules which are not Gorenstein-projective.
The algebra which we construct will be denoted by $\Lambda_c$, with $c = e-3$.
The algebras $\Lambda_0$ have been exhibited already in [RZ1] and [RZ2]
(and the general case is a straightforward generalization).

We need to assume that the base field $k$ contains 
an element $q \in k$ with infinite multiplicative order.
Thus, let $c\ge 0.$ We define $\Lambda = \Lambda_c$
by generator and relations.
The algebra $\Lambda = \Lambda_c$ is generated by
$x, y, z, u_1,\dots, u_c$, subject to the
relations:
$$
\gather
 x^2,\ y^2,\ z^2,\ yz,\ xy+qyx,\ xz-zx,\ zy-zx,\cr
  xu_i-u_ix,\ yu_i,\ u_iy,\ zu_i,\ u_iz,\ u_iu_j,
\endgather
$$
for all $1\le i,j\le c.$ We obtain a short local algebra
of Hilbert type $(3+c,2+c)$ say with radical $J$, such that
$yx, zx, u_1x, \dots, u_cx$ is a basis of $J^2 = \soc{}\Lambda J = \soc J_\Lambda$.

We may visualize (the coefficient quiver of) ${}_\Lambda J$ as follows:
$$
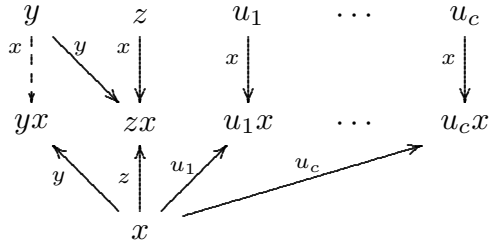
{\beginpicture
    \setcoordinatesystem units <1.2cm,1.2cm>
\put{$y$} at 0 1
\put{$z$} at 1 1
\put{$u_1$} at 2 1
\put{$\cdots$} at 3 1
\put{$u_c$} at 4 1
\put{$yx$} at 0 0
\put{$zx$} at 1 0
\put{$u_1x$} at 2 0
\put{$\cdots$} at 3 0
\put{$u_cx$} at 4 0
\put{$x$} at 1 -1
\arr{1 0.8}{1 0.2}
\arr{2 0.8}{2 0.2}
\arr{4 0.8}{4 0.2}
\arr{0.2 0.8}{0.8 0.2}
\arr{0.8 -0.8}{0.2 -0.2}
\arr{1 -0.8}{1 -0.2}
\arr{1.2 -0.8}{1.8 -0.2}
\arr{1.4 -0.85}{3.6 -0.2}
\setdashes <1mm>
\arr{0 0.8}{0 0.2}
\multiput{$\ss x$} at -.15 0.7  0.85 0.7  1.85 0.6  3.85 0.6 /
\put{$\ss y$} at .45 0.7
\put{$\ss y$} at .25 -.5
\put{$\ss z$} at .85 -.5
\put{$\ss u_1$} at 1.4 -.4
\put{$\ss u_c$} at 2.55 -.38
\endpicture}
$$
using the usual convention that a solid arrow $v \rightarrow v'$
labeled say by $x$ means that
$xv = v',$ a dashed arrow $v \dashrightarrow v'$ labeled by $x$
means that $xv$ is a non-zero multiple of $v'$ (in our case, $xy = -qyx$).
Here, the middle layer with the vertices
$yx, zx, u_1x, \dots, u_cx$ is the basis of $J^2$, as mentioned already.

We are interested in the modules $M(\alpha)$ with $\alpha\in k$ with
basis $v, v', v'', v_1,\dots, v_c$,
such that $xv = \alpha v',\ yv = v',\ zv = v'', u_iv = v_i$, for all
$1\le i \le c$
and such that
 $v', v'', v_1,\dots, v_c$ are annihilated by all generators.
$$
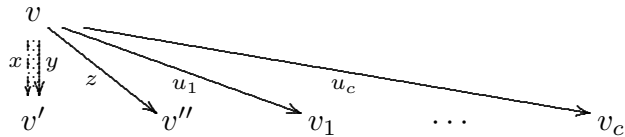
{\beginpicture
    \setcoordinatesystem units <.8cm,1cm>
    \put{\beginpicture
    \put{$v$} at 0 1.1
    \put{$v'$\strut} at 0 -.1
    \put{$v''$\strut} at 2 -.1
    \put{$v_1$\strut} at 4 -.1
    \put{$\cdots$\strut} at 5.8 -.1
    \put{$v_c$\strut} at 8 -.1
    \arr{0.08 0.8}{0.08 0.2}
    \arr{.2 .95}{1.7 0}
    \arr{.4 .95}{3.7 0}
    \arr{.7 .95}{7.7 0}
    \setdashes <1mm>
    \arr{-.08 0.8}{-.08 0.2}
    \put{$\ssize x$\strut} at   -.25 .6
    \put{$\ssize y$\strut} at    .25 .6
    \put{$\ssize z$\strut} at  .8 0.4
    \put{$\ssize u_1$\strut} at  2.1 0.35
    \put{$\ssize u_c$\strut} at  4.3 0.35
    \setshadegrid span <.5mm>
    \vshade -.15 0.2 0.81  .15 0.2 0.81 /
    \endpicture} at 0 0
    \endpicture}
    $$
The modules $M(\alpha)$ with $\alpha\in k$
are pairwise non-isomorphic indecomposable $\Lambda$-modules of length $3+c$.
As in [RZ1] one shows:
	\medskip
{\bf (1)} {\it The module $M(0)$ is Gorenstein-projective and
$\Omega$-periodic with period $1$.} In particular,
there are non-zero minimal exact complexes of projective
modules of type I.
	\medskip
{\bf (2)} {\it The module $M(q)$ is semi-Gorenstein-projective and
not torsionless.}
	\medskip
{\bf (3)} {\it The module $M(1)$ is $\infty$-torsionfree and
$\Omega M(1)$ has a simple direct summand.}
Therefore $P_\bullet(M(1))$ is a
non-zero minimal exact complex of projective
modules of type II.
	\medskip
Of course, also further properties of $\Lambda_0$ shown in [RZ1]
carry over to the algebras $\Lambda_c$ with arbitrary $c\ge 0.$
Here, we only
want to stress that for 
{\it any $a \ge 2$, there does exist
a short local algebra $A$, namely $A = \Lambda_{a-2}$, of Hilbert type $(a+1,a)$
which has modules $M,M',M''$ of length $a+1$
such that $M$ is Gorenstein-projective,
$M'$ is semi-Gorenstein-projective and not torsionless, and $M''$ is
$\infty$-torsionfree, with $\Omega M''$ having a simple direct summand.}
	\bigskip\bigskip
{\bf 12. The Auslander-Reiten conjecture (proof of Theorem 1.5)}
	\medskip
{\bf 12.1. Preliminary considerations}
	\smallskip
We need some preliminary considerations
(they are well-known, see for example Iyama [I], Section 2.1, and also
[M2]). If $M, N$ are modules, $\underline{\Hom}(\Omega M,N) = \Hom(M,N)/\Cal P(M,N)$,
where $\Cal P(M,N)$ denotes the set of homomorphisms $M \to N$ which factor through
a projective module. 
	\medskip
{\bf Lemma 12.1.}
{\it Let $\Ext^1(Z,A) = 0.$ Then, for any module $N$, we have}
$$
\align
 \Ext^1(Z,N) &\simeq \underline{\Hom}(\Omega Z,N), \tag{a} \cr
 \underline{\Hom}(Z,N)&\simeq \underline{\Hom}(\Omega Z,\Omega N).\tag{b}
\endalign
$$
	\medskip
Proof.
Let $0 \to X @>u>>  PZ \to Z \to 0$ be exact, where $PZ$ is a
projective cover of $Z$. Thus $X = \Omega Z.$

(a) We get the exact sequence
$$
 \Hom(PZ,N) @>\Hom(u,N)>> \Hom(X,N) @>\delta >> \Ext^1(Z,N) @>>> 0
$$
Of course, the image of $\Hom(u,N)$ always lies in $\Hom(X,N)_{\add A}$ (the set of
homomorphisms $X \to N$ which factor through $\add A$). Since $\Ext^1(Z,A) = 0$,
the map $u$ is a left $\add(A)$-approximation, thus any homomorphism $X \to N$
which factors through $\add(A)$ factors through $u\colon X \to PZ$. This shows that
the image of $\Hom(u,N)$ is equal to $\Hom(X,N)_{\add A}$. By definition,
$\underline{\Hom}(X,N) = \Hom(X,N)/\Hom(X,N)_{\add A}$, thus $\delta$
yields an isomorphism $\underline{\Hom}(X,N) \simeq \Ext^1(Z,N).$

(b) Let $0 \to \Omega N \to PN \to N \to 0$ be exact. Any map $f\colon Z \to N$ lifts
to a map $f'\colon PZ \to PN$ and thus yields by restriction a map $f''\colon X \to \Omega N.$
If $f$ factors though $\add A$, then $f''$ factors also through $\add A$. In this way,
we obtain an additive map $\eta\colon \underline{\Hom}(Z,N) \to \underline{\Hom}(X,\Omega N)$.
Since $u$ is a left $\add(A)$-approximation, the map $\eta$ is bijective.
$\s$
	\medskip
{\bf Lemma 12.2.} {\it If $\Ext^i(Z,A) = 0$ for $i=1,2$, then, for any module $N$}
$$
 \Ext^1(Z,N) \simeq \Ext^1(\Omega Z,\Omega N).
$$
	\medskip
Proof.
Since $\Ext^1(Z,A) = 0$, we have $\Ext^1(Z,N) \simeq \underline{\Hom}(\Omega Z,N)$.
Since $\Ext^1(\Omega Z,A) = 0$, we have
$\underline{\Hom}(\Omega Z,N)\simeq \underline{\Hom}(\Omega^2 Z,\Omega N)$ and
$\underline{\Hom}(\Omega^2 Z,\Omega N)\simeq \Ext^1(\Omega Z,\Omega N)$.
$\s$
	\bigskip
{\bf Corollary 12.3.}  {\it If $M$ is semi-Gorenstein-projective, and $N$ is an arbitrary
module, we have
$\Ext^i(M,N) \simeq \Ext^i(\Omega M,\Omega N)$,  for all $i\ge 1$.}
	\medskip
Proof. We apply Lemma 12.2 to $\Omega^{i-1}M$ and see:
$$
 \Ext^i(M,N) \simeq \Ext^1(\Omega^{i-1}M,N) \simeq
 \Ext^1(\Omega^{i}M,\Omega N) \simeq \Ext^i(\Omega M,\Omega N).
$$
\vglue-.5truecm
$\s$
	\bigskip
{\bf 12.2. The non-vanishing of $\Ext^i(M,M)$ for all $i\ge 1$}
	\smallskip
{\bf Proposition 12.4.} {\it Let $A$ be a short local algebra which is not self-injective, and let $M$ be a non-projective semi-Gorenstein-projective module.
Then $\Ext^i(M,M)\neq 0$ for {\bf all} $i\ge 1$.}
	\medskip
Proof. We can assume that $M$ is indecomposable, then also all the modules
$\Omega^i M$ are indecomposable with $i\ge 0.$
Let $(e,a)$ be the Hilbert-type
of $A$. Let $t = t(\Omega M)$.
According to Theorem 1.4, we have $a = e-1 \ge 2$ and $\bdim \Omega^i M =
(t,at)$ for all $i\ge 1.$
We have for $i\ge 1$
$$
 \Ext^i(M,M) = \Ext^i(\Omega M,\Omega M) = \Ext^1(\Omega^iM,\Omega M)
$$
where we use Corollary 12.3.
Now $\bdim \Omega^iM = (t,at) = \bdim \Omega M$, thus both modules $\Omega^iM$ and
$\Omega M$ are regular modules, see Section A.2 in Appendix A. 
Since $e\ge 3$, it follows that
$\Ext^1_{L(e)}(\Omega^iM,\Omega M) \neq 0.$  But then also
$\Ext^1_{A}(\Omega^iM,\Omega M) \neq 0,$ since $L(e)$ is a factor algebra of $A$.
(In general, of $B$ is a factor algebra of $A$, and $M',M''$ are $B$-modules, then
$\Ext^1_B(M',M'')$ can be considered as a subset of $\Ext^1_A(M',M'')$.)
$\s$
	\bigskip
{\bf 12.3. Proof of Theorem 1.5} 
	\smallskip

Let $A$ be a short local algebra and
let $M$ be a non-projective semi-Gorenstein-projective module. 
First, we consider the case that $A$ is not injective. According to Proposition 12.4
we have $\Ext^t(M,M) \neq 0$ for all $t\ge 1.$ In particular, $\Ext^1(M,M) \neq 0.$ 
Second, let $A$ be self-injective. Then we also have $\Ext^1(M,M) \neq 0$, now according
to Hoshino, see the first part of Proposition A.5 in Appendix A.
$\s$
	\bigskip\bigskip
{\bf 13. The Main Lemma, revisited}
	\medskip
{\bf 13.1. Main Lemma in the case $J^2 = \soc {}_AA$}
	\smallskip
{\bf Lemma 13.1.}
{\it Let $A$ be a short local algebra with $J^2 = \soc {}_AA.$
Let $M$ be a module of Loewy length at most $2$. Let $\Omega M = X\oplus S^w$ with
$X$ bipartite and $w\in \Bbb N.$ Then}
$$
 \bdim \Omega M = \omega^e_a \bdim M + (w,-w).
$$
	\medskip
Proof. Let $M' = \Omega M$ and take an exact sequence $0 \to M' \to P \to M \to 0$ with
$P$ projective and with an inclusion map $M' \to P$. Let
$U = J^2P.$ As in the proof of 5.4, we see that $JM' \subseteq U \subseteq \soc M'$
and that
$$
 \bdim M'  = \omega^e_a\bdim M +(w,-w).
$$
where $w = |U/JM'|.$

Now $J^2 = \soc{}_AA = \soc {}_AJ$ means that ${}_AJ^2$ is bipartite, thus also $JP$ is
bipartite. Therefore $M'\subseteq JP$ implies that $\soc M' \subseteq \soc JP = J^2P = U,$
and therefore $U = \soc M'.$

Write $M' = X\oplus W$ with $X$ bipartite and $W$ semisimple. Then $U = \soc M' =
JX\oplus W$,
and $JM' = JX\oplus JW = JX.$ Altogether, we get
$U = JM'\oplus W$. It follows that $w = |U/JM'| = |W|.$
Thus, $W$ is isomorphic to $S^w$ and therefore $M' = X\oplus W = X\oplus S^w$
with $X$ bipartite.
$\s$
	\medskip
{\bf 13.2. Consequences}
	\smallskip
Recall that a module $M$ of Loewy length at most
2 is said to be aligned (see Section 5.5), provided $\bdim \Omega M = \omega^e_a \bdim M$.
	\medskip
{\bf Corollary 13.2.} {\it Let $A$ be a short local algebra with $J^2 = \soc {}_AA$.
Then a module $M$ of Loewy length at most $2$ is aligned if and only if $\Omega M$ is bipartite.}
	\medskip
Proof. Let $M$ be a module of Loewy length at most 2.
We have seen in Corollary 5.6 that if $\Omega M$ is bipartite, then $M$ is aligned.
For the converse, we need the assumption that $J^2 = \soc {}_AA$.
By Lemma 13.1,
we know that $\Omega M = X\oplus S^w$ with $X$ bipartite and
$\bdim \Omega M = \omega^e_a \bdim M + (w,-w).$ If $M$ is aligned, then
$\bdim \Omega M = \omega^e_a \bdim M$, thus $w = 0$, and therefore $\Omega M$ is bipartite.
$\s$
 	\medskip
Using Lemma 13.1, we are able to improve Theorem 1.3 in the case $J^2 = \soc {}_AA$.
	\medskip
{\bf Corollary 13.3.} {\it
Let $A$ be a short local algebra of Hilbert type $(e,e-1)$
which is not self-injective and assume that
$J^2 = \soc {}_AA$.

Let $P_\bullet = (P_i,d_i)_i$ be a non-zero minimal exact complex of projective modules of type {\rm II}, let $M_i$ be the image of $d_i$ and $t_i = t(P_i) = t(M_i).$
As we know, there is $v\in \Bbb Z$ with $t_{v+1} > t_v = t_{v-1}$. Let $t = t_v$.
Then all the modules $M_i$ with $i\le v$ are bipartite, whereas $M_{v+1}$ is
not bipartite.}
	\medskip
Proof. By Theorem 1.3, we know that $M_{v+1}$ is not bipartite
and that $\bdim M_i = (t,at)$ for all $i\le v$.
Suppose that $M_i$ is not bipartite,
say $M_i = U\oplus S^w$ with $U$ bipartite and $w\ge 1.$
Let $M = M_{i-1}.$ According to 13.1, we have
$\bdim M_i = \bdim \Omega M = \omega^e_a\bdim M +(w,-w).$ Thus $t(M_i) = t+w > t$
and therefore $i > v.$
$\s$
	\bigskip
{\bf Remark 13.4.} Let us return to the Main Lemma 5.4 itself.
Let $M$ be a module of Loewy length at most 2.
If we use covering theory, the number $w$ provided by the Main Lemma 5.4
can be understood well. Thus, let
$\widetilde A$ be a $\Bbb Z$-cover of $A$ (we assume that the set of vertices of
the quiver of $\widetilde A$ is $\Bbb Z$, and that the
arrows go from $i$ to $i+1$, for all $i$).
Let $\pi$ be the push-down functor. Let $\widetilde M$
be a module with $\pi(\widetilde M) = M$, such that $\top \widetilde M$
is a direct sum of copies of $S(0)$ (we recall the
definition of $\widetilde M$ in Section A.2 in Appendix A).
Then $\Omega \widetilde M= U \oplus S(2)^w\oplus S(1)^{w'},$
with $U$ being bipartite (and having support equal to $\{1,2\}$ provided $U\neq 0$).
It follows that $\Omega M = \pi(\Omega \widetilde M) =
\pi(U) \oplus S^{w+w'}$. Here we see the number $w$ which is mentioned
in the Main Lemma 5.4. 
If we consider $\Omega \widetilde M$ as a representation of the
$e$-Kronecker quiver with vertices $1,2$, then $S(2)^w$ is a maximal direct summand
of $\Omega\widetilde M$ which is semisimple and projective, whereas $S(1)^{w'}$
is a maximal direct summand
of $\Omega\widetilde M$ which is semisimple and injective.
	\bigskip\medskip
{\bf 14. Algebras without non-projective reflexive modules
and without non-zero minimal exact complexes of projective modules}
	\medskip
{\bf Proposition 14.1.} {\it Let $e \ge 2$. For any $0 \le a \le e^2$, there
exists a short local algebra of Hilbert type $(e,a)$
such that any reflexive module is projective
and such that the only minimal exact complex of projective modules is the zero complex.}
	\medskip

Proof. Let $E$ be a vector space of dimension $e$ say with basis $x_1,\dots, x_e$
and let $T$ be the truncated
tensor algebra $T = k\oplus E \oplus (E\otimes E)$. Of course, $T$ is a short
local algebra with $J(T) = E\oplus (E\otimes E)$ and $J(T)^2 = E\otimes E,$ thus
$e(T) = e,\ a(T) = e^2.$

Let $0\le a \le e^2.$ We will choose a suitable subspace $U \subseteq
E\otimes E$ with $\dim U = e^2-a$ and define $A = T/U.$ Then $J(A) = J(T)/U$.
Always, $J(A) = J(T)/U$ will be decomposable, thus Theorem 1.1 asserts that $A$
has no non-projective reflexive modules.

If $a = 0,$ then we have to take $U = E\otimes E$ and obtain $A = L(e)$.
Since $e \ge 2$, $J(A) = E$ is a semisimple $A$-module of length $e$, thus decomposable.

Let $E'$ be the subspace of $E$ with basis $x = x_1,$ and
$E''$ the subspace generated by $x_2,\dots,x_e$. Thus
$E = E'\oplus E''$.

If $e \le a$, then $E\otimes E'$ has dimension $e(e-1) \ge e^2-a$, thus
there is a subspace $U\subseteq E\otimes E''$ of dimension $e^2 - a$.
Then, for $A=T/U$, we have  $J(A) = J'\oplus J''$, where $J' = E'\oplus (E\otimes E')$
and $J'' = E''\oplus (E\otimes E'')/U)$ are non-zero submodules of ${}_AJ(A)$,
thus ${}_AJ(A)$ is decomposable. Note that
$\dim J(A)^2 = \dim(E\otimes E')+\dim (E\otimes E'')/U = e + (e(e-1)-(e^2-a)) = a.$

Finally, let $1\le a < e$. Let $U'$
be the subspace of $E\otimes E$ with basis $x_{a+1}\otimes x, \dots,x_e\otimes x$,
and let $U'' = E\otimes E''.$ Let $U = U'\oplus U''.$
By abuse of notation, we will denote the
residue class of $z\in T$ modulo $U$ just by $z$ again.
We note that ${}_AJ(A)$ is the direct sum
of the local module $N$ generated by $x = x_1$ (with basis $x,\ x_1\otimes x,\dots,
x_a \otimes x$, thus $\bdim N = (1,a)$)
and a semisimple module with basis $x_2,\dots, x_e$, thus $J(A) \simeq
N\oplus S^{e-1}.$ In particular, ${}_AJ(A)$ is again decomposable.

We claim that the
only minimal exact complex of projective $A$-modules is the zero complex.
According to Theorem 1.3, we only have to look at the case $a = e-1.$
Note that $J(A)$ has the basis $x_1,\dots,x_e;\ x_1\otimes x,\dots, x_{a}\otimes x.$

The only indecomposable modules cogenerated
by ${}_AJ(A)$ are $N$ and $S$ (namely, the annihilator $C$ of ${}_AJ(A)$ is
the ideal generated by $J^2$ and the element $x_{e}$, thus
$A'' = A/C$ is of the form $L(a)$, and
${}_{A''}N$ is the indecomposable projective $L(a)$-module).

We have $\Omega S = {}_AJ(A) = N\oplus S^{e-1}$.
And we have $\Omega N = S^e$
(namely,
the map $f\colon {}_AA \to N$ with $f(1) = x$
maps $x_i$ to $x_i\otimes x$, thus its kernel has basis $x_1\otimes x,\dots, x_{a}\otimes x$
and $x_e$, thus $\Omega N$ is of the form $S^e.$)

Assume now that $P_\bullet$ is a minimal exact complex of projective modules
and that $M$ is one of the images. Then $M$ is torsionless of Loewy length at most $2$,
thus of the form $M = N^u\oplus S^v$ for some natural numbers $u,\ v \ge 0$.
We have $t(M) = u+v.$ Since
$$
 \Omega M = \Omega (N^u\oplus S^v) = S^{eu}\oplus N^v\oplus S^{(e-1)v},
$$
we have $t(\Omega M) = eu+v+(e-1)v = e(u+v).$ It follows that $t(P_{i+1}) = et(P_i)$
for all $i\in \Bbb Z.$ Since $e\ge 2$, this is only possible if
$t(P_i) = 0$ for all $i\in \Bbb Z$, thus $P_\bullet$ is the zero complex.
$\s$
	\medskip
{\bf Remark.} The assumption $e \ge 2$ is necessary, since all short
local algebras with $e = 1$ are self-injective and not semisimple
(thus, the simple module is non-projective and reflexive and occurs as an image
in a minimal exact complex of projective modules).
	\bigskip\bigskip
{\bf 15. Algebras with $\mho$-paths of length 2 and 3}
	\medskip
The existence of an $\mho$-path of length 2 
means the existence of a non-projective reflexive
module; the existence of an $\mho$-path of length 3 
means the existence of a non-projective 
$3$-torsionfree module, thus of a non-projective module $M$ such that both
$M$ and $\mho M$ are reflexive modules. 
	\medskip
{\bf 15.1. Algebras with $\mho$-paths of length 2}
	\smallskip
{\bf Proposition 15.1.}
{\it Let $1\le a \le e-1.$ There exists an (even commutative) short local
algebra $A$ of Hilbert type $(e,a)$ with a reflexive module of Loewy length $2$
with dimension vector
$(1,a)$.}
	\medskip
Proof. Let $c = e-a-1$. Let $A$ be the commutative algebra with
generators
$$
 x,\ y_1,\dots,y_a,\ z_1,\dots, z_c,
$$
and relations
$$
 x^2,\ xz_j,\ y_iy_{i'},\ y_iz_j,\ z_j^2-xy_a,\ z_jz_{j'},
$$
for all $i,i'\in \{1,\dots,a\}$ and all $j,j'\in\{1,\dots,c\}$ with $j'\neq j$.
The elements $xy_1,\dots,xy_a$ form a basis of the vector space
$J^2 = \soc {}_AA = \soc A_A$. For $a = c = 2,$ the module ${}_AJ$ looks as follows
$$
{\beginpicture
    \setcoordinatesystem units <1cm,1.5cm>
\put{$x$} at 0 1
\put{$y_1$} at 1 1
\put{$y_2$} at 2 1
\put{$z_1$} at 3 1
\put{$z_2$} at 4 1

\put{$xy_1$} at 1 0
\put{$xy_2$} at 2 0

\arr{2 0.8}{2 0.2}
\arr{1 0.8}{1 0.2}

\arr{0.05 0.8}{0.8 0.2}
\arr{0.2 0.8}{1.8 0.2}
\multiput{$\ss x$} at 1.8 .7  1.2 .7  /
\put{$\ss y_1$} at 0.2 .5
\put{$\ss y_2$} at 0.55 .8
\arr{2.8 0.8}{2.2 0.2}
\arr{3.8 0.8}{2.4 0.2}
\put{$\ss z_1$} at 2.45 .7
\put{$\ss z_2$} at 3.25 .7
\endpicture}
$$

Let $M = Ax.$ Then $M$ is a module with Loewy length $2$ and $\bdim M = (1,a)$.
Let us show that the embedding $\iota\colon Ax \to {}_AA$ is a left $\add({}_AA)$-approximation.

First, consider an element $m = \alpha x + \sum \beta_iy_i + \sum\gamma_jz_j$ with
coefficients $\alpha,\ \beta_i,\ \gamma_j\in k$ and assume that there is a
surjective map $Ax \to Am.$ We have $xm = \sum \beta_ixy_i$.
Since the element $x$ annihilates $Ax$, we must have $xm = 0$, thus $\beta_i = 0$
for all $i$. We have $z_jm = \gamma_jxy_a$. Since the element $z_j$ annihilates $x$,
we must have $\gamma_j = 0.$ It follows that $m = \alpha x.$
This shows that for any homomorphism $f\colon Ax \to {}_AJ$, there is a scalar $\alpha\in k$
such that $f-\alpha\iota$ maps into $J^2.$

Second, we show that all the maps
$g\colon Ax \to {}_AJ^2$ factor through $\iota$. Let $g(m) = \sum \delta_i xy_i$ with $\delta_i
\in k.$
Let $g'$ be the right multiplication on ${}_AA$ with $\sum\delta_iy_i$
Since
$$
 g'\iota(m) = g'(x) = x\sum\delta_iy_i =
  \sum \delta_i xy_i = g(m),
$$
it follows that $g'\iota = g.$ Altogether, we see that
$\iota$ is a left $\add({}_AA)$-approximation.

It remains to show that the factor module $\mho M = {}_AA/Ax$ is cogenerated by
${}_AJ.$ Now ${}_AA/Ax$ maps onto $Ax$ as well as onto all the modules $Az_j$
with $1\le j\le c$ and the intersection of the kernels of these maps is zero.
This shows that ${}_AA/Ax$ can be embedded into $Ax\oplus \bigoplus_j Az_j$.
$\s$
	\medskip
Note that the element $x$ constructed in the proof is a Conca element
of $A$, as defined in [AIS] (see Section 10.3, and also [RZ3]).
	\bigskip
{\bf 15.2. An algebra with an $\mho$-path of length 3}
	\smallskip
{\bf Proposition 15.2.} 
{\it There exists an (even commutative) short local algebra $A$ of Hilbert type $(6,2)$
with a non-projective $3$-torsionfree module $M$ having dimension vector $(2,2)$.}
	\medskip
Proof. 
Let $A$ be the commutative local algebra with generators $x_1,y_1,z_1,x_2,y_2,z_2$, and
with the following relations: all squares of the generators (these are 6 relations),
all products of pairs of generators with different indices (these are 9 relations), as
well as the four additional relations
$$
 y_1z_1,\ y_2z_2,\ x_1y_1 - x_2y_2,\ x_1z_1 - x_2z_2.
$$
Altogether, we have 19 relations. The ideal $J^2$ has the basis $x_1y_1, x_1z_1,$
thus the Hilbert type of $A$ is $(6,2)$.
	\medskip
We visualize $J$ as follows:
$$
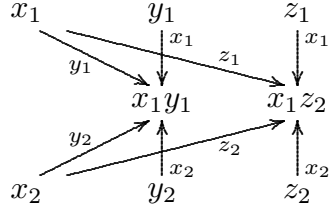
{\beginpicture
    \setcoordinatesystem units <1.5cm,1cm>
\put{$x_1$} at 0 0
\put{$y_1$} at 1 0
\put{$z_1$} at 2 0
\put{$x_1y_1$} at 1 -1
\put{$x_1z_2$} at 2 -1
\put{$x_2$} at 0 -2
\put{$y_2$} at 1 -2
\put{$z_2$} at 2 -2
\arr{0.1 -0.2}{0.9 -0.8}
\arr{0.3 -0.2}{1.9 -0.8}
\put{$\ss x_1$} at 1.15 -.3
\put{$\ss x_1$} at 2.15 -.3
\put{$\ss y_1$} at 0.4 -.6
\put{$\ss z_1$} at 1.5 -.5
\arr{1 -0.2}{1 -0.8}
\arr{2 -0.2}{2 -0.8}

\arr{0.1 -1.8}{0.9 -1.2}
\arr{0.3 -1.8}{1.9 -1.2}
\put{$\ss x_2$} at 1.15 -1.75
\put{$\ss x_2$} at 2.15 -1.75
\put{$\ss y_2$} at 0.4 -1.4
\put{$\ss z_2$} at 1.5 -1.5
\arr{1 -1.8}{1 -1.2}
\arr{2 -1.8}{2 -1.2}

\endpicture}
$$
and we may mention that $x_1$ and $x_2$ are Conca elements.
	\medskip
Let $M = Ax_1+Ax_2 \subset A$. Thus $\bdim M = (2,2).$
	\medskip
Claim: {\it The module $M$ is reflexive with $\mho M = M^*.$}
	\medskip
Proof. Let $f\colon M \to A$ be a homomorphism.

An easy calculation shows that there is $\lambda\in k$ such that
$f(c) -\lambda c \in \soc A$ for all $c\in M.$
[Namely, let $f(x_1) \equiv \alpha_1x_1+\beta_1y_1+\gamma_1z_1 + \alpha_2x_2+\beta_2y_2+\gamma_2z_2$
modulo $\soc A$, for some scalars $\alpha_i, \beta_i,\gamma_i$ ($i=1,2$).
Then $0 = f(x_1x_1) = x_1f(x_1) =
\beta_1x_1y_1 + \gamma_1x_1z_1$ shows that $\beta_1 = 0 = \gamma_1.$
Second,
 $0 = f(x_2x_1) = x_2f(x_1) =
\beta_2x_1y_1 + \gamma_2x_1z_1$ shows that $\beta_2 = 0 = \gamma_2.$
Third, $0 = f(y_2x_1) = y_2f(x_1)$ shows that $\alpha_2 = 0.$
Altogether, we see that $f(x_1) \equiv \lambda x_1$ with $\lambda = \alpha_1.$
Similarly, there is $\lambda'\in k$ with $f(x_2) \equiv \lambda' x_2$
But we also have $0 = f(x_1y_1 - x_2y_2) = y_1f(x_1)-y_2f(x_2) =
(\lambda-\lambda')x_1y_1$
thus $\lambda = \lambda'$ and therefore $f(c) \equiv \lambda c$ for all $c\in M.$]

If we use in addition that $\soc A \subset M$, we see that $f$ is the restriction to $M$ of an endomorphism of $A$.

In this way, we see that the inclusion map $M \to A$ is a (minimal)
left $(\add A)$-approximation. As a consequence, we have $\mho M = A/M.$ Now $A/M$ is the
algebra $C = L(4)$ with radical generators $y_1,z_1,y_2,z_2$. The monomorphism
$u\colon A/M \to A^2$ defined by $u(1) = (x_1,x_2)$ shows that $\mho M = A/M$ is torsionless,
therefore $M$ is reflexive.

Since $M$ is reflexive, and $\bdim M = (2,2),$ Proposition 6.1 asserts that $\bdim M^* = (1,4),$
thus $M^*$ is a local module. Since $M$ is annihilated by $x_1,x_2$, also $M^*$ is
annihilated by $x_1,x_2$, thus $M^*$ is the free module of rank 1 over the algebra $C$
(actually, the calculations presented above yield a direct way to see that $M^* = C$).
It follows that  $M^* = \mho M.$

If $M$ is reflexive, also $M^*$ is reflexive. Thus we see that in our case $\mho M = M^*$
is reflexive. Since both modules $M$ and $\mho M$ are reflexive, $M$ is $3$-torsionfree.
$\s$
	\medskip
Here is the $\mho$-path
with the dimension vectors $\bdim M,\ \bdim \mho M,\ \bdim \mho^2 M$ mentioned below
(note that the module
$\mho^3M$ has Loewy length 3).
$$
{\beginpicture
\setcoordinatesystem units <1.7cm,.7cm>
\put{\beginpicture
\put{$M$} at 0 0
\put{$\mho M$} at 1 0
\put{$\mho^2 M$} at 2 0
\put{$\mho^3 M,$} at 3 0
\setdashes <1mm>
\arr{0.7 0}{0.2 0}
\arr{1.6 0}{1.3 0}
\arr{2.6 0}{2.3 0}
\put{$\bdim $} at -1 -1
\put{$(2,2)$} at 0 -1
\put{$(1,4)$} at 1 -1
\put{$(2,11)$} at 2 -1
\endpicture} at 0 0
\endpicture}
$$
Since $M \in \Cal Z(\frac12),$ we have $\mho M \in \Cal Z(2)$, and
$\mho^2M \in \Cal Z(\frac{11}4),$  as asserted in Proposition 6.1.
	\bigskip\medskip
{\bf 16. Final remarks}
	\medskip
{\bf 16.1. The torsionless modules for a short local algebras} 
	\smallskip
The modules we have been interested in are mainly torsionless modules,
namely syzygy
modules; therefore we often have restricted
the attention to the $A$-modules of Loewy length at most 2, thus
to $L(e)$-modules, or, better, to the factor category $\mod L(e)/\add{}_AA$ (here,
we factor out the ideal of $\mod L(e)$ given by all maps which factor through a projective
$A$-module). Of course, the syzygy functor $\Omega_A$ has also to be taken into
account; it is an endo-functor of the category $\mod L(e)/\add{}_AA$.

Note that the syzygy modules in $\mod A$ are the modules cogenerated by $W = {}_AJ.$
This means: We start with an $L(e)$-module $W$
(namely the radical $W = {}_AJ$ of $A$) and look at the category $\sub W$ of all
$L(e)$-modules cogenerated by $W$, as well as at the endo-functor $\Omega_A$ of
$\sub W/\add{}_AA$.

In dealing with $L(e)$-modules $M$, the main invariant is the dimension vector
$\bdim M$; it is a pair of natural numbers, thus an element of $\Bbb Z^2$.
Here, $\Bbb Z^2$ is the Grothendieck group of the $L(e)$-modules with
respect to the exact sequences of the form $0 \to JM \to M \to M/JM \to 0,$ where
$M$ is any $L(e)$-module (equivalently, given an $L(e)$-module, we may consider the
corresponding $K(e)$-module $\widetilde M$, see Section A.2 in Appendix A,
and take as $\bdim M$
the usual dimension vector of $\widetilde M$).
As we have mentioned, the main tool in this paper
has been the transformation $\omega^e_a$ on $\Bbb Z^2$, since it describes
for the modules $M$ in $\sub W$ the dimension vector
$\bdim \Omega_A M$ in terms of $\bdim M$, at least roughly.
The transformation $\omega^e_a$ plays a role quite similar to the usual use of $\omega^e_1$
(or better of $(\omega^e_1)^2$)
in the representation theory of the $e$-Kronecker quiver (where $(\omega^e_1)^2$
describes the change
of the dimension vectors of indecomposable non-projective modules when we apply the Auslander-Reiten
translate $\tau$). A decisive difference if of course the fact that $\omega^e_1$ is invertible,
whereas, for $a\ge 2$, $\omega^e_a$ is not invertible over $\Bbb Z.$
	\medskip
{\bf 16.2. Auslander-Reiten theory and homological behavior}
	\smallskip
We want to stress that the Auslander-Reiten-quiver of an algebra $A$ 
does not determine the homological behavior of $\mod A$, see 
for example the short local self-injective algebras with $e=2$ as discussed
in Sections A.8 -- A.11 in Appendix A:
For all self-injective short local $k$-algebras with $e = 2,$
 the isomorphism classes of the indecomposable modules
 are indexed by the same set: namely, there are the  
 indecomposable $L(e)$-modules and there is one additional module,
 the projective-injective indecomposable module $P$.
The Auslander-Reiten quivers coincide: always $P$ is
 inserted at the same place.
But the homological behavior may be completely different, as the
structure of the $\mho$-components shows. 
The operator $\Omega_A$ yields an arbitrary M\"obius transformation
on the projective line $\Bbb P^1(k)$ and this transformation is not displayed by
the Auslander-Reiten quiver of $A$.
	\medskip
{\bf 16.3. Projective coresolutions} 
	\smallskip
Part of the paper
has been devoted to the study of acyclic minimal complexes of projective
modules, thus to the study of minimal projective coresolutions (of a module without
non-zero projective direct summands): Note that a minimal
projective coresolution determines uniquely an acyclic minimal complex of projective
modules 
and any acyclic minimal complex of projective modules is obtained in this way.
As we have seen, a minimal projective coresolution of a module seldom does exist.
Also, if it exists, then it may not be unique
(see for example the module $M(0,0,1)$ mentioned in [RZ2], 1.7). However, if it
exists, then its structure may be very restricted: If $A$ is a short local algebra,
and $P_0 \to P_{-1} \to P_{-2}\to \cdots$ is a non-zero minimal projective coresolution
of some module,
let $t_i = t(P_i)$. Then either $t_i = t_{i-1}$ for $i\ll 0$ (and $a = e-1$) or else
$t_{i+1}+t_{i-1} = et_i$ for all $i\ll 0$ (and $A$ is self-injective with $a = 1$),
see Theorem 1.3 and Proposition A.7 in Appendix A.
	\bigskip\medskip
{\bf Appendix A. Radical-square-zero algebras and self-injective algebras}
	\medskip
This Appendix aims to describe the categories $\mod A$ where $A$ is a
short local algebra which is self-injective (equivalently, Gorenstein, see Remark 3.3)
or has radical-square zero.
We start in A.2
with the radical-square zero $k$-algebra $A = L(e)$ (with radical $J$ of
dimension $e$ and $A/J = k$). In order to do so, we look
in A.1 at a related algebra,
the path algebra $K(e)$ of the $e$-Kronecker quiver.
	\bigskip
{\bf A.1. The structure of $\mod K(e)$}
	\medskip
We denote by $K(e)$ the {\it $e$-Kronecker quiver} with $e$ arrows (or its path algebra):
$$
{\beginpicture
    \setcoordinatesystem units <2cm,.6cm>
\multiput{$\circ$} at 0 0  1 0 /
\put{$\ss 0$} at 0 -.5
\put{$\ss 1$} at 1 -.5
\arr{0.2 0}{0.8 0}
\put{$\langle e \rangle$} at 0.5 .4
\endpicture}
$$
(here and also later, we will depict a set of $e$ arrows with same source and same sink
by a single arrow endowed with the symbol $\langle e \rangle$).
A representation (or module) $V$ of $K(e)$ will be
written in the form $V = (V_0,V_1;\ \phi\colon k^e\otimes V_0 \to V_1).$
There are two simple representations, namely $S(0) = (k,0;0)$ and $S(1)= (0,k;0)$.

The Grothendieck group of $\mod K(e)$ (with respect to exact sequences) is $\Bbb Z^2$.
Given a representation $V$ of $K(e)$, the corresponding element in the Grothendieck group
is the dimension vector $\bdim V = (\dim V_0,\dim V_1)$ of $V$.
On $\Bbb Z^2,$ we consider the quadratic form $q(x,y) = x^2 + y^2 - exy$.
This form $q$ is positive definite, if $e = 1$, it is positive semidefinite, if $e = 2$
and indefinite, if $e\ge 3$. If $e \neq 2,$ there is no non-zero pair $(x,y)$
with $q(x,y) = 0.$

We have
$q(\bdim V) = \dim\End(V) - \dim\Ext^1(V,V)$ for every module $V$ (see [R1]);
more generally, given modules $V,V'$ with $\bdim V = \bdim V'$, we have
$q(\bdim V) = \dim\Hom(V,V') - \dim\Ext^1(V,V')$.
We can use $q$ in order to distinguish between the
regular indecomposable and the non-regular indecomposable modules:
An indecomposable module $V$ is {\it regular,} provided
$\Ext^1(V,V) \neq 0$, and this happens if and only if $q(\bdim V) \le 0$.
The remaining
indecomposable modules are the indecomposable modules with $q(\bdim V) = 1$
and then $\dim \End(V)  = 1$. An element $(x,y)\in \Bbb Z^2$ is said to be a {\it real
root of} $q$ provided $q(x,y) = 1$ and an {\it imaginary root} provided $q(x,y)\le 0.$
If $V$ is a regular indecomposable module, then
there exists an indecomposable module $V'$ with $\bdim V' = \bdim V$ such that
$V$ and $V'$ are not isomorphic. If $V$ is indecomposable  with $q(\bdim V) = 1,$
then any indecomposable module $V'$ with $\bdim V' = \bdim V$ is isomorphic to $V$.
A non-regular indecomposable module $V$ with $\bdim V = (x,y)$ is said to be
{\it preprojective}
provided $x<y$, otherwise it is said to be {\it preinjective} (and then $x > y).$

For $e =1$, there are just 3 indecomposable representations, namely $S(1), P(0),
S(0)$, with $\bdim S(1) = (0,1), \bdim P(1) = (1,1)$ and $\bdim S(0) = (1,0).$

We assume now that $e\ge 2$.
The indecomposable preprojective modules can be labeled $P_0, P_1, P_2,\dots$, with
$P_0 = S(1)$, $P_1$
the indecomposable projective representation corresponding to the vertex $0$
(thus $\bdim P_1 = (1,e)$) and
$\bdim P_{i+1} = e\bdim P_i - \bdim P_{i-1}$ for $i\ge 1$.
Similarly, the indecomposable preinjective modules can be labeled $Q_0 = S(0), Q_1,
Q_2,\dots$; with $Q_0 = S(0)$, $Q_1$ the indecomposable injective representation
corresponding to the vertex $1$ (thus $\bdim Q_1 = (e,1)$)
and $\bdim Q_{i+1} = e\bdim Q_i - \bdim Q_{i-1}$ for $i\ge 1$. If we define
$b_n$ for $n\ge -1$ recursively by $b_{-1} = 0,\ b_0 = 1$ and
$b_{n+1} = eb_{n}-b_{n-1}$ for $n\ge 0$,
then $\bdim P_n = (b_{n-1},b_n)$ and $\bdim Q_n = (b_n,b_{n-1})$
(for example, for $e = 3$, the sequence $b_{-1},\, b_0,\, b_1 \dots$ is just
the sequence of the even-index
Fibonacci numbers $0,\,1,\,3,\,8,\,21,\,55,\,144,\dots$).
An explicit formula for the numbers $b_n$ due to Avramov, Iyengar and
\c Sega will be exhibited in  Appendix B.
	\medskip
The global structure of $\mod K(e)$ can be seen by looking at the Auslander-Reiten quiver of $K(e)$.
It has the following shape:
$$
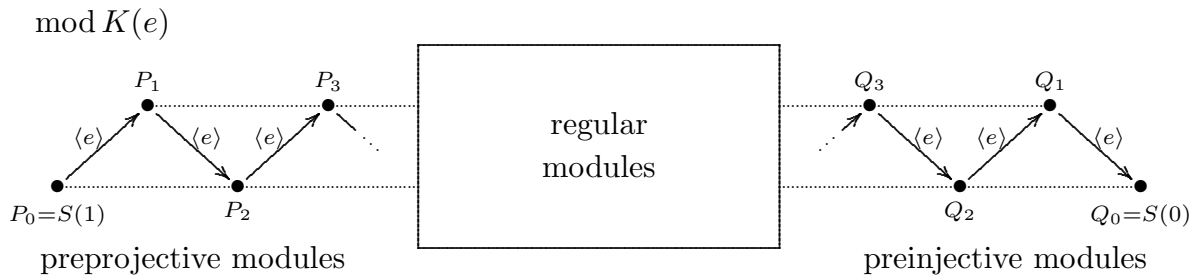
{\beginpicture
    \setcoordinatesystem units <1cm,.9cm>
\multiput{} at -5.7 -1.3  4 2.3 /
\put{$\mod K(e)$} at -5.5 2.
\setdots <.5mm>
\plot -6 0  -2 0 /
\plot -5 1  -2 1 /
\plot 6 0  2 0 /
\plot 5 1  2 1 /
\put{preprojective modules} at -4.5 -.9
\put{preinjective modules} at 4.5 -.9
\setsolid
\plot -2 -.75  2 -.75  2 1.75  -2 1.75 -2 -.75 /
\multiput{$\bullet$} at -6 0  -5 1  -4 0  -3 1  /
\multiput{$\bullet$} at 6 0  5 1  4 0  3 1  /
\put{regular} at 0 0.75
\put{modules} at 0 0.25
\put{$\ss P_0 = S(1)$} at -6 -0.35
\put{$\ss Q_0 = S(0)$} at 6 -0.35
\put{$\ss P_1$} at -5 1.3
\put{$\ss P_2$} at -4  -.3
\put{$\ss P_3$} at -3 1.3
\put{$\ss Q_1$} at 5 1.3
\put{$\ss Q_2$} at 4 -.3
\put{$\ss Q_3$} at 3 1.3
\arr{-5.9 0.1}{-5.1 0.9}
\arr{-4.9 .9}{-4.1 0.1}
\arr{-3.9 0.1}{-3.1 0.9}
\plot -2.9 0.9  -2.7 0.7 /
\arr{5.1 0.9}{5.9 0.1}
\arr{4.1 0.1}{4.9 .9}
\arr{3.1 0.9}{3.9 0.1}
\arr{2.7 0.7}{2.9 0.9}
\setdots <1mm>
\plot -2.6 0.6  -2.4 0.4 /
\plot 2.6 0.6  2.4 0.4 /

\multiput{$\ss \langle e\rangle$} at -5.65 0.6  -4.35 0.6  -3.65 0.6 /
\multiput{$\ss \langle e\rangle$} at 5.65 0.6  4.35 0.6  3.65 0.6 /
\endpicture}
$$

There are two Auslander-Reiten components of non-regular modules: the
preprojective component (seen on the left) and the preinjective component
(seen on the right). 
The Auslander-Reiten components of regular modules are homogeneous tubes for $e = 2$,
and are of the form $\Bbb Z \Bbb A_\infty$ for $e \ge 3.$

Non-zero maps between preprojective modules (and between
preinjective modules, respectively) go from left to right. Also, there are no
non-zero maps from a regular module to a preprojective module, and no
non-zero maps from a preinjective module to a preprojective or a regular module.

	\medskip
{\bf History.} Here are at least some hints.
The representations of $K(2)$ are
called  Kronecker modules, since they have been classified by
Kronecker in 1890. We will give a brief survey on related investigations at the
end of Section A.9.

The fact that there are just 3 indecomposable representations of $K(1)$
is a basic statement of elementary linear algebra.

The representation theory
of $K(e)$ with $e\ge 3$ has attracted a lot of interest in the last 40 years,
but is still very mysterious.

The algebras $K(e)$ with $e = 1,\ e = 2,$ and $e\ge 3$ are typical 
representation-finite, tame, and wild algebras, respectively. 
One expects that any one-parameter family of indecomposable modules of a tame
algebra is related to the regular modules of $K(2)$, and that any wild algebra 
has a full subcategory which is related to the regular representations of $K(3).$
	\bigskip

{\bf A.2. The push-down functor $\pi\colon \mod K(e) \to \mod L(e)$}
	\smallskip
We recall that $L(e)$ is the local $k$-algebra with radical $J$ such that
$J^2 = 0$, $\dim J = e$ and $L(e)/J = k.$
We assume here that $|J| = e \ge 2$ and identify $J = k^e$.

We denote by $\pi\colon \mod K(e) \to \mod L(e)$
the push-down functor: It sends $V = (V_0,V_1;\ \phi\colon k^e\otimes V_0 \to V_1)$ to the
representation
$$
 \pi V =
 \pi(V_0,V_1;\ \phi\colon k^e\otimes V_0 \to V_1) = \Bigl(V_0\oplus V_1; \bmatrix 0 & 0 \cr
                                                      \phi& 0 \endbmatrix\Bigr).
$$
Under the functor $\pi$, the two simple representations of $K(e)$ are sent to the
unique simple $L(e)$-module $S$.
The indecomposable $K(e)$-modules of length at least 2
correspond under $\pi$ bijectively to the indecomposable $L(e)$-modules of
length at least 2, thus to the indecomposable bipartite $L(e)$-modules.
We have $\bdim \pi V = \bdim V$ for any $K(e)$-module $V$ without a simple projective
direct summand.

Conversely, given an $L(e)$-module $M$, we denote by $\widetilde M$
the $K(e)$-module
$$
 \widetilde M = (\top M,\rad M;\  \overline\mu\colon  J\otimes \top M \to \rad M),
$$
where $\overline\mu$ is induced by the multiplication map $\mu\colon J\otimes M \to M$ 
(note that
$J\otimes \rad M$ is contained in the kernel of $\mu$ and that the image of $\mu$
is $\rad M.$). We have $\bdim \widetilde M = \bdim M$ for any $L(e)$-module $M$.

We have $\pi\widetilde M \simeq M$ for any
$L(e)$-module $M$, and conversely, we have
$\widetilde {\pi V} \simeq V$ for any $K(e)$-module $V$ without a simple projective direct summand.
Altogether we see: {\it $\pi$ and \ \ $\widetilde{}$\ \ provide
inverse bijections between isomorphism classes as follows$:$}
$$
{\beginpicture
\setcoordinatesystem units <2cm,1cm>
\put{$\left\{ \matrix \text{\rm indecomposable} \cr
                          \text{\rm $K(e)$-modules $V$}\cr
                   \text{\rm different from $S(1)$}
                   \endmatrix\right\}$} at -1 0
\put{$\left\{ \matrix \text{\rm indecomposable} \cr
                          \text{\rm $L(e)$-modules}\cr
                   \endmatrix\right\}$} at 3 0
\arr{0.5 0.1}{1.5 0.1}
\arr{1.5 -.1}{0.5 -.1}
\put{$\widetilde{}$} at 1 -.4
\put{$\pi$} at 1 .3
\endpicture}
$$
	\medskip
An indecomposable $L(e)$-module $M$ will be said to be {\it regular} provided
$\widetilde M$ is a regular $K(e)$-module.
The Auslander-Reiten quiver for $L(e)$ is
obtained from the Auslander-Reiten quiver of $K(e)$
by identifying the vertices $S(1)$ and $S(0)$ in order to obtain the vertex $S$.
$$
{\beginpicture
    \setcoordinatesystem units <1cm,1cm>
\multiput{} at -4 -2  4 4 /
\setdots <.5mm>
\ellipticalarc axes ratio 2:1 120 degrees from -4 0  center at 0 0
\ellipticalarc axes ratio 2:1 -60 degrees from -4 0  center at 0 0
\ellipticalarc axes ratio 2:1 60 degrees from  4 0  center at 0 0
\ellipticalarc axes ratio 2:1 -60 degrees from 4 0  center at 0 0

\ellipticalarc axes ratio 2:1 73 degrees from -4 1  center at 0 1
\ellipticalarc axes ratio 2:1 -60 degrees from -4 1  center at 0 1
\ellipticalarc axes ratio 2:1 60 degrees from  4 1  center at 0 1
\ellipticalarc axes ratio 2:1 -73 degrees from 4 1  center at 0 1

\multiput{$\bullet$} at 0 -2  1 -.95  2 -1.75  2.9 -.4  -1 -.95  -2 -1.75  -2.9 -.4  /
\setsolid
\arr{.9 -1.05}{0.1 -1.9}
\arr{1.9 -1.65}{1.1 -1.05}
\arr{2.85 -.5}{2.1 -1.65}
\arr{3.2 -.7}{2.95 -.45}
\arr{-0.1 -1.9}{-.9 -1.05}
\arr{-1.1 -1.05}{-1.9 -1.65}
\arr{-2.1 -1.65}{-2.85 -.5}
\plot -2.95 -.45  -3.2 -.7 /
\multiput{$\ss \langle e\rangle$} at  0.4 -1.35  1.6 -1.2  2.35  -1.0   /
\multiput{$\ss \langle e\rangle$} at  -.4 -1.35  -1.6 -1.2  -2.35  -1.0   /

\plot -2 1  2 1  2 3.5  -2 3.5 -2 1 /
\put{regular} at 0 2.5
\put{modules} at 0 2

\put{$S$} at 0 -2.3
\put{$\ss\pi P_1$} at -1.1 -.6
\put{$\ss\pi P_2$} at -2 -2.05
\put{$\ss\pi P_3$} at -2.8 -.1
\put{$\ss\pi Q_1$} at 1.1 -.6
\put{$\ss\pi Q_2$} at 2 -2.05
\put{$\ss\pi Q_3$} at 2.8 -.1
\put{$\mod L(e)$} at -5 3

	\setdots <0.3mm>
\plot -4 0 -4 1 /
\plot 4 0 4 1 /
\endpicture}
$$
	\bigskip
{\bf Proposition A.1 (Homomorphisms).} {\it If $M, M'$ are $L(e)$-modules, then $\pi$
yields an injective map
$$
 \Hom_{K(e)}(\widetilde M, \widetilde{M'}) @>\pi>> \Hom_{L(e)}(M,M')
$$
and}
$$
 \Hom_{L(e)}(M,M') = \pi\Hom_{K(e)}(\widetilde M, \widetilde{M'}) \oplus
 \Hom_k(\top M,\rad M').
$$
	\medskip
Proof. 
It is easy to show this directly. But one also may invoke the general
covering theory as developed by Gabriel and his students. We use the $\Bbb Z$-cover $Q$
of $L(e)$ with vertex set $\Bbb Z$, with $e$ arrows $z \to z\!+\!1$ \ for all $z\in \Bbb Z$
and with all paths of length 2 as relations. We identify the full subquiver of $Q$
with vertices $0,1$ with $K(e)$.

If $V$ is a representation of $Q$ and $j\in \Bbb Z$,
let $V[j]$ be the shifted representation with $V[j]_i = V_{i+j}.$ The push-down functor
$\pi$ can be extended to a functor $\pi\colon \mod Q \to \mod L(e)$ and covering theory
asserts that $\pi$ yields a bijection between
$\bigoplus_{j\in \Bbb Z} \Hom_Q(V,V[j])$ and $\Hom_{L(e)}(\pi V,\pi V').$

It remains to consider the indecomposable representations $V,V'$ of $Q$ which are either
bipartite with support $\{0,1\},$ or equal to $S(0).$ For example, if both $V,V'$ are
bipartite with support in $\{0,1\},$ then $\Hom_Q(V,V'[1]) = \Hom_k(V_0,V'_1) =
\Hom_{k}(\top V,\rad V').$
$\s$
	\bigskip
{\bf Proposition A.2 (Extensions).} {\it Let $M$ be an indecomposable regular $L(e)$-module.
Then
$$
 \Ext^1_{L(e)}(M,M) \neq 0.
$$
If $e\ge 3$ and
$M, M'$ are indecomposable regular $L(e)$ modules with
$\bdim M = \bdim M'$, then }
$$
 \Ext^1_{L(e)}(M,M') \neq 0.
$$
	\medskip
Proof. We have
$$
 \dim\End(\widetilde M) - \dim\Ext^1_{K(e)}(\widetilde M,\widetilde M) =
 q(\bdim M) \le 0.
$$
Thus $\End(\widetilde M) \neq 0$ implies that
$\Ext^1_{K(e)}(\widetilde M,\widetilde M) \neq 0$. Of course, a non-split
self-extension of $\widetilde M$ yields under $\pi$ a non-split
self-extension of the $L(e)$-module $M$.

The second assertion is shown in the same way, now using that for $e\ge 3$ we
have $q(\bdim M) < 0.$
$\s$
	\bigskip
{\bf Proposition A.3 (Solid modules).} {\it Let $M$ be an $L(e)$-module. The
following conditions are equivalent:
\item{\rm (i)} $M$ is solid.
\item{\rm(ii)} $M\neq 0$ and $\End M = k\cdot 1_M+\{\phi\in \End M\mid
 \image \phi \subseteq \rad M \subseteq \Ker \phi\}$.
\item{\rm(iii)} $\dim\End M = 1+ |\top M|\cdot|\rad M|.$
\item{\rm(iv)} $\End(\widetilde M) = k.$
\smallskip
\noindent
If these conditions are satisfied, then $M$ is indecomposable.}
	\medskip
Proof. (i) $\implies$ (ii). Assume that $M$ is solid.
An endomorphism of $M$ which does not vanish on $\soc M$
has to be invertible. In particular, $M$ has to be indecomposable: Namely, if $M= M'\oplus M''$
is a direct sum decomposition, then the projection onto $M'$ maps $\soc M'$ onto itself and
vanishes on $\soc M''.$ Thus, either $M = S$ or else $M$ is bipartite.
If $\phi$ is an endomorphism of $M$ and its restriction to $\soc M$ is the scalar multiplication
by $\lambda\in k$, then $\phi-\lambda 1_M$ maps $M$ into $\rad M$. This shows that $\End(M) = k\cdot 1_M
\oplus \Hom(\top M,\rad M)$, thus (ii) is satisfied.

(ii) $\implies$ (iii) is trivial. The implication (iii) $\implies$ (iv) is a direct consequence
of the Proposition A.1.

(iv) $\implies$ (i). Since $\widetilde M$ is indecomposable, also $M$ is indecomposable.
If $M = S$, then clearly $M$ is solid. Thus, we can assume that $M$ is bipartite. 
Proposition A.1
shows that any endomorphism $\phi$ is of the from $\phi= \lambda\cdot 1_M+\phi'$,
where $\soc M = \rad M \subseteq \Ker(\phi')$. This shows that the restriction of $\phi$
to $\soc M$ is the scalar multiplication by $\lambda.$
$\s$
	\medskip
{\bf Proposition A.4 (Modules without self-extensions).}
{\it Let $e\ge 2$. Let $M$ be an indecomposable $L(e)$-module.
The following conditions are equivalent.
\item{\rm (i)} $M$ is isomorphic to $\pi P_i$ or $\pi Q_i$ for some $i\ge 1$,
\item{\rm(ii)} $M$ is not simple and $q(\bdim M) = 1$.
\item{\rm(iii)} $\Ext^1_{L(e)}(M,M) = 0.$}
	\medskip
Proof. An indecomposable $K(e)$-module $V$ satisfies $q(\bdim V) = 1$ if and only if $V$ is
preprojective or preinjective. This yields the equivalence of (i) and (ii).

(iii) $\implies$ (i):
If $M$ is regular, then Proposition A.1 asserts that
$\Ext^1_{L(e)}(M,M) \neq 0$. If $M = S$, then, of course,
$\dim\Ext^1_{L(e)}(M,M) = e > 0.$ This shows that an indecomposable module
$M$ with $\Ext^1_{L(e)}(M,M) = 0$
is isomorphic to $\pi P_i$ or $\pi Q_i$ for some $i\ge 1$,

(i) $\implies$ (iii). Let $M$ be a bipartite module with $\bdim M = (x,y)$.
We define $g(M) = \dim\End(M) - 1 - xy.$ Since $xy = |\top M|\cdot|\rad M|,$ we see that
$g(M) \ge 0$ and that $g(M) = 0$ if and only if $M$ is solid.

The projective cover of $M$ is isomorphic to $L(e)^x$, and $\Omega M$ is semi-simple,
namely isomorphic to $S^z$ with $z = ex-y.$ We apply $\Hom(-,M)$ to the exact sequence
$0 \to S^z \to L(e)^x \to M \to 0$ and obtain the exact sequence
$$
  0 \to \Hom(M,M) \to \Hom(L(e)^x,M) \to \Hom(S^z,M) \to \Ext^1_{L(e)}(M,M) \to 0.
$$
We have $\dim\Hom(M,M) = xy+1+g(M)$, $\dim\Hom(L(e)^x,M) = x(x+y)$ and finally
$\dim\Hom(S^z,M) = zy =
(ex-y)y.$ Thus
$$
\align
 \dim\Ext^1_{L(e)}(M,M) &= xy+1+g(M) -x(x+y) + (ex-y)y \cr
   &= 1+g(M)-x^2 + exy - y^2 = 1 -q(x,y)+g(M).
\endalign
$$
If $M$ is isomorphic to $\pi P_i$ or $\pi Q_i$ for some $i\ge 1,$ then $q(x,y) = 1$
and $M$ is solid, thus $g(M) = 0,$  and therefore
$\Ext^1_{L(e)}(M,M) = 0.$
$\s$
	\bigskip
{\bf Historical remark.}
The algebra $K(e)$ is obtained from $L(e)$ by a process which has been called
``separation of a node'' by Martinez [MV1] (a {\it node} is a simple module $S$
which never occurs as a composition factor of $\rad M/(\rad M \cap \soc M)$,
for any module $M$; if the algebra is given by a quiver with relations,
then a vertex $v$ is a {\it node} iff the composition of any arrow ending in $v$ with
any arrow starting in $v$ is a relation). It seems that the first systematic separation of nodes was used in Gabriel's paper [Gb]:
He showed that using separation of the nodes,
the representations of a radical-square-zero
algebra over an algebraically closed field can be obtained from the representations
of a corresponding hereditary algebra (note that for a radical-square-zero algebra,
all simple modules are nodes). The separation of nodes yields algebras which
are stably equivalent, as later described in Auslander-Reiten-Smal\o{} [ARS, Chapter X].
	\bigskip
\vfill\eject
{\bf A.3. The self-injective short local algebras $A$ with $e\ge 2$}
	\medskip
Let $A$ be a self-injective short algebra with $e\ge 2$. We obtain
the Auslander-Reiten quiver for $A$ from the Auslander-Reiten quiver of $A/J^2$
by inserting  the vertex $A$.
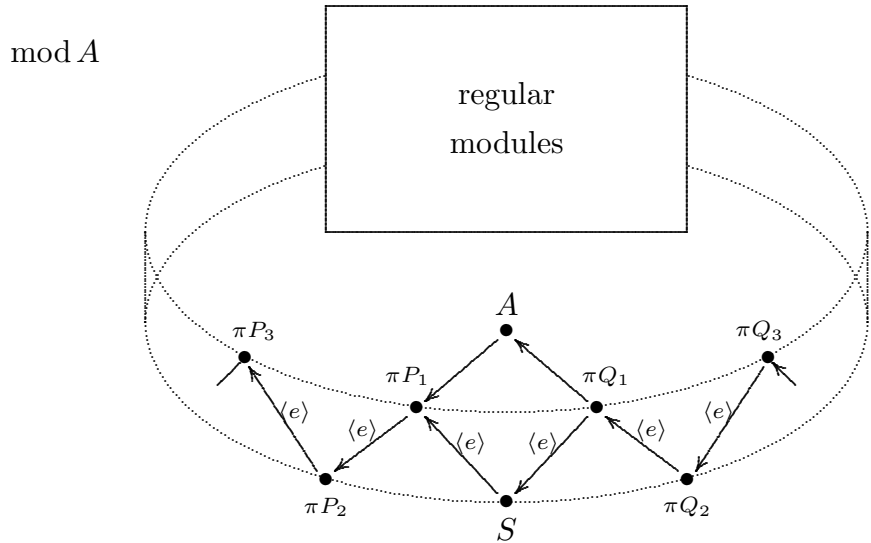
$$
{\beginpicture
    \setcoordinatesystem units <1cm,1cm>
\multiput{} at -4 -2  4 4 /
\setdots <.5mm>
\ellipticalarc axes ratio 2:1 120 degrees from -4 0  center at 0 0
\ellipticalarc axes ratio 2:1 -60 degrees from -4 0  center at 0 0
\ellipticalarc axes ratio 2:1 60 degrees from  4 0  center at 0 0
\ellipticalarc axes ratio 2:1 -60 degrees from 4 0  center at 0 0

\ellipticalarc axes ratio 2:1 120 degrees from -4 1  center at 0 1
\ellipticalarc axes ratio 2:1 -60 degrees from -4 1  center at 0 1
\ellipticalarc axes ratio 2:1 60 degrees from  4 1  center at 0 1
\ellipticalarc axes ratio 2:1 -60 degrees from 4 1  center at 0 1

\multiput{$\bullet$} at 0 -2  1 -.95  2 -1.75  2.9 -.4  -1 -.95  -2 -1.75  -2.9 -.4  0 -.1 /
\setsolid
\arr{.9 -1.05}{0.1 -1.9}
\arr{1.9 -1.65}{1.1 -1.05}
\arr{2.85 -.5}{2.1 -1.65}
\arr{3.2 -.7}{2.95 -.45}
\arr{-0.1 -1.9}{-.9 -1.05}
\arr{-1.1 -1.05}{-1.9 -1.65}
\arr{-2.1 -1.65}{-2.85 -.5}
\plot -2.95 -.45  -3.2 -.7 /

\arr{.9 -0.88}{0.1 -0.2}
\arr{-.1 -0.2}{-.9 -0.88}
\plot -2 1  2 1  2 3.5  -2 3.5 -2 1 /
\put{regular} at 0 2.5
\put{modules} at 0 2

\put{$S$} at 0 -2.3
\put{$\ss\pi P_1$} at -1.1 -.6
\put{$\ss\pi P_2$} at -2 -2.05
\put{$\ss\pi P_3$} at -2.8 -.1
\put{$\ss\pi Q_1$} at 1.1 -.6
\put{$\ss\pi Q_2$} at 2 -2.05
\put{$\ss\pi Q_3$} at 2.8 -.1
\multiput{$\ss \langle e\rangle$} at  0.4 -1.35  1.6 -1.2  2.35  -1.0   /
\multiput{$\ss \langle e\rangle$} at  -.4 -1.35  -1.6 -1.2  -2.35  -1.0   /
\put{$A$} at 0 0.2
\put{$\mod A$} at -5 3
	\setdots <0.3mm>
\plot -4 0 -4 1 /
\plot 4 0 4 1 /
\endpicture}
$$
The modules $\pi P_i$ with $i\ge 1$ are the indecomposable $A$-modules
which are different from ${}_AA$ and
preprojective in the sense of Auslander-Smal\o{} [AS].
The modules $\pi Q_i$ with $i\ge 1$ are the indecomposable $A$-modules
which are different from ${}_AA$ and
preinjective in the sense of Auslander-Smal\o{}.
	\medskip
Finally, let us present the Auslander-Reiten quiver of the triangulated category
$\underline{\mod}\ A$.
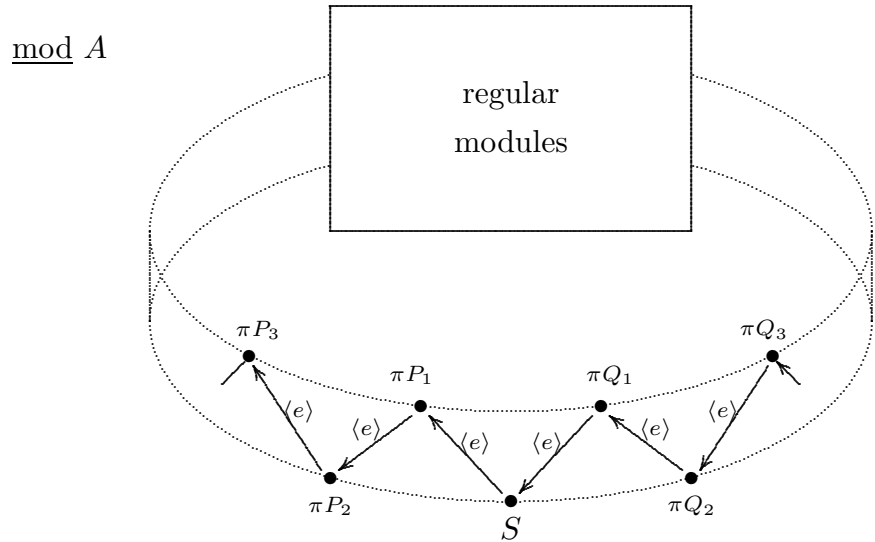
$$
{\beginpicture
    \setcoordinatesystem units <1cm,1cm>
\multiput{} at -4 -2  4 4 /
\setdots <.5mm>
\ellipticalarc axes ratio 2:1 120 degrees from -4 0  center at 0 0
\ellipticalarc axes ratio 2:1 -60 degrees from -4 0  center at 0 0
\ellipticalarc axes ratio 2:1 60 degrees from  4 0  center at 0 0
\ellipticalarc axes ratio 2:1 -60 degrees from 4 0  center at 0 0

\ellipticalarc axes ratio 2:1 120 degrees from -4 1  center at 0 1
\ellipticalarc axes ratio 2:1 -60 degrees from -4 1  center at 0 1
\ellipticalarc axes ratio 2:1 60 degrees from  4 1  center at 0 1
\ellipticalarc axes ratio 2:1 -60 degrees from 4 1  center at 0 1

\multiput{$\bullet$} at 0 -2  1 -.95  2 -1.75  2.9 -.4  -1 -.95  -2 -1.75  -2.9 -.4  /
\setsolid
\arr{.9 -1.05}{0.1 -1.9}
\arr{1.9 -1.65}{1.1 -1.05}
\arr{2.85 -.5}{2.1 -1.65}
\arr{3.2 -.7}{2.95 -.45}
\arr{-0.1 -1.9}{-.9 -1.05}
\arr{-1.1 -1.05}{-1.9 -1.65}
\arr{-2.1 -1.65}{-2.85 -.5}
\multiput{$\ss \langle e\rangle$} at  0.4 -1.35  1.6 -1.2  2.35  -1.0   /
\multiput{$\ss \langle e\rangle$} at  -.4 -1.35  -1.6 -1.2  -2.35  -1.0   /

\plot -2.95 -.45  -3.2 -.7 /

\plot -2 1  2 1  2 3.5  -2 3.5 -2 1 /
\put{regular} at 0 2.5
\put{modules} at 0 2

\put{$\ss\pi P_1$} at -1.1 -.6
\put{$\ss\pi P_2$} at -2 -2.05
\put{$\ss\pi P_3$} at -2.8 -.1
\put{$\ss\pi Q_1$} at 1.1 -.6
\put{$\ss\pi Q_2$} at 2 -2.05
\put{$\ss\pi Q_3$} at 2.8 -.1

\put{$S$} at 0 -2.3
\put{$\underline{\mod}\ A$} at -5 3
	\setdots <0.3mm>
\plot -4 0 -4 1 /
\plot 4 0 4 1 /
\endpicture}
$$
	\bigskip
\vfill\eject
{\bf A.4. The cases $e = 1$}
	\medskip
If $A$ is a self-injective short algebra with $e = 1$, then either $a = 0$ or $a = 1$.
In both cases, $A$ is uniserial, thus its module category is well understood.

It may be of interest to draw the four relevant pictures in the case $e = a = 1$,
so that one may
compare them with the pictures for $e\ge 2$ exhibited above.
Note that the last three categories shown below (the categories $\mod L(1), \ \mod A,$
and $\underline{\mod}\, A$) 
live (again) on a cylinder.
For a unified presentation, we also show $\mod K(1)$ as embedded into a cylinder ---
a rather unusual display of a single triangle.
Always, the dashed vertical lines are lines which have to be identified.
The indecomposable representation of $K(1)$ of length 2 is denoted by $I$. Of course,
if $A$ is a short local algebra of Hilbert type $(1,1),$ then $J = \rad A = \pi I.$

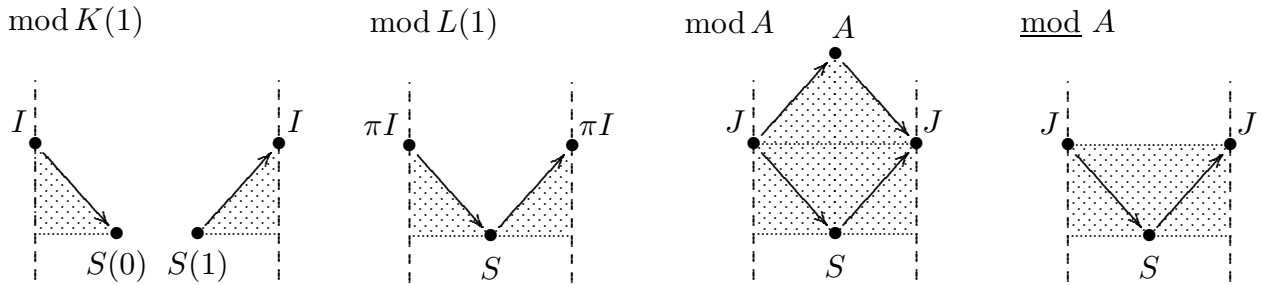
$$
{\beginpicture
    \setcoordinatesystem units <.9cm,1cm>
\setshadegrid span <.6mm>
\put{\beginpicture
\multiput{} at 0 2 /
\multiput{$\bullet$} at 0 1  1 0  2 0  3 1  /
\arr{0.1 0.9}{0.9 0.1}
\arr{2.1 0.1}{2.9 0.9}
\setdots <.5mm>
\plot 0 0  1 0 /
\plot 2 0  3 0 /
\setdashes <1mm>
\plot 0 -.5  0 1.7 /
\plot 3 -.5  3 1.7 /
\put{$\mod K(1)$} at .5 2.35
\put{$S(1)$} at 2 -.35
\multiput{$I$} at -.2 1.25  3.2 1.25 /
\put{$S(0)$} at 1 -.35
\vshade  0 0 1 <z,,,> 1 0 0 /
\vshade  2 0 0 <z,,,> 3 0 1 /
\endpicture} at -.1 0
\put{\beginpicture
\multiput{} at 0 2 /
\multiput{$\bullet$} at -1 1  0 0  1 1 /
\arr{-.9 0.9}{-.1 0.1}
\arr{.1 0.1}{.9 0.9}
\setdots <.5mm>
\plot -1 0  1 0 /
\setdashes <1mm>
\plot -1 -.5  -1 1.7 /
\plot 1 -.5  1 1.7 /
\put{$\mod L(1)$} at -.7 2.35
\put{$S$} at -0 -.35
\multiput{$\pi I$} at -1.33 1.25  1.3 1.25 /
\vshade -1 0 1 <,z,,> 0 0 0 <z,,,> 1 0 1 /
\endpicture} at 4 0
\put{\beginpicture
\multiput{} at 0 2 /
\multiput{$\bullet$} at -1 1  0 0  1 1  0 2 /
\arr{-.9 0.9}{-.1 0.1}
\arr{.1 0.1}{.9 0.9}
\arr{-.9 1.1}{-.1 1.9}
\arr{.1 1.9}{.9 1.1}
\setdots <.5mm>
\plot -1 0  1 0 /
\plot -1 1  1 1 /
\setdashes <1mm>
\plot -1 -.5  -1 1.7 /
\plot 1 -.5  1 1.7 /
\put{$\mod A$} at -1.3 2.35
\put{$S$} at -0 -.35
\put{$A$} at  0.1 2.3
\multiput{$J$} at -1.23 1.25  1.2 1.25 /
\vshade -1 0 1 <,z,,> 0 0 2 <z,,,> 1 0 1 /

\endpicture} at 8 0
\put{\beginpicture
\multiput{} at 0 2 /
\multiput{$\bullet$} at -1 1  0 0  1 1 /
\arr{-.9 0.9}{-.1 0.1}
\arr{.1 0.1}{.9 0.9}
\setdots <.5mm>
\plot -1 0  1 0 /
\plot -1 1  1 1 /
\setdashes <1mm>
\plot -1 -.5  -1 1.7 /
\plot 1 -.5  1 1.7 /
\put{$\underline{\mod}\ A$} at -1 2.35
\put{$S$} at -0 -.35
\multiput{$J$} at -1.23 1.25  1.2 1.25 /
\vshade -1 0 1  1 0 1 /
\endpicture} at 12 0
\endpicture}
$$
	\bigskip\medskip
{\bf A.5. Extensions of modules over self-injective algebras}
	\medskip
In the following proposition, the first assertion is due to Hoshino [Ho1, Theorem 3.4].
	\medskip
{\bf Proposition A.5.} {\it Let $A$ be a self-injective short local algebra.}
	\smallskip
{\bf (a)} {\bf (Hoshino [Ho1])} {\it 
 If $M$ is a non-projective module, then $\Ext^1(M,M)\neq 0.$}
	\smallskip
{\bf (b)} {\it 
If $e\neq 2$ and $M, M'$ are
non-projective indecomposable modules with $\bdim M = \bdim M'$.
Then $\Ext^1(M,M') \neq 0.$}
	\medskip
Proof. Since $M,M'$ are non-projective indecomposable modules, they have
Loewy length at most 2. Since there are non-projective modules,
we must have $e\ge 1$ and thus
$\Ext^1(S,S) \neq 0$, where $S$ is the simple $A$-module.

If $e = 1$, see Section A.4: 
Either $M$ is simple, thus $M' \simeq M$ and $\Ext^1(M,M') \neq 0$,
or else $a = 1$
and $M$ is of length $2$. Then again $M' \simeq M$ and there is an
exact sequence $0 \to M \to {}_AA\oplus S
\to M \to 0$, which shows that $\Ext^1(M,M) \neq 0,$ thus $\Ext^1(M,M') \neq 0.$

Thus, we can assume that $e\ge 2$.
If $M$ and $M'$ are regular, then $\Ext^1_{L(e)}(M,M') \neq 0,$ see Proposition A.2.
Since  there is a non-split 
exact sequence $0 \to M' \to M'' \to M \to 0$ in $\mod L(e)$, this sequence is also
a non-split exact sequence in $\mod A$, therefore
$\Ext^1_A(M,M') \neq 0.$
If $M$ is not regular, then $\bdim M = \bdim M'$ implies that $M\simeq M'$ and
$M$ belongs to the orbit of $S$ under $\Omega$ and $\Omega^{-1}.$
The Corollary 12.3 asserts that $\Ext^1(M,M) \simeq \Ext^1(S,S)\neq 0$
(since $A$ is self-injective, all modules are semi-Gorenstein-projective). 
$\s$
	\bigskip
\vfill\eject
{\bf A.6. The BGP-functors}
	\medskip
We want to show that for a self-injective short local algebra $A$ of Hilbert-type
$(e,1)$,
the syzygy functor $\Omega = \Omega_A$ corresponds to a BGP-reflection functor for
the $K(e)$-modules, as considered in [DR].
	
A BGP-functor $\sigma_\mu$
for the representations of $K(e)$
starts with two $k$-$k$-bimodules
${}_0W_1,\ _1W_0$ of dimension $e$ and a non-degenerate bilinear form
$\mu\colon _0W_1\otimes {}_1W_0\to k$. By definition,
$$
 \sigma_\mu(V_0,V_1;\ \phi\colon {}_1W_0\otimes V_0\to V_1\,) =
(\Ker \phi, \ \phi'\colon {}_0W_1\otimes \Ker\phi \to V_0),
$$
where $\phi'$ is the composition
$$
 {}_0W_1\otimes \Ker\phi @>1\otimes u>> {}_0W_1\otimes_1W_0\otimes V_0 @>\mu\otimes1>>
  k\otimes V_0 = V_0,
$$
with $u\colon \Ker\phi \to {}_1W_0\otimes V_0$ the canonical inclusion map.
We have $\sigma_\mu(S(1)) = 0$.
Let $\mod_0 K(e)$ (and $\mod_1K(e)$)
be the full subcategory of all $K(e)$-modules without simple
projective (and injective, respectively) direct summands.
The restriction of $\sigma_\mu$
to $\mod_0K(e)$ is an equivalence
$\mod_0K(e) \to \mod_1 K(e).$
If we denote the
matrix $\omega^e_1$ just by $\sigma,$ then $\bdim \sigma_\mu M = \sigma \bdim M$,
for any indecomposable $K(e)$-module $M$ which is not simple projective.
	\medskip
If $M$ is indecomposable and not
isomorphic to $S(1)$, then $\bdim \sigma_\mu M = \sigma \bdim M$. It
follows that for $e\ge 2$, we have
$$
\align
 \sigma_\mu P_i &= \left\{\matrix P_{i-1} &\text{if} & i \ge 1,\cr
                                    0     &\text{if} & i = 0, \endmatrix \right. \cr
 \sigma_\mu Q_i &= Q_{i+1} \quad \text{for all $i\ge 0$.}
\endalign
$$
	\medskip
Now we fix a self-injective algebra $A$ of Hilbert-type $(e,1)$ and an embedding
of $k^e$ as a complement of $J^2$ in $J$, thus we identify $J/J^2$ with $W = k^e$.
Let ${}_1W_2 = {}_2W_1 = W$ and take as bilinear form $\mu\colon W \otimes W \to k$
the multiplication map $J/J^2 \otimes J/J^2 \to J^2 = k.$ Since $A$ is self-injective, $\mu$
is non-degenerate and we write $\sigma_A = \sigma_\mu$.
For any $A$-module $M$, let $\Omega_A M$ be its first syzygy module.
	\medskip
{\bf Proposition A.6.} {\it Let $A$ be a self-injective short local algebra with $e(A) = e.$  
Let $M$ be in $\mod_0K(e)$. Then the $A$-module $\pi(\sigma_A M)$ is isomorphic to
$\Omega_A \pi(M)$.}
	\medskip
We have to exclude $S(1)$, since
$\pi(\sigma_A S(1))) = 0,$ whereas $\Omega_A \pi S(1) =
\Omega_A S = {}_AJ$. 
	\medskip
Proof. Let us start with the $A$-module
$M = \pi(T,\ \phi\colon W\otimes T \to JM)$, where $T = \top M = M/JM$ (thus,
we identify $M$ with $T\oplus JM$, this is the right column in the following diagram).
Its projective cover is $PM =
A\otimes T = (k\oplus W \oplus J^2)\otimes T$ (this is the middle column) with
canonical map $p = \left[\smallmatrix 1 & 0 & 0\cr
                                     0 & \phi & 0 \endsmallmatrix\right]\colon  PM\to M$.
This yields $\Omega_AM$ (namely the left column) as the kernel of
$p$.
Altogether, we deal with five exact sequences of vector spaces
(displayed in the upper five rows),
organized in two commutative diagrams. In this way, we obtain
the exact sequence of representations of $K(e)$ exhibited as the lowest row:
$$
\CD
0 @>>>   0       @>>>   k\otimes T  @>1>>        T  @>>> 0  \cr
0 @>>>  W\otimes 0       @>>>    W\otimes T @>1>>     W\otimes T  @>>> 0  \cr
@. @VVV               @VV1V                @VV\phi V \cr
0 @>>>  \Ker(\phi)  @>u>>    W\otimes T @>\phi>> JM  @>>> 0 \cr
0 @>>> W\otimes \Ker(\phi)  @>1\otimes u>> W\otimes W\otimes T
                                       @>1\otimes\phi>> W\otimes JM  @>>> 0 \cr
@. @V(\mu\otimes 1)(1\otimes u)VV               @VV\mu\otimes 1V        @VVV \cr
0 @>>>  J^2\otimes T  @>1>>  J^2\otimes T     @>>> 0 @>>> 0 \cr\cr
0 @>>>   \Omega_AM       @>>>    PM @>p>>        M  @>>> 0  \cr
\endCD
$$
	\medskip
\noindent
There is
the following commutative diagram of functors:
$$
\CD
 \mod_0 K(e) @>\tsize\sigma_A>> \mod K(e) \cr
   @VV\iota\pi V                    @VV\iota\pi V   \cr
 \mod A/\add(A)  @>\tsize\Omega_A>> \mod A/\add(A)
\endCD
$$
where $\iota\colon \mod L(e) \to \mod A$ is the canonical embedding.
	\medskip
A detailed study of the operation of $\Omega_A$ on the set of indecomposable modules
of length 2 in case $e(A) = 2$ will be given at the end of Section A.8.
	\bigskip
{\bf Historical remark.} Reflection functors for quivers were introduced by
Bernstein-Gelfand-Ponomarev [BGP] and play an important role in the representation
theory of quivers. (The printer has asked us to define BGP as used in the head line of A.6:
these are the initionals
of the authors of the paper [BGP].)
They have been generalized to species in [DR]. As we have seen
above, this generalization is actually also of interest for quivers, for example
for the $e$-Kronecker quiver $K(e)$, since one avoids in this way the use of a fixed
basis of the arrow space. But we should stress that the account given here
deviates from the usual convention (say used in [BGP] and [DR]) which is based on changing the orientation of arrows. Indeed, the BGP-reflection functors considered in [BGP] and [DR]
send a representation of the $e$-Kronecker quiver $\circ
@>\langle e\rangle >>\circ$ to a representation of
the quiver $\circ @<\langle e\rangle <<\circ$
(with opposite orientation). In contrast,
we relabel the vertices in order to obtain $\sigma_\mu$ as an endo-functor of $\mod K(e).$
As a consequence, the change of the dimension vector under $\sigma_\mu$ is described
by the product $\sigma$ of the usual BGP-reflection matrix $\left[\smallmatrix 1 & 0 \cr
                    e & -1 \endsmallmatrix \right]$
and the matrix
$\left[\smallmatrix 0 & 1 \cr
                    1 & 0 \endsmallmatrix \right]$
(corresponding to the exchange of the coordinate axes):
$$
 \sigma
  = \bmatrix 0 & 1 \cr
                    1 & 0 \endbmatrix\cdot
  \bmatrix 1 & 0 \cr
                    e & -1 \endbmatrix
 = \bmatrix e & -1 \cr
                    1 & 0 \endbmatrix
 =    \omega^e_1.
$$
	\bigskip
{\bf A.7. The $\mho$-quiver}
	\medskip
We are going to analyze the $\mho$-quivers of the short local algebras $A$
which have radical square zero or which are self-injective. 
	\medskip
{\bf A.7.1. The radical-square-zero algebras $L(e)$ with $e\ge 2$.}
First, let us look at those algebras which are not self-injective.
These are the algebras $A = L(e)$ with $e \ge 2.$ Here, 
the simple module $S$ is the only non-projective indecomposable module which is
torsionless, and $\Omega S$ is not isomorphic to $S$. 
Thus, all but one $\mho$-components are of type $\Bbb A_1$, the remaining
one is the $\mho$-component
containing $S$ and this component is of the form $\Bbb A_2$:
$$
{\beginpicture
    \setcoordinatesystem units <2cm,1cm>
\put{$[S]$} at 0 0
\put{$[\mho S]$} at 1 0
\setdashes <1mm>
\arr{0.7 0}{0.2 0}
\endpicture}
$$
The corresponding $\mho$-sequence is $0 \to S \to {}_AA^e \to \mho S \to 0,$ and
$\bdim \mho S = (e,e^2-1).$
	\bigskip
Now, we look at the self-injective algebras. 
	\medskip
{\bf A.7.2. The $\mho$-quiver of a self-injective algebra.}
If $A$ is any finite-dimensional algebra, then $A$ is self-injective algebra, iff all modules are torsionless, iff any module $M$ satisfies $\Ext^1(M,A) = 0$;
thus iff any vertex of the $\mho$-quiver of $A$ is the end of an arrow, iff
any vertex of the $\mho$-quiver of $A$ is the start of an arrow; thus iff any
$\mho$-component of $A$ is either of the form $\widetilde {\Bbb A}_n$ with $n\ge 0$
or of the form $\Bbb Z.$
For a self-injective algebra, the operator $\mho$ coincides with $\Sigma,$
where $\Sigma M$ is the cokernel of the canonical map $M \to IM$, where $IM$
is the injective envelope of $M$. It is usual also to write in this case
$\Omega^{-1}M = \Sigma M = \mho M$, since for $M$ indecomposable and not
projective, we have $\Sigma \Omega M = M = \Omega\Sigma M.$
	\medskip
{\bf A.7.3. The radical-square-zero algebra $A = L(1).$}
In this case, $S$ is the only non-projective indecomposable module, thus
the $\mho$-quiver has just one component, namely the loop 
$$
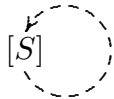
{\beginpicture
    \setcoordinatesystem units <1.5cm,1.5cm>
\multiput{} at 0 0.4  0 -.4 /
\put{$[S]$} at 0 0
\setdashes <1mm>
\circulararc 300 degrees from 0 -.15 center at  0.3  0
\setsolid
\arr{0.04 0.2}{0.01 0.15}
\endpicture}
$$

The remaining self-injective short local algebras $A$ have $a(A)\neq 0$, thus
$a(A) = 1.$ 
	\medskip
{\bf A.7.4. The case $e = 1$ and $a= 1$.}
Let $S[2]$ be the indecomposable module of length $2$, thus
$S,\ S[2]$ are the only non-projective indecomposable modules and $\Omega S = \Sigma S 
= S[2].$
$$
{\beginpicture
    \setcoordinatesystem units <2cm,1cm>
\put{$[S]$} at 0 0
\put{$[S[2]]$} at 1.1 0
\setdashes <1mm>
\setquadratic
\plot 0.2 0.2  0.5 0.3  0.8 0.2 /
\plot 0.2 -.2  0.5 -.3  0.8 -.2 /
\setsolid
\arr{0.22 0.215}{0.2 0.2}
\arr{0.8 -.215}{0.82 -.2}
\endpicture}
$$
	\medskip
For the cases $e \ge 2$, we will use A.6, in order to 
describe the $\mho$-quiver of $A$.
	\medskip
{\bf A.7.5. The regular modules.} By definition, these are the indecomposable $A$-modules
of the form $M = \pi X$, where $X$ is a regular $K(e)$-module. An $\mho$-component
which contains a regular module $M$ contains only regular modules and looks as follows:
$$
{\beginpicture
    \setcoordinatesystem units <2.2cm,1cm>
\put{$M$} at 0 0
\put{$\Omega_AM$} at -1 0
\put{$\Omega_A^2M$} at -2 0
\put{$\Omega_A^{-1}M$} at 1 0
\put{$\Omega_A^{-2}M$} at 2 0
\put{$\cdots$} at -3 0
\put{$\cdots$} at 3 0
\setdashes <.5mm>
\arr{-0.2 0}{-0.7 0}
\arr{-1.3 0}{-1.7 0}
\arr{-2.3 0}{-2.7 0}
\arr{0.65 0}{0.2 0}
\arr{1.65 0}{1.35 0}
\arr{2.7 0}{2.35 0}
\put{$\bdim M$\strut} at 0 -.6
\put{$\sigma \bdim M$\strut} at -1 -.6
\put{$\sigma^2\bdim  M$\strut} at -2 -.6
\put{$\sigma^{-1}\bdim M$\strut} at 1 -.6
\put{$\sigma^{-2}\bdim M$\strut} at 2 -.6
\endpicture}
$$
(below any module, we show the corresponding dimension vector). In general, such
an $\mho$-component is of type $\Bbb Z$.
{\it Only for $e = 2$, $M$ may be $\Omega_A$-periodic,} and then, of course,
we deal with an $\mho$-component of type $\widetilde A_n$ for some $n\ge 0.$
See Sections A.8 -- A.11 
for further discussion of the case $e = 2.$ To repeat: {\it If $e\ge 3,$
then all $\mho$-components containing regular modules are of type $\Bbb Z$} (and, as
we will see next, also the only additional component is of type $\Bbb Z$).
	\medskip 
{\bf A.7.6. The non-regular modules.} 
For all self-injective short local algabras $A$ with $e\ge 2$, there is in addition 
the $\mho$-component 
containing the simple module $S.$ It is always of type $\Bbb Z$ and
consists of $S$ and the
modules $\pi P_i$ and $\pi Q_i$ with $i\ge 1.$
We have $\pi Q_i = \Omega_A^{\;i}S$ and $\pi P_i = \Omega_A^{-i}S$;
in particular, we have $\pi Q_1 = {}_AJ,$ and $\pi P_1 = {}_AA/J^2.$
$$
{\beginpicture
    \setcoordinatesystem units <2.2cm,1cm>
\put{$S$} at 0 0
\put{$\pi Q_1$} at -1 0
\put{$\pi Q_2$} at -2 0
\put{$\pi P_1$} at 1 0
\put{$\pi P_2$} at 2 0
\put{$\cdots$} at -3 0
\put{$\cdots$} at 3 0
\setdashes <.7mm>
\arr{-0.2 0}{-0.7 0}
\arr{-1.3 0}{-1.7 0}
\arr{-2.3 0}{-2.7 0}
\arr{0.65 0}{0.2 0}
\arr{1.65 0}{1.35 0}
\arr{2.7 0}{2.35 0}

\put{$\ss\bmatrix b_2\cr b_1\endbmatrix$} at -2 -.8
\put{$\bmatrix b_1\cr b_0\endbmatrix$} at -1 -.8
\put{$\bmatrix 1\cr 0\endbmatrix$} at 0 -.8
\put{$\bmatrix b_0\cr b_1\endbmatrix$} at 1 -.8
\put{$\bmatrix b_1\cr b_2\endbmatrix$} at 2 -.8

\put{$\cdots$} at -3 -.8
\put{$\cdots$} at 3 -.8

\setdots<.8mm>
\plot 0.5 -0.3  0.5 -1.3 /
\endpicture}
$$
(again, we show below any module the corresponding dimension vector).
Since for $i\ge 0$, we have $\Omega^iS = \pi Q_i$ and $\bdim \pi Q_i =
\bdim Q_i = (b_i,b_{i-1})$, we see that
$$
 \beta_i(S) = b_i
$$
for all $i\ge 0$. This means that
{\it the numbers $b_i$ for $i\ge 0$ are just the Betti numbers of $S$.}
	
In the display of the $\mho$-component of $S$ we have inserted a dotted vertical
line between the dimension vectors of $S$ and of $\pi P_1.$ This separation line
should stress that
$\Omega(\pi P_1) = S$, whereas
$\sigma(\bdim \pi P_1) = \sigma \left[\smallmatrix b_0 \cr b_1 \endsmallmatrix\right] =
\left[\smallmatrix 0 \cr 1 \endsmallmatrix\right]
\neq  \left[\smallmatrix 1 \cr 0 \endsmallmatrix\right] =  \bdim S$.
There is just one $\mho$-sequence which is
not bipartite, namely the sequence starting in $S$ (as
mentioned already in 2.4(a)):
$$
 0 \to S \to {}_AA \to \pi P_1 \to 0.
$$
It is this sequence which is marked by the dotted separation line.
	\bigskip
{\bf A.7.7. The $\mho$-paths of length $2$.} All but two $\mho$-paths of length 2 are 
controlled by a single formula which relates the tops of the modules involved:
	\medskip
{\bf Proposition A.7.} {\it Let $A$ be a self-injective short local algebra. 
Let $M_{1} \leftarrow M_0 \leftarrow M_{-1}$ be an $\mho$-path.
If $M_0$ is not isomorphic to $S$ nor to $\pi P_1 = {}_AA/J^2$, then}
$$
  t(M_{1}) + t(M_{-1}) = e t(M_0).
$$
	\medskip
Proof. The only $\mho$-sequence which is not bipartite is the sequence
$0 \to S \to A \to A/S \to 0$. Thus, if If $M_0$ is not isomorphic to
$S$ nor to $A/J^2$,
then both sequences $0 \to M_1 \to P(M_0) \to M_0 \to 0$ and
$0 \to M_0 \to P(M_{-1}) \to M_{-1} \to 0$ are bipartite. Let $\bdim M_{-1} = (t,s)$.
Then $\bdim M_0 = (et-s,t)$ and $\bdim M_1 = (e(et-s)-t,et-s)$. Since
$t(M_1) = e(et-s)-t,\ t(M_0) = et-s,\ t(M_{-1}) = t$, we see that
$ t(M_{1}) + t(M_{-1}) = e M_0.$
$\s$
	\medskip
There are the two remaining $\mho$-paths $M_1\leftarrow M_0 \leftarrow M_{-1}$
with $M_0 = S$ and $M_0 = {}_AA/J^2.$ Both are part of the $\mho$-component which
contains $S$. This $\mho$-component has been displayed above. Let us show again the
relevant part:
$$
{\beginpicture
    \setcoordinatesystem units <2cm,1cm>
\put{$S$} at 0 0
\put{$\pi Q_1$} at -1 0
\put{$\pi P_1$} at 1 0
\put{$\pi P_2$} at 2 0
\setdashes <.7mm>
\arr{-.3 0}{-.7 0}
\arr{0.65 0}{0.2 0}
\arr{1.65 0}{1.35 0}
\put{$={}_AJ$} at -1 -.4
\put{$={}_AA/J^2$} at 1 -.4
\put{$=\mho({}_AA/J^2)$} at 2 -.4
\setdashes<1mm>
\setdots<.8mm>
\plot 0.5 -.5  0.5 -.1 /
\plot 0.5 .5  0.5 .1 /
\endpicture}
$$
and recall that
$t({}_AJ) = b_1 = e,\ t(S) = b_0 = 1,\ t({}_AA/J^2) = b_0 = 1,\
t(\mho({}_AA/J^2)) = b_1 = e$.
	\medskip
{\bf Corollary A.8.} {\it Let $e \ge 2.$
Let $P_\bullet$ be a minimal exact complex of projective modules
and let  $t_i = t(P_i)$.
If all images of $P_\bullet$ are bipartite, then
$$
 t_{i-1}+t_{i+1} = e t_i \tag{$*$}
$$
for all $i\in \Bbb Z$. If $S$ is the image of $P_0 \to P_{-1}$, then $(*)$
holds for all $i\notin\{0,-1\}$ and $t_{-1} = t_0 = 1,\ t_{-2} = t_{1} = e.$} $\s$
	\bigskip
{\bf Historical Remark.}
For a self-injective algebra $A$, the $\mho$-quiver just
depicts the graph of the operation $\Omega$ on the set of isomorphism classes of
indecomposable non-projective modules, thus it visualizes a basic concept which has been
used since the early days of homological algebra.
	\bigskip\bigskip

{\bf A.8. Self-injective algebras with $e = 2:$ The modules of length $2$}
	\medskip
Let $A$ be a self-injective short local algebra with $e = e(A) = 2.$ In Sections A.8 and A.9, we are going
to survey some properties of the regular modules. We start with the indecomposable
modules of length $2$; they always are regular. 
	\medskip
{\bf Lemma A.9.} {\it Let $A$ be a self-injective short local algebra with $e = 2.$
If $M$ is an indecomposable
module of length $2$, also $\Omega M$ and $\Sigma M (= \mho M)$ are indecomposable modules of length 2.}
	\medskip
Proof. An indecomposable module $M$ of length 2 is local, thus its projective
cover is a free module of rank 1 and therefore $\Omega M$ has length 2, again.
Since $\Omega M$ is a submodule of $PM = {}_AA$ and ${}_AA$ has simple socle, also
$\Omega M$ has simple socle, thus $\Omega M$ is indecomposable. This shows that
$\Omega M$ is an indecomposable module of length 2. Similarly, one shows
that $\Sigma M$ is an indecomposable module of length $2$. $\s$
	\medskip
{\bf Corollary A.10.}  {\it 
Let $A$ be a self-injective short local algebra with $e = 2$.
If $M$ is an indecomposable
module of length $2$, then all modules in the $\mho$-component containing $M$
are indecomposable modules of length $2$.}
	\medskip
If $M$ is indecomposable and of length 2, then $\Ext^1_{L(2)}(M,M) \neq 0,$ 
according to Proposition A.1 in Appendix A. Therefore also $\Ext^1_A(M,M) \neq 0.$
The (uniquely defined) exact
sequence $0 \to M \to M' \to M \to 0$ with $J^2(M') = 0$ will be called
the {\it Kronecker extension} of $M$ (and we write $M' = M[2]$).
Also $\Ext^1(M,\Omega M) \neq 0$, since there is the exact sequence
$0 \to \Omega M \to PM \to M \to 0,$ and this is an $\mho$-sequence. 
These two kinds of extensions, the Kronecker
extension and the $\mho$-extension, are the basic data for dealing with indecomposable
modules of length 2. 
	\medskip
{\bf Proposition A.11.} {\it Let $A$ be a self-injective short local algebra with $e = 2$.
Let  $M, M'$ be indecomposable modules of length $2$. Then $\Ext^1(M,M') = 0$ iff
$M'\not\simeq M$ and $M'\not\simeq \Omega M.$}
	\medskip
Proof. We assume that $M'$ is not isomorphic to $M$ or $\Omega M$, and 
have to show that $\Ext^1(M,M') = 0.$ Let $\epsilon\colon  0 \to M' \to Y \to M \to 0$ be a
non-split exact sequence. If $J^2Y = 0$, then this is an exact sequence of
$L(e)$-modules, thus the Kronecker extension of $M$. Assume now that 
$J^2Y \neq 0.$ Then $Y$ has an indecomposable direct summand isomorphic to ${}_AA$.
Since both ${}_AA$ and $Y$ have length 4, we see that $Y = {}_AA$. Thus $Y \to M$
is a projective cover of $M$ and
and $\epsilon$ is the $\mho$-extension of $M$. In particular, $M' = \Omega M.$
$\s$
	\medskip
{\bf Corollary A.12.} {\it Let $A$ be a self-injective short local algebra with $e = 2$.
Let $M$ be an indecomposable module of length $2$ and $t \ge 0.$
Then $\Ext^{t+1}(M,M) = 0$ iff $\Sigma^tM$ is not isomorphic to $M$ or $\Omega M$.}
	\medskip
Proof. We have $\Ext^{t+1}(M,M) = \Ext^1(M,\Sigma^t M).$ $\s$ 
	\medskip
{\bf Corollary A.13.} {\it Let $A$ be a self-injective short local algebra with $e = 2$.
Let $M$ be an indecomposable module of length $2$ and $t \ge 0.$
Then $\Ext^{t}(M,M) = 0$ for all $t\ge 2$ 
iff the $\mho$-component containing $M$ is of type $\Bbb Z$.} $\s$
	\medskip
{\bf Corollary A.14.} {\it Let $A$ be a self-injective short local algebra with $e = 2$.
Let $M$ be an indecomposable module of length $2$.
The following conditions are equivalent.
\item{\rm(i)} $\Ext^2(M,M) = 0.$
\item{\rm(ii)} There is $i\ge 1$ with $\Ext^i(M,M) = 0.$
\item{\rm(iii)} The $\mho$-component containing $M$ has cardinality at least $3$.
\item{\rm(iv)} $\Omega^2M \not\simeq M.$\par}
	\medskip
Proof. (i) implies (ii) is trivial. 

(ii) implies (iii): Assume that the $\mho$-component
containing $M$ has cardinality at most $2$. Then the modules belonging to the
component are $M$ and $\Omega M$. Thus, for any $i\ge 1$, the module $\Sigma^{i-1}M$
is isomorphic to $M$ or to $\Omega M$. Thus, for any $i \ge 1$, the group $\Ext^i(M,M)=
\Ext^{1}(M,\Sigma^{i-1}M)$ is equal to $\Ext^1(M,M)$ or to $\Ext^1(M,\Omega M)$ and
both groups are non-zero. This contradicts (ii).

(iii) implies (iv) is trivial.

(iy) implies (i). 
We assume that $\Omega^2M \not\simeq M.$ Then clearly also $\Omega M \not\simeq M$.
Now $\Omega^2M \not\simeq M$ implies that $\Omega M \not \simeq \Sigma M$
and $\Omega M \not\simeq M$ implies that $M \not \simeq \Sigma M$. According to
Corollary 1 (with $t = 1$), we have $\Ext^2(M,M) = 0.$
$\s$
	\medskip
If $M$ is an indecomposable non-projective module with $\Omega^2M \simeq M$,
then the $\mho$-component  containing $M$ has cardinality at most 2: 
it consists of the two modules $M,\ \Omega M$ which may or may not be isomorphic.
Here are the two cases: on left the case that $\Omega M \simeq M$;
on the right, the case that $\Omega M \not\simeq M \simeq \Omega^2 M.$

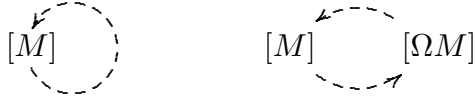
$$
{\beginpicture
    \setcoordinatesystem units <1.5cm,1.5cm>
\put{\beginpicture
\multiput{} at 0 0.4  0 -.4 /  
\put{$[M]$} at 0 0
\setdashes <1mm>
\circulararc 300 degrees from 0 -.15 center at  0.3  0
\setsolid
\arr{0.04 0.2}{0.01 0.15}
\endpicture} at 0 0 
\put{\beginpicture
 \put{$[M]$} at 0 0
\put{$[\Omega M]$} at 1.1 0
\setdashes <1mm>
\setquadratic
\plot 0.2 0.2  0.5 0.3  0.8 0.2 /
\plot 0.2 -.2  0.5 -.3  0.8 -.2 /
\setsolid
\arr{0.22 0.215}{0.2 0.2}
\arr{0.8 -.215}{0.82 -.2}
\endpicture} at 2.5 0
\endpicture}
$$
In the left case $\Omega M \simeq M$, the vector space $\Ext^1(M,M)$ is
2-dimensional, a basis of $\Ext^1(M,M)$ is given by the Kronecker extension and the
$\mho$-extension. In the right case $\Omega M \not\simeq M \simeq \Omega^2 M$, 
the vector space $\Ext^1(M,M)$ is 1-dimensional with basis the Kronecker extension.
Also 
the vector space $\Ext^1(M,\Omega M)$ is 1-dimensional, it has the $\mho$-extension
as basis. 
	\medskip
Of special interest are the self-injective short local algebras with $e = 2$ which
have $\mho$-components of simple regular modules of cardinality at least $3$.
Such $\mho$-components don't exist for commutative algebras, as
Huneke-\c Sega-Vraciu [HSV] have shown:
{\it If $A$ is a commutative self-injective short local algebra, then 
any indecomposable non-projective module $M$ satisfies
$\Ext^i(M,M) \neq 0$ for all $i\ge 1.$} (Note that we have seen in Lemma 10.4 that for 
any commutative, self-injective, short local algebra $A$, any $\mho$-component 
consisting of local modules has cardinality at most $2$.)
	\medskip
Let us look at the operation of $\Omega_A$ on the set of indecomposable modules of length 2.
Let $V$ be a vector space. Non-zero elements $x,y$ of $V$
are called {\it equivalent} provided
there is $\lambda\in k$ (necessarily non-zero) with $y = \lambda x.$ 
We denote by $\Bbb P(V)$ the set of equivalence classes of non-zero elements of $V.$
If $\dim V = n+1,$ then $\Bbb P(V)$ is called the {\it $n$-dimensional projective space}.
We write $\Bbb P^1(k)$ instead of $\Bbb P(k^2)$ and call this the {\it projective line}
over $k$.

If $A$ is a self-injective short local algebra with $e(A) = 2,$ we may identify
$\Bbb P(J/J^2) = \Bbb P^1(k)$ with the set of indecomposable modules of length 2. 
Namely, let us fix  a generating set  $x_0,x_1$ of
${}_AJ$, as well as a non-zero element $z\in J^2.$ 
If $(\alpha_0,\alpha_1)\in k^2$ is a non-zero pair, then 
$M(\alpha_0,\alpha_1) = A(\alpha_0x_0+\alpha_1x_1)$ 
is an indecomposable module of length $2$, all are obtained in this way, and 
$M(\alpha_0,\alpha_1) \simeq M(\alpha'_0,\alpha'_1)$ iff 
$(\alpha_0,\alpha_1)$ and $(\alpha'_0,\alpha'_1)$ are equivalent.

Note that $\Omega = \Omega_A$ provides a permutation of the set of 
(isomorphism classes of) indecomposable modules of length 2, thus of $\Bbb P(J/J^2)
= \Bbb P^1(k)$. We are going to describe
this permutation. The multiplication in $A$ yields a (non-degenerate) bilinear form 
$\mu\colon J/J^2\otimes J/J^2 \to J^2$, say given by the $(2\times 2)$-matrix $B$:
$$
   (\alpha'_0x_0+\alpha'_1x_1)\cdot (\alpha_0x_0+\alpha_1x_1) = 
    \mu((\alpha'_0,\alpha'_1),(\alpha_0,\alpha_1))\cdot z = 
   (\alpha'_0,\alpha'_1)B(\alpha_0,\alpha_1)^t\cdot z
$$
We define a linear map $\eta_B\colon k^2 \to k^2$ by
$$
 \eta_B(\alpha_0,\alpha_1) = (\alpha_0,\alpha_1)\left(\smallmatrix 0&1\cr
                     -1&0\endsmallmatrix\right) B^{-1},
$$
of course, $\eta_B$ is invertible. 
	\medskip
{\bf Lemma A.15.} {\it For any pair non-zero $(\alpha_0,\alpha_1)$ in $k^2$, we have}
$$
  \Omega_A M(\alpha_0,\alpha_1) = M(\eta_B(\alpha_0,\alpha_1)).
$$

Proof. Let $x = \alpha_0x_0+\alpha_1x_1$ and
$y = \alpha'_0x_0+\alpha'_1x_1$, where $(\alpha'_0,\alpha'_1)
= \eta_B(\alpha_0,\alpha_1).$ 
We have
$$
\align
 yx &= \mu((\alpha'_0,\alpha'_1),(\alpha_0,\alpha_1)) 
   =\mu(\eta_B(\alpha_0,\alpha_1),(\alpha_0,\alpha_1)) \cr
   &= 
 (\alpha_0,\alpha_1)\left(\smallmatrix 0&1\cr
                     -1&0\endsmallmatrix\right) B^{-1} B\, (\alpha_0,\alpha_1)^t = 0.
\endalign
$$
But if $x,y$ are elements in $J\setminus J^2$ with $yx = 0$, then there is
the following exact sequence
$$
 {}_AA @>\rho(y)>> {}_AA @>\rho(x)>> {}_AA,
$$
where $\rho(y)$ and $\rho(x)$ denote the right multiplication by $y$, and by $x$,
respectively. The image of $\rho(y)$ is $Ay$, the image of $\rho(x)$ is $Ax$, thus 
$\Omega Ax = Ay.$ $\s$
	\medskip
It is clear that any invertible linear transformation
$\eta\colon k^2 \to k^2$ induces a permutation of $\Bbb P^1(k)$, sending the equivalence class of $(\alpha_0,\alpha_1)$ to the equivalence class of $\eta(\alpha_0,\alpha_1)$.
Such a permutation is called a {\it M\"obius transformation.} 
Thus, we have shown the first part of the following assertion:
	\medskip
{\bf Proposition A.16.} 
{\it The operation of $\Omega_A$ on the set $\Bbb P^1(k)$ of indecomposable modules 
of length $2$ is given by a M\"obius transformation, and
 any M\"obius transformation of $\Bbb P^1(k)$ occurs in this way.}
	\medskip
Proof. It remains to show that 
any M\"obius transformation is of the form $\eta_B$ for some invertible matrix 
$B$. Let $\mu$ be the non-degenerate bilinear form on $k^2$ given by the matrix $B$.
We just need a self-injective short local algebra such that the
multiplication map $J/J^2\otimes J/J^2 \to J^2$ is given by $\mu$. 
Then $\Omega_A$ will operate on the set of indecomposable modules 
of length $2$ by $\eta_B$. 

But for any non-degenerate bilinear form $\mu$ on an $e$-dimensional vector space $W$,
we can define the algebra
$A_\mu$ with underlying vector space $k\oplus W\oplus k$ and with multiplication 
$$
   (c,w,d)(c',w',d') = (cc',cw'+c'w,cd'+c'd+\mu(w,w')),
$$
where $c,c',d,d'\in k$ and $w,w'\in W$. Of course, $A_\mu$ is a self-injective
short local algebra with $e(A_\mu) = e$. And by definition, the
multiplication map $J/J^2\otimes J/J^2 \to J^2$ is just $\mu\colon W\otimes W \to k$.
$\s$
	\bigskip
{\bf A.9. Self-injective algebras with $e = 2$: The regular modules}
	\medskip
As we have mentioned, the aim of Sections A.8 and A.9 is to survey properties of
the regular modules for a self-injective short local algebra $A$ with $e = e(A) = 2.$
Until now, we were dealing just with the indecomposable
modules of length $2$ (they always are regular). 

In our case $e(A) = 2$, an indecomposable $A$-module $M$ is regular provided 
provided it has Loewy length at most 2 and $\bdim M = (m,m)$ for some $m$.
A regular module $M$ is said to be
{\it simple regular} provided the only proper submodule of $M$ which is regular,
is the zero module. Of course, the indecomposable modules of length 2 are simple regular.
We should stress that in case
the base field $k$ is algebraically closed, the indecomposable modules of length 2
are the only simple regular modules! However, if $k$ is not algebraically closed, then
there are additional simple regular modules.
	\medskip
We denote by $\Cal R$ the full subcategory of all modules which are direct sums
of regular indecomposable modules.
Also, we denote by $\widetilde{\Cal R}$ the 
the full subcategory of $\mod K(2)$ consisting of all $K(2)$-modules which are
direct sums of indecomposable regular $K(2)$-modules.
The push-down functor $\pi\colon \mod K(2) \to \mod L(2)$
provides a bijection between the indecomposable objects in $\widetilde{\Cal R}$
and the indecomposable objects in $\Cal R.$ Actually, there is an ideal $\Cal I$
of the category $\mod L(2)$, namely the class of all maps with semisimple image,
such that $\widetilde{\Cal R}/(\Cal I \cap \widetilde{\Cal R}) = \Cal R.$ 
To phrase this differently: if $X,Y$ are regular $L(2)$-modules, then 
$\Hom_{L(2)}(X,Y)$ is just the direct sum $\Hom_{K(2)}(\widetilde X,\widetilde Y)
\oplus \Cal I(X,Y)$, and $\Cal I(X,Y)$ corresponds to the set of 
linear maps $\top X \to \soc Y.$
	
It is important to know that the subcategory $\widetilde{\Cal R}$ of $\mod K(2)$
is an abelian subcategory and the embedding functor is exact. In particular,
$\widetilde{\Cal R}$ is a hereditary length category. The simple objects of 
the abelian
category $\widetilde{\Cal R}$ are called the {\it simple regular} $K(e)$-modules.
It is clear that an $L(2)$-module $R$ is simple regular 
iff $R = \pi\widetilde R$ for some simple regular $K(2)$-module. 

Any indecomposable regular $K(2)$-module $M$ has a unique 
Jordan-H\"older sequence in $\widetilde{\Cal R}$, and all the factors
are isomorphic. We write $\widetilde R[t]$ 
for the indecomposable regular $K(2)$-module with
a filtration with 
$t$ factors of the form $\widetilde R$.
In this way, we obtain a bijection between the
isomorphism classes of the indecomposable objects of $\widetilde{\Cal R}$ and
the pairs $\widetilde R,t$, where $\widetilde R$ is simple regular and $t\in \Bbb N_1.$ 
Using the pushdown functor $\pi$, we obtain
a bijection between the
isomorphism classes of the indecomposable objects of $\Cal R$ and
the pairs $R,t$, where $R$ is simple regular and $t\in \Bbb N_1$
(the number $t$ is called the {\it regular length} of $R[t]$).  
	\medskip
{\bf Proposition A.17.} 
{\it Let $A$ be a self-injective short local algebra with $e = 2.$ Let
$R$ be a simple regular $A$-module and $t\ge 1$.}
\item{(a)} {\it $\Omega R$ is simple regular.}
\item{(b)} {\it $\Omega(R[t]) = (\Omega R)[t]$.}
\item{(c)} {\it The type of the $\mho$-component of $R[t]$ is the same as the
 type of the $\mho$-component of $R$.}
	\medskip
Proof of (a). Since $R$ is regular, its dimension vector is of the form $(m,m)$, thus
also $\bdim \Omega R = (m,m)$ and therefore $\Omega R$ is regular. Let $U$ be a
proper regular submodule of $\Omega R$. Let $V = R/U$, this is also a regular module.
If $\bdim U = (u,u)$ and $\bdim V = (v,v),$ then we have $u+v = m.$ Starting with
injective envelopes of $U$ and $V$, 
the horseshoe lemma asserts that there is an injective module $I$ such that
$\Sigma \Omega R \oplus I$ has a submodule of the form $\Sigma U$ with factor module
$\Sigma V.$ Since $\bdim U + \bdim V = \bdim R,$ we see that $I = 0,$ thus
$\Sigma U$ is a submodule of $R$ with factor module $\Sigma V$. By assumption, $V \neq 0$,
thus $\Sigma V \neq 0$. It follows that $\Sigma U$ is a proper regular submodule of $R$.
Since $R$ is simple regular, we have $\Sigma U = 0,$ thus $U = 0.$ This shows that
$R$ is simple regular.
	\smallskip
Proof of (b). We use induction on $t$. Let $t\ge 2.$ There is an exact sequence
$0 \to R[t-1] \to R[t] \to R \to 0.$ The horseshoe lemma asserts that 
there is a projective
module $P$ such that $\Omega(R[t])\oplus P$ is an extension of 
$\Omega(R[t-1])$ 
by $\Omega R$, and by induction we have $\Omega(R[t-1]) = (\Omega R)[t-1]$. 
Now $\Omega(R[t])$ is indecomposable, thus $\Omega(R[t]) = 
(\Omega R)[s]$ for some $1\le s \le t.$ Assume that $s < t.$
We apply $\Sigma$. On the one hand, we have $\Sigma\Omega(R[t])= R[t]$,
on the other hand, we have $\Sigma((\Omega R)[s]) = \Sigma\Omega(R[s]) = R[s]$,
where we use again induction. But for $s < t$, the module $R[t]$ is not isomorphic to 
$R[s]$. It follows that $s = t,$ thus $\Omega(R[t]) = (\Omega R)[t]$.
	\smallskip 
Proof of (c). The type of the $\mho$-component of an indecomposable module $M$ 
is $\widetilde {\Bbb A}_n$ iff $n \ge 0$ is minimal
with $\Omega^n M \simeq M$ (and $\Bbb Z$, if there is no $n$ of this kind). 
According to (b) we have $\Omega^n (R[t]) \simeq R[t]$ iff 
$\Omega^n R \simeq R.$
$\s$ 
	\bigskip
{\bf Historical Remark for A.8 and A.9.}
The representations of $K(2)$ have been classified by
Kronecker in 1890, 
completing earlier partial results by Jordan and Weierstrass,
as mentioned for example in  [ARS].
This classification plays an
important role in many parts of mathematics. A standard reference for the
matrix approach (in the language of matrix pencils)
is Gantmacher's book on matrix theory [Gm].
There is the
equivalent theory of coherent sheaves over the projective line, where the usual
reference is the splitting theorem of Grothendieck (but one should be aware that
this result can be traced
back to  Hilbert (1905), Plemelj (1908), and G. D. Birkhoff (1913), see [OSS]).
	
The category $\mod L(2)$ plays also a prominent role in modular
representation theory, since the group algebra $kG$ of the Klein four group $G$ in
characteristic 2 is a self-injective short local algebra with $e = 2$ (namely 
$kG = A_1$). The usual
reference are papers by Bashev (1961) 
and Heller-Reiner (1961), see Benson [B].
	\bigskip
{\bf A.10. Self-injective algebras with $e = 2$: Normal forms} 
	\medskip
Let us now assume that $k$ is algebraically
closed. In this case, it is easy to determine normal forms for the self-injective 
short local algebras $A$ with $e(A) = 2$.
	\smallskip
\item{$\bullet$} There are the algebras $A_q = k\langle x,y\rangle/\langle x^2, y^2, xy+qyx\rangle$ with $q\in k^*;$ note that $A_q$ is isomorphic to $A_{q^{-1}}.$
\item{$\bullet$} In addition, there is the algebra 
$A_0 = k\langle x,y\rangle/\langle x^2, y^2 +xy, 
  y^2-yx\rangle.$ 
	\medskip
Sketch of proof. Let $A$ be a self-injective short local $k$-algebra with $e(A) = 2.$ First
one shows that there is always an element $x\in J\setminus J^2$ with $x^2 = 0.$ Then one
takes  $z\in J$ such that $x,z$ is a generating set for ${}_AJ.$ 
Since $A$ is self-injective, the elements $xz$ and
$zx$ have to be non-zero, thus there is $q\in k^*$ with $xz+qzx = 0.$ If we have 
$z^2 = 0,$ then $A \simeq A_q.$ Thus, let $z^2 \neq 0.$ Then there is 
$\alpha\in k^*$ with $z^2 = \alpha zx.$
If $q\neq 1$, then $y = z+(q-1)^{-1}\alpha x$
satisfies $y^2 = 0$ and we deal again with the previous case (with $z$ replaced by $y$).
Thus, the case $z^2 \neq 0$ and $q = 1$ remains: here, 
the elements $x$ and $y = \alpha^{-1}z$ satisfy the defining relations for $A_0$.
$\s$
	\medskip 
If $A = A_q$, with $q\in k^*$, then 
   the modules $Ax$ and $Ay$ are $\Omega$-periodic of 
   period 1. 
The $\Omega$-orbit of a module $A(x+\alpha y)$ with $\alpha\neq 0$ 
   is the set $\{A(x+q^t\alpha y) \mid t\in \Bbb Z\}$; its cardinality is equal to the
   multiplicative order $o(q)$ of $q$. 
The algebra $A_1$ is the exterior algebra in 2 generators; in this case
all regular indecomposable modules are $\Omega$-periodic of period 1, thus the
corresponding M\"obius transformation is the identity.
The remaining algebras $A_q$ (with $q\in k^*$ 
   and $q  \neq 1$) have precisely two $\Omega$-orbits of cardinality 1 which consist
   of indecomposable modules of length 2, namely the orbits of $Ax$ and $Ay.$ This means
that the corresponding M\"obius transformation has precisely two fixed points.
All other orbits of indecomposable modules of length 2
have cardinality $o(q)$. (If $o(q) = \infty$, then $A_q$ is just
   a quantum exterior algebra in 2 generators as discussed in A.11.)

For $A = A_0$, the module $Ax$ is $\Omega$-periodic of period 1.
   The $\Omega$-orbit of the module $A(y+\alpha x)$ with $\alpha\in k$ 
   is the set $\{A(y+(\alpha+t)x) \mid t\in \Bbb Z\}$; thus its cardinality is equal to the
   characteristic $\ch k$ of $k$. 
The corresponding M\"obius transformation has precisely one fixed points (these M\"obius
transformations are often called {\it parabolic}).
	\medskip
If $A$ is commutative, then $A = A_{-1}$ or else the characteristic of $k$
is $2$ and $A = A_0.$ In both cases, the corresponding M\"obius transformation has order 2.
	
Of course, all the algebras $A_q$ with $q\in k^*$ are special biserial.
If the characteristic of $k$ is 2, also $A_0$ is special biserial, namely
isomorphic to $A\langle x,y\rangle/\langle xy, yx, x^2-y^2\rangle.$ 
If the characteristic of $k$ is different from $2$, then $A_0$ is (biserial, but)
not special biserial. +
	\bigskip
{\bf A.11. Example: The quantum exterior algebra in 2 generators}
	\medskip 
The quantum exterior algebra $A$ in two variables $x,y$
is the $k$-algebra generated by $x,y$ with the relations
$x^2,\ y^2,\ xy+qyx,$ where $q\in k^* = k\setminus\{0\}$ 
has infinite (multiplicative) order.
Note that the
elements $1,\ x,\ y,$ and $yx$ form a basis for $A$.
	\medskip
We consider the left ideals $M_\alpha = A(x+\alpha y)$ with $\alpha\in k;$ these
are indecomposable modules of length 2. 
	\medskip
{\bf Lemma A.18.} {\it Let $A$ be the quantum exterior algebra $A$ in two variables.
If $\alpha\in k,$, then}
$$
 \Omega M_\alpha = M_{q\alpha}.
$$
	\medskip
Proof. If $a\in A$, let
$\rho_a\colon {}_AA \to {}_AA$ be the right multiplication by $a$. Of course, the image 
of $\rho_a$ is the left ideal $Aa$. 
The relations show that $(x+q\alpha y)(x+\alpha y) = 0.$ This implies that the
composition 
$$
 {}_AA @>\rho_{x+q\alpha y}>> {}_AA @>\rho_{x+\alpha y}>>  {}_AA 
$$
is zero. The image of the left map is $M_{q\alpha},$ the image of the right map
is $M_\alpha.$ It follows that $\Omega M_\alpha = M_{q\alpha}.$
$\s$
	\medskip
Let $A$ be the quantum exterior algebra in two variables $x,y$. 
The $\mho$-component containing $M_1$ looks as follows:
$$
{\beginpicture
    \setcoordinatesystem units <1.5cm,1cm>
\put{$M_1$\strut} at 0 0
\put{$M_q$\strut} at -1 0
\put{$M_{q^2}$\strut} at -2 0
\put{$M_{q^{-1}}$\strut} at 1 0
\put{$M_{q^{-2}}$\strut} at 2.05 0
\put{$\cdots$} at -3 0
\put{$\cdots$} at 3 0
\setdashes <.7mm>
\plot -0.2 0  -0.7 0 /
\arr{-1.3 0}{-1.7 0}
\arr{-2.3 0}{-2.7 0}
\arr{0.65 0}{0.2 0}
\arr{1.7 0}{1.35 0}
\arr{2.75 0}{2.45 0}
\setsolid
\arr{-0.65 0}{-0.7 0}
\arr{-1.65 0}{-1.7 0}
\arr{-2.65 0}{-2.7 0}
\arr{0.25 0}{0.2 0}
\arr{1.4 0}{1.35 0}
\arr{2.5 0}{2.45 0}

\endpicture}
$$
Since the modules $M_{q^i}$ with $i\in \Bbb Z$ are pairwise non-isomorphic, we see
that this $\mho$-component is of type $\Bbb Z.$ Thus we see:
	\medskip 
{\bf Proposition A.19.} 
{\it If $A$ is the quantum exterior algebra in $2$ generators, then
there exists a two-dimensional indecomposable module $M$} (namely $M = M_1$) {\it 
with $\mho$-component of type $\Bbb Z$. Thus}
$$
   \Ext^i(M,M) = 0\quad\text{\it for all $i\ge 2$, whereas}\quad \Ext^1(M,M) \neq 0.
$$
$\s$
	
{\bf Corollary A.20 (Smal\o{}\ [Sm]).} 
{\it If $A$ is the quantum exterior algebra in $2$ generators, 
there are indecomposable modules $M$ and $N_i$ with $i\in \Bbb N_1$ such
that $\Ext^i(M,N_i) \neq 0$ and  $\Ext^j(M,N_i) = 0$ for all $j > i.$}
	\medskip
Proof. Let $M$ be a 2-dimensional 
indecomposable module with $\mho$-component of type $\Bbb Z.$ Let
$N_i = \Omega^{i-1}M.$ Then 
$$
 \Ext^i(M,N_i) = \Ext^i(M,\Omega^{i-1}M) = \Ext^1(M,\Sigma^{i-1}\Omega^{i-1}M) =
 \Ext^1(M,M) \neq 0.
$$
Also, for $j > i$, we have
$$
 \Ext^j(M,N_i) = \Ext^{j-i+1}(M,\Sigma^{i-1}\Omega^{i-1}M) = \Ext^{j-i+1}(M,M) = 0,
$$
since $j-i+1 \ge 2$. 
$\s$
	\medskip
Recall that Auslander had conjectured that for every module $M$ there exists a
bound $b(M)$ with the following property: if $N$ is a module with
$\Ext^j(M,N) = 0$ for $j \gg 0$, then 
$\Ext^j(M,N) = 0$ for $j > b(M).$ Corollary A.20 shows:
	\medskip
{\bf Corollary A.21 (Smal\o \ [Sm]).} {\it The quantum exterior algebra 
in two variables is a
counter-example to the Auslander conjecture.} $\s$
	\medskip
The first counter-example for the Auslander conjecture 
was given by Jorgensen-\c Sega [JS1]. 
	\bigskip
{\bf A.12. Koszul modules}
	\smallskip
The paper [RZ3] will draw the attention to Koszul modules as
defined by Herzog and Iyengar [HI], see also [AIS].
If $A$ is a short local algebra, then
{\it an $A$-module $M$ of Loewy length at most $2$ is
a Koszul module if and only if all the modules $\Omega^nM$ with $n\ge 0$ are
aligned}, see [RZ3].

Since for a self-injective algebra $A$, any
$A$-module is Gorenstein-projective, the minimal projective resolutions of all
indecomposable non-projective modules are displayed by the $\mho$-quiver.
It follows:
	\medskip
{\bf Proposition A.22 ([Sj, MV2, AIS]).}
{\it Let $A$ be a self-injective short local algebra with $e\ge 2$.
If $M$ is indecomposable, then $M$ is Koszul
if and only if $M$ is not preprojective in the sense of Auslander-Smal\o{}
(thus not of the form $\pi P_1, \pi P_2,\dots$).} $\s$
	\medskip
Let us add:
	\medskip
{\bf Proposition A.23.}
{\it Let $A$ be a self-injective short local algebra. If $e\ge 2$, then
the simple module $S$ is a Koszul module, and
for any module $M$, there exists $m\ge 0$ such that $\Omega^mM$ is Koszul.
If $e = 1$, and $a = 1$, then the only Koszul modules are the projective modules.}
	\medskip
Proof. We can assume that $M$ is an indecomposable module.
First, let $e\ge 2$ and assume that $M$ is not
Koszul, then $M = \pi P_m$ for some $m\ge 1$ and therefore
$\Omega^m (\pi P_m) = S$ is Koszul. If $e = 1$, and $a = 1$,
then $A$ is uniserial, thus
$M$ is isomorphic to $k$, ${}_AJ$ or ${}_AA$, and, of course, the modules $k$
and ${}_AJ$ are not Koszul.
$\s$
	\medskip
{\bf Historical Remark.} The Koszul modules over a self-injective
short local algebra have been determined by Sj\"odin [Sj],
Mart\'\i nez-Villa [MV2] and Avramov-Iyengar-\c Sega [AIS].
We hope that our outline of the general setting explains what is considered
as a surprising behavior in [AIS].

Already in 1979, Sj\"odin [Sj] has looked for
indecomposable non-projective modules $M$ at the
power series $P^A_M = \sum_{n\ge 0} \beta_n(M)T^n$ (called the {\it Poincar\'e series}
of $M$). He showed that for a self-injective short local algebra $A$, the series
$P^A_M$ is rational (this follows from the fact that
$\Omega^m M$ is Koszul for some $m \ge 0$).
	\bigskip\medskip
\vfill\eject
{\bf Appendix B. A formula of Avramov-Iyengar-\c Sega}
	\medskip
{\bf B.1.  The sequences $b(e,a)_n$}
	\smallskip
Let $e, a $ be real numbers.
We define recursively the sequence $b_n = b(e,a)_n$ with $n\ge -1$ as follows:
$b_{-1} = 0,$ $b_0 = 1$ and
$$
  b_{n+1} = eb_{n}-ab_{n-1}, \tag{$*$}
$$
for $n\ge 0.$ In this paper, we are interested it is the case that $e, a$ 
are natural numbers and $a\le e^2$. 
Namely, if $A$ is a short local algebra with Hilbert type $(e,a)$, then $e,a$ are
natural numbers with $a\le e^2$, and the recursion rule $(*)$
has popped up in Section 5.6, when dealing with a module $M$ such that
both $M$ and $\Omega M$ are aligned.

As a consequence, we see the relevance of the numbers $b_n = b(e,a)_n$:
{\it We have $\beta_n(S) = b_n$ for all $0 \le n \le N$
if and only if the modules $\Omega^n S$ with $0\le n < N$ are aligned.} As we have mentioned
in A.12 (with reference to [RZ3]),
the module $S$ is a Koszul module in the sense of [HI] iff
all the modules $\Omega^nM$ with $n\ge 0$ are aligned.
Thus {\it $S$ is a Koszul module iff $\beta_n(S) = b_n$ for all $n\ge 0$}
(and then $\bdim \Omega^n S = (b_n,b_{n-1})$).
	\medskip
The paper [AIS] aimed to provide a concise formula for the numbers $b(e,1)_n$ with $e\ge 3$,
but the formula presented there was slightly distorted and
usually did not even give integers. We are indebted to
Avramov, Iyengar and \c Sega for communicating to us a proper revision
and to allow us to include it here.
	\bigskip
{\bf B.2. The formula of Avramov, Iyengar, \c Sega}
	\smallskip
{\bf Theorem B.1 (Avramov, Iyengar, \c Sega).} {\it If $a < \frac14e^2$, then
for all $n\ge 0$}
$$
 b(e,a)_n =  \frac 1{2^n} \sum_{j=0}^{\lfloor \frac {n}2\rfloor}
 \binom{n+1}{2j+1} (e^2-4a)^{j}e^{n-2j}.
$$

Proof (Avramov, Iyengar, \c Sega):
Since we assume
that $a < \frac14e^2$, the roots of the polynomial $1-eT+aT^2$ are real numbers, and do not
coincide. The roots are
$$
 \lambda = \frac{e - q}2,\quad \text{and} \quad \rho = \frac{e+ q}2,
  \quad \text{where}\quad q = \sqrt{e^2-4a}\ > \ 0.
$$
Starting with the factorization
$$
 1-eT+aT^2 = (1-\rho T)(1-\lambda T),
$$
we may look at the power series expansion
of the rational function $(1-eT+aT^2)^{-1}$:
$$
\frac 1{1-eT+aT^2} \
  =\  \frac 1{(\rho -\lambda)}\left(\frac \rho {1-\rho T}
    - \frac{\lambda}{1-\lambda T}\right) \
  = \ \frac{1}{q}\; \sum_{n\ge 0} (\rho^{n+1}-\lambda^{n+1})\,T^n
$$
Of course, we have
$$
\frac 1{1-eT+aT^2} = \sum_{n\ge 0} b(e,a)_n T^n,
$$
therefore
$$
  b(e,a)_n =  \tfrac1q (\rho^{n+1}-\lambda^{n+1}).
$$
The binomial expansions of $\rho^{n+1}$ and $\lambda^{n+1}$ yield
$$
\align
\rho ^{n+1} - \lambda^{n+1}
  &= \sum_{i=0}^{n+1}\binom{n+1}i \frac 1{2^{n+1}}
  \left(e^{n+1-i}q^i - (-1)^i e^{n+1-i}
  q^i\right)\\
  &= \frac 1{2^n} \sum_{j=0}^{\lfloor \frac {n}2\rfloor} \binom{n+1}{2j+1}
  q^{2j+1}e^{n-2j}
\endalign
$$
Altogether, one gets that
$$
\align
b(e,a)_n =  \tfrac1q (\rho^{n+1}-\lambda^{n+1}) &=
\frac 1{2^n}\sum_{j=0}^{\lfloor \frac {n}2\rfloor} \binom{n+1}{2j+1} q^{2j}e^{n-2j},\cr
& =
\frac 1{2^n}\sum_{j=0}^{\lfloor \frac {n}2\rfloor} \binom{n+1}{2j+1} (e^2-4a)^{j}e^{n-2j}.
\endalign
$$
\vglue-.5cm
$\s$

	\medskip
Note that the formula exhibited above is already of interest in the case $e = 3$ and
$a = 1$.
In this case the numbers $b_n = b(3,1)_n$ are just the even-index Fibonacci numbers
(see Section A.1 in Appendix A).
	\bigskip\medskip
{\bf Acknowledgement.} The authors would like to thank L. L. Avramov, S. B. Iyengar
and L. M. \c Sega for providing the formula exhibited in Appendix B, and to
D. Jorgensen and R. Takahashi
for valuable help concerning the literature on commutative rings.
They are indebted to a referee for several suggestions concerning the
presentation of the results. Supported by NSFC 12131015, 11971304.

	\bigskip\bigskip
{\bf References}
	\medskip
\item{[AGP]} L. L. Avramov, V. N. Gasharov, I. V. Peeva. Complete intersection dimension.
  Inst. Hautes Etudes Sci. Publ. Math. 86 (1997), 67--114.
\item{[AIS]} L. L. Avramov, S. B. Iyengar, L. M. \c Sega. Free resolutions
  over short local rings. J. London Math. Soc. 78 (2008), 459--476.
\item{[AR]} M. Auslander, I. Reiten. On a generalization of the Nakayama
  conjecture. Proc. Amer. Math. Soc. 52 (1975), 69--74.
\item{[ARS]}  M. Auslander, I. Reiten, S. O. Smal\o. Representation Theory
  of Artin Algebras. Cambridge Studies in Advanced Math. 36.
  Cambridge University Press (1995).
\item{[AS]} M. Auslander, S. O. Smal\o. Preprojective modules over artin algebras.
  J. Algebra 66 (1980), 61--122.
\item{[B]} D. J. Benson. Representations and cohomology. I. Cambridge University Press
  (1991).
\item{[BGP]} I. N. Bernstein, I. M. Gelfand, V. A. Ponomarev.
   Coxeter functors and Gabriel's theorem.
   Uspechi Mat. Nauk. 28, 19--33, (1973),  Russian Math. Surveys 28, 17--32 (1973).
\item{[BH]} W. Bruns, J. Herzog. Cohen-Macaulay Rings. Cambridge Studies
   in advanced mathematics 39. Cambridge University Press (1993).
\item{[CV]} L. W. Christensen, O. Veliche. Acyclicity over local rings with
  radical cube zero. Illinois J. Math. 51 (2007), 1439--1454.
\item{[DR]} V. Dlab, C. M. Ringel.
  Indecomposable representations of graphs and algebras.
  Memoirs of the American Mathematical Society. Providence, RI.
  Vol. 173 (1976).
\item{[Gb]} P. Gabriel. Unzerlegbare Darstellungen I. manuscripta math. 6 (1972),
  71--103.
\item{[Gm]} F. R. Gantmacher. Matrizentheorie.
  Springer-Verlag Berlin (2013).
\item{[GP]} V. N. Gasharov, I. V. Peeva. Boundedness versus periodicity over
  commutative local rings. Trans. Amer. Math. Soc. 320 (1990), 569--580.
\item{[HI]} J. Herzog, S. Iyengar. Koszul modules. J. Pure Appl. Algebra
  201 (2005), 154--188.
\item{[Ho1]} M. Hoshino. Modules without self-extensions and Nakayama's conjecture.
  Arch. Math. 43 (1982), 111--137.
\item{[Ho2]} M. Hoshino. On algebras with radical cube zero,
  Arch Math 52 (1989), 226--232.
\item{[HJS]} M. T. Hughes, D. A. Jorgensen, L. M. \c Sega.
  Acyclic complexes of finitely generated free modules over local ring.
  Math. Scand. 105 (2009), 85--98.
\item{[HSV]}
  C. Huneke, L. M. \c Sega, A. N. Vraciu.  Vanishing of Ext and
  Tor over some Cohen-Macaulay local rings. Illinois J. Math. 48 (2004),
  no. 1, 295--317. 
\item{[I]} O. Iyama. Auslander correspondences. Adv. Math. 210
  (2007), 51--82.
\item{[JS1]} D. A. Jorgensen, L. M. \c Sega. 
  Nonvanishing cohomology and classes of Gorenstein rings. Adv. Math.
  188(2) (2004), 470--490.
\item{[JS2]} D. A. Jorgensen, L. M. \c{S}ega. Independence of the
  total reflexivity conditions for modules.  Algebras and
  Representation Theory 9(2) (2006), 217--226.
\item{[L]} J. Lescot. Asymptotic properties of Betti numbers of modules over
  certain rings. J. Pure Appl. Algebra 38 (1985), 287--298.
\item{[M1]} R. Marczinzik. 
  Simple reflexive modules over Artin algebras. J. Algebra Appl. 18 (2019), no. 10, 1950193
\item{[M2]} R. Marczinzik. On weakly Gorenstein algebras.
  arXiv:1908.04738
\item{[MV1]} R. Mart\'\i nez-Villa. Algebras stably equivalent to l-hereditary. In:
  Representation theory II, Springer LNM 832 (1980), 396--431.
\item{[MV2]} R. Mart\'\i nez-Villa. Applications of Koszul algebras: The preprojective
   algebras. In: Canadian Math. Soc. Conference Proceedings. 18 (1996), 487--504
\item{[OSS]} C. Okonek, M. Schneider, H. Spindler. Vector bundles on complex
  projective spaces. Progress in Mathematics. Birkh\"auser (1980).
\item{[R1]} C. M. Ringel. Representations of $k$-species and bimodules.
   J. Algebra 41 (1976), 269--302.
\item{[R2]} C. M. Ringel.  Simple reflexive modules over finite-dimensional algebras.
    J. Algebra Appl.
    20, Issue 09 (2021), Article 2150166.
    \newline
    DOI: 10.1142/S0219498821501668
\item{[RZ1]}  C. M. Ringel, P. Zhang. Gorenstein-projective and
  semi-Gorenstein-projective modules.
  Algebra \& Number Theory 14-1 (2020), 1--36.
  dx.doi.org/10.2140/ant.2020.14.1
\item{[RZ2]}  C. M. Ringel, P. Zhang. Gorenstein-projective and
  semi-Gorenstein-projective modules. II.
  J. Pure Appl. Algebra 224 (2020) 106248. \newline
  https://doi.org/10.1016/j.jpaa.2019.106248
\item{[RZ3]}  C. M. Ringel, P. Zhang. Koszul modules (and the $\Omega$-growth of modules)
  over short local algebras. J. Pure Appl. Algebra 225 (2021) 106772 \newline
  https://doi.org/10.1016/j.jpaa.2021.106772 
\item{[Sj]} G. Sj\"odin. The Poincar\'e series of modules over a local Gorenstein ring with
  $\bold m^3 = 0.$ Preprint 2. Mathematiska Institutionen. Stockholms Universitet (1979).
\item{[Sm]} S. O. Smal\o. Local limitations of the $\Ext$ functor do not exist,
   Bull. London Math. Soc. 38 (2006), 97--98.
\item{[V]} O. Veliche. Construction of modules with finite homological dimension.
  J. Algebra 250 (2002), 427--449.
\item{[Y]} Y. Yoshino. Modules of G-dimension zero over local
  rings with the cube of maximal ideal being zero. In: Commutative algebra,
  singularities and computer algebra (Sinaia, 2002), NATO Sci. Ser.
  II, Kluwer, Dordrecht (2003), 255--273. 
	\bigskip
{\baselineskip=1pt
\rmk
C. M. Ringel\par
Fakult\"at f\"ur Mathematik, Universit\"at Bielefeld \par
POBox 100131, D-33501 Bielefeld, Germany  \par
ringel\@math.uni-bielefeld.de
\smallskip

P. Zhang \par
School of Mathematical Sciences, Shanghai Jiao Tong University \par
Shanghai 200240, P. R. China.\par
pzhang\@sjtu.edu.cn\par}

	\vfill\eject

\bye